\documentclass[preprint,12pt]{elsarticle}

\usepackage{graphicx}
\usepackage{multicol}
\usepackage{mathtools}
\usepackage{amssymb}
\usepackage{tikz}
\usepackage{multirow}
\usepackage{graphicx}       
\usepackage{multicol}       
\usepackage[bottom]{footmisc}
\usepackage{subcaption}
\usepackage{pgfplots}
\usepackage{newtxtext}       
\usepackage[varvw]{newtxmath}       
\usepackage{todonotes}
\usepackage[overload]{empheq}
\usepackage{ulem}
\usepackage{amsmath}
\usepackage{pgfplots}
\usepackage{tikz}
\usepackage{pgfplots}
\usepackage{wrapfig}

\usetikzlibrary{spy,calc}
\usetikzlibrary{decorations.pathreplacing}
\usepgfplotslibrary{groupplots}
\usepgfplotslibrary{colorbrewer}
\usepgfplotslibrary{fillbetween}
\usetikzlibrary{colorbrewer}
\usetikzlibrary{external}

\newcommand*{\affmark}[1][*]{\textsuperscript{#1}}

\definecolor{grey}{RGB}{230.,230.,230}
\definecolor{lightgrey}{RGB}{240.,240.,250}
\definecolor{lightblue}{RGB}{80, 200,255}
\definecolor{whiteblue}{RGB}{120, 220,255}
\definecolor{darkblue}{RGB}{20, 200, 200}
\definecolor{lighturkis}{RGB}{42, 220,255}
\definecolor{whitegrey}{RGB}{200, 200,220}

\newcommand{\uB}{\boldsymbol{u}}
\newcommand{\FB}{\boldsymbol{F}}
\newcommand{\vB}{\boldsymbol{v}}
\newcommand{\xB}{\boldsymbol{x}}
\newcommand{\nB}{\boldsymbol{n}}

\newcommand{\jB}{\boldsymbol{j}}

\newcommand{\QB}{\boldsymbol{Q}}
\newcommand{\SB}{\boldsymbol{S}}
\newcommand{\sB}{s^{\#}}
\newcommand{\sInt}{\boldsymbol{s}_{\Gamma}}
\newcommand{\qB}{\boldsymbol{q}}
\newcommand{\PB}{\boldsymbol{P}}
\newcommand{\JB}{\boldsymbol{J}}

\newcommand{\smacro}{\boldsymbol{\Gamma}^s} 
\newcommand{\smicro}{\boldsymbol{\Gamma}^m}

\pgfplotsset{
	compat=1.16,
	cycle list/Dark2,
	siunitxlabels/.style={
		/pgfplots/typeset ticklabel/.code={{\pgfmathparse{\tick}$\num[zero-decimal-to-integer]{\pgfmathresult}$}},
	},
	colormap/hp_1/.style={colormap={hp_1}{
			rgb=(0.403921568627451,0.509803921568627,0.705882352941176);
			rgb=(0.376470588235294,0.737254901960784,0.83921568627451);
			rgb=(0.87843137254902,0.905882352941176,0.92156862745098);
			rgb=(0.176470588235294,0.176470588235294,0.176470588235294);
		}
	},
	colormap/XRay/.style={colormap={XRay}{
			rgb=(1.0, 1.0, 1.0);
			rgb=(0.0, 0.0, 0.0);
		}
	},
	colormap/erdciceFireL/.style={colormap={erdciceFireL}{
			rgb=(0.870485, 0.913768, 0.832905);
			rgb=(0.586919, 0.887865, 0.934003);
			rgb=(0.31583, 0.776442, 0.867858);
			rgb=(0.18302, 0.632034, 0.787722);
			rgb=(0.117909, 0.484134, 0.713825);
			rgb=(0.0507239, 0.335979, 0.654741);
			rgb=(0.0, 0.209874, 0.511832);
			rgb=(0.0, 0.114689, 0.28935);
			rgb=(0.0157519, 0.00332021, 4.55569e-08);
			rgb=(0.312914, 0.0, 0.0);
			rgb=(0.520865, 0.0, 0.0);
			rgb=(0.680105, 0.15255, 0.0025996);
			rgb=(0.785109, 0.339479, 0.000797922);
			rgb=(0.857354, 0.522494, 0.0);
			rgb=(0.910974, 0.699774, 0.0);
			rgb=(0.951921, 0.842817, 0.478545);
			rgb=(0.881371, 0.912178, 0.818099);
			}
	},
	colormap/erdciceFireH/.style={colormap={erdciceFireH}{
			rgb=(4.05432e-07, 0.0, 5.90122e-06);
			rgb=(0.0, 0.120401, 0.302675);
			rgb=(0.0, 0.216583, 0.524574);
			rgb=(0.0552475, 0.345025, 0.6595);
			rgb=(0.128047, 0.492588, 0.720288);
			rgb=(0.188955, 0.641309, 0.792092);
			rgb=(0.327673, 0.784935, 0.873434);
			rgb=(0.60824, 0.892164, 0.935547);
			rgb=(0.881371, 0.912178, 0.818099);
			rgb=(0.951407, 0.835621, 0.449279);
			rgb=(0.904481, 0.690489, 0.0);
			rgb=(0.85407, 0.510864, 0.0);
			rgb=(0.777093, 0.33018, 0.00088199);
			rgb=(0.672862, 0.139087, 0.00269398);
			rgb=(0.508815, 0.0, 0.0);
			rgb=(0.299417, 0.000366289, 0.000547829);
			rgb=(0.0157519, 0.00332021, 4.55569e-08);
			}
	},
	colormap/inferno/.style={colormap={inferno}{
			rgb=(0.001462, 0.000466, 0.013866);
			rgb=(0.002267, 0.001270, 0.018570);
			rgb=(0.003299, 0.002249, 0.024239);
			rgb=(0.004547, 0.003392, 0.030909);
			rgb=(0.006006, 0.004692, 0.038558);
			rgb=(0.007676, 0.006136, 0.046836);
			rgb=(0.009561, 0.007713, 0.055143);
			rgb=(0.011663, 0.009417, 0.063460);
			rgb=(0.013995, 0.011225, 0.071862);
			rgb=(0.016561, 0.013136, 0.080282);
			rgb=(0.019373, 0.015133, 0.088767);
			rgb=(0.022447, 0.017199, 0.097327);
			rgb=(0.025793, 0.019331, 0.105930);
			rgb=(0.029432, 0.021503, 0.114621);
			rgb=(0.033385, 0.023702, 0.123397);
			rgb=(0.037668, 0.025921, 0.132232);
			rgb=(0.042253, 0.028139, 0.141141);
			rgb=(0.046915, 0.030324, 0.150164);
			rgb=(0.051644, 0.032474, 0.159254);
			rgb=(0.056449, 0.034569, 0.168414);
			rgb=(0.061340, 0.036590, 0.177642);
			rgb=(0.066331, 0.038504, 0.186962);
			rgb=(0.071429, 0.040294, 0.196354);
			rgb=(0.076637, 0.041905, 0.205799);
			rgb=(0.081962, 0.043328, 0.215289);
			rgb=(0.087411, 0.044556, 0.224813);
			rgb=(0.092990, 0.045583, 0.234358);
			rgb=(0.098702, 0.046402, 0.243904);
			rgb=(0.104551, 0.047008, 0.253430);
			rgb=(0.110536, 0.047399, 0.262912);
			rgb=(0.116656, 0.047574, 0.272321);
			rgb=(0.122908, 0.047536, 0.281624);
			rgb=(0.129285, 0.047293, 0.290788);
			rgb=(0.135778, 0.046856, 0.299776);
			rgb=(0.142378, 0.046242, 0.308553);
			rgb=(0.149073, 0.045468, 0.317085);
			rgb=(0.155850, 0.044559, 0.325338);
			rgb=(0.162689, 0.043554, 0.333277);
			rgb=(0.169575, 0.042489, 0.340874);
			rgb=(0.176493, 0.041402, 0.348111);
			rgb=(0.183429, 0.040329, 0.354971);
			rgb=(0.190367, 0.039309, 0.361447);
			rgb=(0.197297, 0.038400, 0.367535);
			rgb=(0.204209, 0.037632, 0.373238);
			rgb=(0.211095, 0.037030, 0.378563);
			rgb=(0.217949, 0.036615, 0.383522);
			rgb=(0.224763, 0.036405, 0.388129);
			rgb=(0.231538, 0.036405, 0.392400);
			rgb=(0.238273, 0.036621, 0.396353);
			rgb=(0.244967, 0.037055, 0.400007);
			rgb=(0.251620, 0.037705, 0.403378);
			rgb=(0.258234, 0.038571, 0.406485);
			rgb=(0.264810, 0.039647, 0.409345);
			rgb=(0.271347, 0.040922, 0.411976);
			rgb=(0.277850, 0.042353, 0.414392);
			rgb=(0.284321, 0.043933, 0.416608);
			rgb=(0.290763, 0.045644, 0.418637);
			rgb=(0.297178, 0.047470, 0.420491);
			rgb=(0.303568, 0.049396, 0.422182);
			rgb=(0.309935, 0.051407, 0.423721);
			rgb=(0.316282, 0.053490, 0.425116);
			rgb=(0.322610, 0.055634, 0.426377);
			rgb=(0.328921, 0.057827, 0.427511);
			rgb=(0.335217, 0.060060, 0.428524);
			rgb=(0.341500, 0.062325, 0.429425);
			rgb=(0.347771, 0.064616, 0.430217);
			rgb=(0.354032, 0.066925, 0.430906);
			rgb=(0.360284, 0.069247, 0.431497);
			rgb=(0.366529, 0.071579, 0.431994);
			rgb=(0.372768, 0.073915, 0.432400);
			rgb=(0.379001, 0.076253, 0.432719);
			rgb=(0.385228, 0.078591, 0.432955);
			rgb=(0.391453, 0.080927, 0.433109);
			rgb=(0.397674, 0.083257, 0.433183);
			rgb=(0.403894, 0.085580, 0.433179);
			rgb=(0.410113, 0.087896, 0.433098);
			rgb=(0.416331, 0.090203, 0.432943);
			rgb=(0.422549, 0.092501, 0.432714);
			rgb=(0.428768, 0.094790, 0.432412);
			rgb=(0.434987, 0.097069, 0.432039);
			rgb=(0.441207, 0.099338, 0.431594);
			rgb=(0.447428, 0.101597, 0.431080);
			rgb=(0.453651, 0.103848, 0.430498);
			rgb=(0.459875, 0.106089, 0.429846);
			rgb=(0.466100, 0.108322, 0.429125);
			rgb=(0.472328, 0.110547, 0.428334);
			rgb=(0.478558, 0.112764, 0.427475);
			rgb=(0.484789, 0.114974, 0.426548);
			rgb=(0.491022, 0.117179, 0.425552);
			rgb=(0.497257, 0.119379, 0.424488);
			rgb=(0.503493, 0.121575, 0.423356);
			rgb=(0.509730, 0.123769, 0.422156);
			rgb=(0.515967, 0.125960, 0.420887);
			rgb=(0.522206, 0.128150, 0.419549);
			rgb=(0.528444, 0.130341, 0.418142);
			rgb=(0.534683, 0.132534, 0.416667);
			rgb=(0.540920, 0.134729, 0.415123);
			rgb=(0.547157, 0.136929, 0.413511);
			rgb=(0.553392, 0.139134, 0.411829);
			rgb=(0.559624, 0.141346, 0.410078);
			rgb=(0.565854, 0.143567, 0.408258);
			rgb=(0.572081, 0.145797, 0.406369);
			rgb=(0.578304, 0.148039, 0.404411);
			rgb=(0.584521, 0.150294, 0.402385);
			rgb=(0.590734, 0.152563, 0.400290);
			rgb=(0.596940, 0.154848, 0.398125);
			rgb=(0.603139, 0.157151, 0.395891);
			rgb=(0.609330, 0.159474, 0.393589);
			rgb=(0.615513, 0.161817, 0.391219);
			rgb=(0.621685, 0.164184, 0.388781);
			rgb=(0.627847, 0.166575, 0.386276);
			rgb=(0.633998, 0.168992, 0.383704);
			rgb=(0.640135, 0.171438, 0.381065);
			rgb=(0.646260, 0.173914, 0.378359);
			rgb=(0.652369, 0.176421, 0.375586);
			rgb=(0.658463, 0.178962, 0.372748);
			rgb=(0.664540, 0.181539, 0.369846);
			rgb=(0.670599, 0.184153, 0.366879);
			rgb=(0.676638, 0.186807, 0.363849);
			rgb=(0.682656, 0.189501, 0.360757);
			rgb=(0.688653, 0.192239, 0.357603);
			rgb=(0.694627, 0.195021, 0.354388);
			rgb=(0.700576, 0.197851, 0.351113);
			rgb=(0.706500, 0.200728, 0.347777);
			rgb=(0.712396, 0.203656, 0.344383);
			rgb=(0.718264, 0.206636, 0.340931);
			rgb=(0.724103, 0.209670, 0.337424);
			rgb=(0.729909, 0.212759, 0.333861);
			rgb=(0.735683, 0.215906, 0.330245);
			rgb=(0.741423, 0.219112, 0.326576);
			rgb=(0.747127, 0.222378, 0.322856);
			rgb=(0.752794, 0.225706, 0.319085);
			rgb=(0.758422, 0.229097, 0.315266);
			rgb=(0.764010, 0.232554, 0.311399);
			rgb=(0.769556, 0.236077, 0.307485);
			rgb=(0.775059, 0.239667, 0.303526);
			rgb=(0.780517, 0.243327, 0.299523);
			rgb=(0.785929, 0.247056, 0.295477);
			rgb=(0.791293, 0.250856, 0.291390);
			rgb=(0.796607, 0.254728, 0.287264);
			rgb=(0.801871, 0.258674, 0.283099);
			rgb=(0.807082, 0.262692, 0.278898);
			rgb=(0.812239, 0.266786, 0.274661);
			rgb=(0.817341, 0.270954, 0.270390);
			rgb=(0.822386, 0.275197, 0.266085);
			rgb=(0.827372, 0.279517, 0.261750);
			rgb=(0.832299, 0.283913, 0.257383);
			rgb=(0.837165, 0.288385, 0.252988);
			rgb=(0.841969, 0.292933, 0.248564);
			rgb=(0.846709, 0.297559, 0.244113);
			rgb=(0.851384, 0.302260, 0.239636);
			rgb=(0.855992, 0.307038, 0.235133);
			rgb=(0.860533, 0.311892, 0.230606);
			rgb=(0.865006, 0.316822, 0.226055);
			rgb=(0.869409, 0.321827, 0.221482);
			rgb=(0.873741, 0.326906, 0.216886);
			rgb=(0.878001, 0.332060, 0.212268);
			rgb=(0.882188, 0.337287, 0.207628);
			rgb=(0.886302, 0.342586, 0.202968);
			rgb=(0.890341, 0.347957, 0.198286);
			rgb=(0.894305, 0.353399, 0.193584);
			rgb=(0.898192, 0.358911, 0.188860);
			rgb=(0.902003, 0.364492, 0.184116);
			rgb=(0.905735, 0.370140, 0.179350);
			rgb=(0.909390, 0.375856, 0.174563);
			rgb=(0.912966, 0.381636, 0.169755);
			rgb=(0.916462, 0.387481, 0.164924);
			rgb=(0.919879, 0.393389, 0.160070);
			rgb=(0.923215, 0.399359, 0.155193);
			rgb=(0.926470, 0.405389, 0.150292);
			rgb=(0.929644, 0.411479, 0.145367);
			rgb=(0.932737, 0.417627, 0.140417);
			rgb=(0.935747, 0.423831, 0.135440);
			rgb=(0.938675, 0.430091, 0.130438);
			rgb=(0.941521, 0.436405, 0.125409);
			rgb=(0.944285, 0.442772, 0.120354);
			rgb=(0.946965, 0.449191, 0.115272);
			rgb=(0.949562, 0.455660, 0.110164);
			rgb=(0.952075, 0.462178, 0.105031);
			rgb=(0.954506, 0.468744, 0.099874);
			rgb=(0.956852, 0.475356, 0.094695);
			rgb=(0.959114, 0.482014, 0.089499);
			rgb=(0.961293, 0.488716, 0.084289);
			rgb=(0.963387, 0.495462, 0.079073);
			rgb=(0.965397, 0.502249, 0.073859);
			rgb=(0.967322, 0.509078, 0.068659);
			rgb=(0.969163, 0.515946, 0.063488);
			rgb=(0.970919, 0.522853, 0.058367);
			rgb=(0.972590, 0.529798, 0.053324);
			rgb=(0.974176, 0.536780, 0.048392);
			rgb=(0.975677, 0.543798, 0.043618);
			rgb=(0.977092, 0.550850, 0.039050);
			rgb=(0.978422, 0.557937, 0.034931);
			rgb=(0.979666, 0.565057, 0.031409);
			rgb=(0.980824, 0.572209, 0.028508);
			rgb=(0.981895, 0.579392, 0.026250);
			rgb=(0.982881, 0.586606, 0.024661);
			rgb=(0.983779, 0.593849, 0.023770);
			rgb=(0.984591, 0.601122, 0.023606);
			rgb=(0.985315, 0.608422, 0.024202);
			rgb=(0.985952, 0.615750, 0.025592);
			rgb=(0.986502, 0.623105, 0.027814);
			rgb=(0.986964, 0.630485, 0.030908);
			rgb=(0.987337, 0.637890, 0.034916);
			rgb=(0.987622, 0.645320, 0.039886);
			rgb=(0.987819, 0.652773, 0.045581);
			rgb=(0.987926, 0.660250, 0.051750);
			rgb=(0.987945, 0.667748, 0.058329);
			rgb=(0.987874, 0.675267, 0.065257);
			rgb=(0.987714, 0.682807, 0.072489);
			rgb=(0.987464, 0.690366, 0.079990);
			rgb=(0.987124, 0.697944, 0.087731);
			rgb=(0.986694, 0.705540, 0.095694);
			rgb=(0.986175, 0.713153, 0.103863);
			rgb=(0.985566, 0.720782, 0.112229);
			rgb=(0.984865, 0.728427, 0.120785);
			rgb=(0.984075, 0.736087, 0.129527);
			rgb=(0.983196, 0.743758, 0.138453);
			rgb=(0.982228, 0.751442, 0.147565);
			rgb=(0.981173, 0.759135, 0.156863);
			rgb=(0.980032, 0.766837, 0.166353);
			rgb=(0.978806, 0.774545, 0.176037);
			rgb=(0.977497, 0.782258, 0.185923);
			rgb=(0.976108, 0.789974, 0.196018);
			rgb=(0.974638, 0.797692, 0.206332);
			rgb=(0.973088, 0.805409, 0.216877);
			rgb=(0.971468, 0.813122, 0.227658);
			rgb=(0.969783, 0.820825, 0.238686);
			rgb=(0.968041, 0.828515, 0.249972);
			rgb=(0.966243, 0.836191, 0.261534);
			rgb=(0.964394, 0.843848, 0.273391);
			rgb=(0.962517, 0.851476, 0.285546);
			rgb=(0.960626, 0.859069, 0.298010);
			rgb=(0.958720, 0.866624, 0.310820);
			rgb=(0.956834, 0.874129, 0.323974);
			rgb=(0.954997, 0.881569, 0.337475);
			rgb=(0.953215, 0.888942, 0.351369);
			rgb=(0.951546, 0.896226, 0.365627);
			rgb=(0.950018, 0.903409, 0.380271);
			rgb=(0.948683, 0.910473, 0.395289);
			rgb=(0.947594, 0.917399, 0.410665);
			rgb=(0.946809, 0.924168, 0.426373);
			rgb=(0.946392, 0.930761, 0.442367);
			rgb=(0.946403, 0.937159, 0.458592);
			rgb=(0.946903, 0.943348, 0.474970);
			rgb=(0.947937, 0.949318, 0.491426);
			rgb=(0.949545, 0.955063, 0.507860);
			rgb=(0.951740, 0.960587, 0.524203);
			rgb=(0.954529, 0.965896, 0.540361);
			rgb=(0.957896, 0.971003, 0.556275);
			rgb=(0.961812, 0.975924, 0.571925);
			rgb=(0.966249, 0.980678, 0.587206);
			rgb=(0.971162, 0.985282, 0.602154);
			rgb=(0.976511, 0.989753, 0.616760);
			rgb=(0.982257, 0.994109, 0.631017);
			rgb=(0.988362, 0.998364, 0.644924);
		}
	},
	colormap/plasma/.style={colormap={plasma}{
			rgb=(0.050383, 0.029803, 0.527975);
			rgb=(0.063536, 0.028426, 0.533124);
			rgb=(0.075353, 0.027206, 0.538007);
			rgb=(0.086222, 0.026125, 0.542658);
			rgb=(0.096379, 0.025165, 0.547103);
			rgb=(0.105980, 0.024309, 0.551368);
			rgb=(0.115124, 0.023556, 0.555468);
			rgb=(0.123903, 0.022878, 0.559423);
			rgb=(0.132381, 0.022258, 0.563250);
			rgb=(0.140603, 0.021687, 0.566959);
			rgb=(0.148607, 0.021154, 0.570562);
			rgb=(0.156421, 0.020651, 0.574065);
			rgb=(0.164070, 0.020171, 0.577478);
			rgb=(0.171574, 0.019706, 0.580806);
			rgb=(0.178950, 0.019252, 0.584054);
			rgb=(0.186213, 0.018803, 0.587228);
			rgb=(0.193374, 0.018354, 0.590330);
			rgb=(0.200445, 0.017902, 0.593364);
			rgb=(0.207435, 0.017442, 0.596333);
			rgb=(0.214350, 0.016973, 0.599239);
			rgb=(0.221197, 0.016497, 0.602083);
			rgb=(0.227983, 0.016007, 0.604867);
			rgb=(0.234715, 0.015502, 0.607592);
			rgb=(0.241396, 0.014979, 0.610259);
			rgb=(0.248032, 0.014439, 0.612868);
			rgb=(0.254627, 0.013882, 0.615419);
			rgb=(0.261183, 0.013308, 0.617911);
			rgb=(0.267703, 0.012716, 0.620346);
			rgb=(0.274191, 0.012109, 0.622722);
			rgb=(0.280648, 0.011488, 0.625038);
			rgb=(0.287076, 0.010855, 0.627295);
			rgb=(0.293478, 0.010213, 0.629490);
			rgb=(0.299855, 0.009561, 0.631624);
			rgb=(0.306210, 0.008902, 0.633694);
			rgb=(0.312543, 0.008239, 0.635700);
			rgb=(0.318856, 0.007576, 0.637640);
			rgb=(0.325150, 0.006915, 0.639512);
			rgb=(0.331426, 0.006261, 0.641316);
			rgb=(0.337683, 0.005618, 0.643049);
			rgb=(0.343925, 0.004991, 0.644710);
			rgb=(0.350150, 0.004382, 0.646298);
			rgb=(0.356359, 0.003798, 0.647810);
			rgb=(0.362553, 0.003243, 0.649245);
			rgb=(0.368733, 0.002724, 0.650601);
			rgb=(0.374897, 0.002245, 0.651876);
			rgb=(0.381047, 0.001814, 0.653068);
			rgb=(0.387183, 0.001434, 0.654177);
			rgb=(0.393304, 0.001114, 0.655199);
			rgb=(0.399411, 0.000859, 0.656133);
			rgb=(0.405503, 0.000678, 0.656977);
			rgb=(0.411580, 0.000577, 0.657730);
			rgb=(0.417642, 0.000564, 0.658390);
			rgb=(0.423689, 0.000646, 0.658956);
			rgb=(0.429719, 0.000831, 0.659425);
			rgb=(0.435734, 0.001127, 0.659797);
			rgb=(0.441732, 0.001540, 0.660069);
			rgb=(0.447714, 0.002080, 0.660240);
			rgb=(0.453677, 0.002755, 0.660310);
			rgb=(0.459623, 0.003574, 0.660277);
			rgb=(0.465550, 0.004545, 0.660139);
			rgb=(0.471457, 0.005678, 0.659897);
			rgb=(0.477344, 0.006980, 0.659549);
			rgb=(0.483210, 0.008460, 0.659095);
			rgb=(0.489055, 0.010127, 0.658534);
			rgb=(0.494877, 0.011990, 0.657865);
			rgb=(0.500678, 0.014055, 0.657088);
			rgb=(0.506454, 0.016333, 0.656202);
			rgb=(0.512206, 0.018833, 0.655209);
			rgb=(0.517933, 0.021563, 0.654109);
			rgb=(0.523633, 0.024532, 0.652901);
			rgb=(0.529306, 0.027747, 0.651586);
			rgb=(0.534952, 0.031217, 0.650165);
			rgb=(0.540570, 0.034950, 0.648640);
			rgb=(0.546157, 0.038954, 0.647010);
			rgb=(0.551715, 0.043136, 0.645277);
			rgb=(0.557243, 0.047331, 0.643443);
			rgb=(0.562738, 0.051545, 0.641509);
			rgb=(0.568201, 0.055778, 0.639477);
			rgb=(0.573632, 0.060028, 0.637349);
			rgb=(0.579029, 0.064296, 0.635126);
			rgb=(0.584391, 0.068579, 0.632812);
			rgb=(0.589719, 0.072878, 0.630408);
			rgb=(0.595011, 0.077190, 0.627917);
			rgb=(0.600266, 0.081516, 0.625342);
			rgb=(0.605485, 0.085854, 0.622686);
			rgb=(0.610667, 0.090204, 0.619951);
			rgb=(0.615812, 0.094564, 0.617140);
			rgb=(0.620919, 0.098934, 0.614257);
			rgb=(0.625987, 0.103312, 0.611305);
			rgb=(0.631017, 0.107699, 0.608287);
			rgb=(0.636008, 0.112092, 0.605205);
			rgb=(0.640959, 0.116492, 0.602065);
			rgb=(0.645872, 0.120898, 0.598867);
			rgb=(0.650746, 0.125309, 0.595617);
			rgb=(0.655580, 0.129725, 0.592317);
			rgb=(0.660374, 0.134144, 0.588971);
			rgb=(0.665129, 0.138566, 0.585582);
			rgb=(0.669845, 0.142992, 0.582154);
			rgb=(0.674522, 0.147419, 0.578688);
			rgb=(0.679160, 0.151848, 0.575189);
			rgb=(0.683758, 0.156278, 0.571660);
			rgb=(0.688318, 0.160709, 0.568103);
			rgb=(0.692840, 0.165141, 0.564522);
			rgb=(0.697324, 0.169573, 0.560919);
			rgb=(0.701769, 0.174005, 0.557296);
			rgb=(0.706178, 0.178437, 0.553657);
			rgb=(0.710549, 0.182868, 0.550004);
			rgb=(0.714883, 0.187299, 0.546338);
			rgb=(0.719181, 0.191729, 0.542663);
			rgb=(0.723444, 0.196158, 0.538981);
			rgb=(0.727670, 0.200586, 0.535293);
			rgb=(0.731862, 0.205013, 0.531601);
			rgb=(0.736019, 0.209439, 0.527908);
			rgb=(0.740143, 0.213864, 0.524216);
			rgb=(0.744232, 0.218288, 0.520524);
			rgb=(0.748289, 0.222711, 0.516834);
			rgb=(0.752312, 0.227133, 0.513149);
			rgb=(0.756304, 0.231555, 0.509468);
			rgb=(0.760264, 0.235976, 0.505794);
			rgb=(0.764193, 0.240396, 0.502126);
			rgb=(0.768090, 0.244817, 0.498465);
			rgb=(0.771958, 0.249237, 0.494813);
			rgb=(0.775796, 0.253658, 0.491171);
			rgb=(0.779604, 0.258078, 0.487539);
			rgb=(0.783383, 0.262500, 0.483918);
			rgb=(0.787133, 0.266922, 0.480307);
			rgb=(0.790855, 0.271345, 0.476706);
			rgb=(0.794549, 0.275770, 0.473117);
			rgb=(0.798216, 0.280197, 0.469538);
			rgb=(0.801855, 0.284626, 0.465971);
			rgb=(0.805467, 0.289057, 0.462415);
			rgb=(0.809052, 0.293491, 0.458870);
			rgb=(0.812612, 0.297928, 0.455338);
			rgb=(0.816144, 0.302368, 0.451816);
			rgb=(0.819651, 0.306812, 0.448306);
			rgb=(0.823132, 0.311261, 0.444806);
			rgb=(0.826588, 0.315714, 0.441316);
			rgb=(0.830018, 0.320172, 0.437836);
			rgb=(0.833422, 0.324635, 0.434366);
			rgb=(0.836801, 0.329105, 0.430905);
			rgb=(0.840155, 0.333580, 0.427455);
			rgb=(0.843484, 0.338062, 0.424013);
			rgb=(0.846788, 0.342551, 0.420579);
			rgb=(0.850066, 0.347048, 0.417153);
			rgb=(0.853319, 0.351553, 0.413734);
			rgb=(0.856547, 0.356066, 0.410322);
			rgb=(0.859750, 0.360588, 0.406917);
			rgb=(0.862927, 0.365119, 0.403519);
			rgb=(0.866078, 0.369660, 0.400126);
			rgb=(0.869203, 0.374212, 0.396738);
			rgb=(0.872303, 0.378774, 0.393355);
			rgb=(0.875376, 0.383347, 0.389976);
			rgb=(0.878423, 0.387932, 0.386600);
			rgb=(0.881443, 0.392529, 0.383229);
			rgb=(0.884436, 0.397139, 0.379860);
			rgb=(0.887402, 0.401762, 0.376494);
			rgb=(0.890340, 0.406398, 0.373130);
			rgb=(0.893250, 0.411048, 0.369768);
			rgb=(0.896131, 0.415712, 0.366407);
			rgb=(0.898984, 0.420392, 0.363047);
			rgb=(0.901807, 0.425087, 0.359688);
			rgb=(0.904601, 0.429797, 0.356329);
			rgb=(0.907365, 0.434524, 0.352970);
			rgb=(0.910098, 0.439268, 0.349610);
			rgb=(0.912800, 0.444029, 0.346251);
			rgb=(0.915471, 0.448807, 0.342890);
			rgb=(0.918109, 0.453603, 0.339529);
			rgb=(0.920714, 0.458417, 0.336166);
			rgb=(0.923287, 0.463251, 0.332801);
			rgb=(0.925825, 0.468103, 0.329435);
			rgb=(0.928329, 0.472975, 0.326067);
			rgb=(0.930798, 0.477867, 0.322697);
			rgb=(0.933232, 0.482780, 0.319325);
			rgb=(0.935630, 0.487712, 0.315952);
			rgb=(0.937990, 0.492667, 0.312575);
			rgb=(0.940313, 0.497642, 0.309197);
			rgb=(0.942598, 0.502639, 0.305816);
			rgb=(0.944844, 0.507658, 0.302433);
			rgb=(0.947051, 0.512699, 0.299049);
			rgb=(0.949217, 0.517763, 0.295662);
			rgb=(0.951344, 0.522850, 0.292275);
			rgb=(0.953428, 0.527960, 0.288883);
			rgb=(0.955470, 0.533093, 0.285490);
			rgb=(0.957469, 0.538250, 0.282096);
			rgb=(0.959424, 0.543431, 0.278701);
			rgb=(0.961336, 0.548636, 0.275305);
			rgb=(0.963203, 0.553865, 0.271909);
			rgb=(0.965024, 0.559118, 0.268513);
			rgb=(0.966798, 0.564396, 0.265118);
			rgb=(0.968526, 0.569700, 0.261721);
			rgb=(0.970205, 0.575028, 0.258325);
			rgb=(0.971835, 0.580382, 0.254931);
			rgb=(0.973416, 0.585761, 0.251540);
			rgb=(0.974947, 0.591165, 0.248151);
			rgb=(0.976428, 0.596595, 0.244767);
			rgb=(0.977856, 0.602051, 0.241387);
			rgb=(0.979233, 0.607532, 0.238013);
			rgb=(0.980556, 0.613039, 0.234646);
			rgb=(0.981826, 0.618572, 0.231287);
			rgb=(0.983041, 0.624131, 0.227937);
			rgb=(0.984199, 0.629718, 0.224595);
			rgb=(0.985301, 0.635330, 0.221265);
			rgb=(0.986345, 0.640969, 0.217948);
			rgb=(0.987332, 0.646633, 0.214648);
			rgb=(0.988260, 0.652325, 0.211364);
			rgb=(0.989128, 0.658043, 0.208100);
			rgb=(0.989935, 0.663787, 0.204859);
			rgb=(0.990681, 0.669558, 0.201642);
			rgb=(0.991365, 0.675355, 0.198453);
			rgb=(0.991985, 0.681179, 0.195295);
			rgb=(0.992541, 0.687030, 0.192170);
			rgb=(0.993032, 0.692907, 0.189084);
			rgb=(0.993456, 0.698810, 0.186041);
			rgb=(0.993814, 0.704741, 0.183043);
			rgb=(0.994103, 0.710698, 0.180097);
			rgb=(0.994324, 0.716681, 0.177208);
			rgb=(0.994474, 0.722691, 0.174381);
			rgb=(0.994553, 0.728728, 0.171622);
			rgb=(0.994561, 0.734791, 0.168938);
			rgb=(0.994495, 0.740880, 0.166335);
			rgb=(0.994355, 0.746995, 0.163821);
			rgb=(0.994141, 0.753137, 0.161404);
			rgb=(0.993851, 0.759304, 0.159092);
			rgb=(0.993482, 0.765499, 0.156891);
			rgb=(0.993033, 0.771720, 0.154808);
			rgb=(0.992505, 0.777967, 0.152855);
			rgb=(0.991897, 0.784239, 0.151042);
			rgb=(0.991209, 0.790537, 0.149377);
			rgb=(0.990439, 0.796859, 0.147870);
			rgb=(0.989587, 0.803205, 0.146529);
			rgb=(0.988648, 0.809579, 0.145357);
			rgb=(0.987621, 0.815978, 0.144363);
			rgb=(0.986509, 0.822401, 0.143557);
			rgb=(0.985314, 0.828846, 0.142945);
			rgb=(0.984031, 0.835315, 0.142528);
			rgb=(0.982653, 0.841812, 0.142303);
			rgb=(0.981190, 0.848329, 0.142279);
			rgb=(0.979644, 0.854866, 0.142453);
			rgb=(0.977995, 0.861432, 0.142808);
			rgb=(0.976265, 0.868016, 0.143351);
			rgb=(0.974443, 0.874622, 0.144061);
			rgb=(0.972530, 0.881250, 0.144923);
			rgb=(0.970533, 0.887896, 0.145919);
			rgb=(0.968443, 0.894564, 0.147014);
			rgb=(0.966271, 0.901249, 0.148180);
			rgb=(0.964021, 0.907950, 0.149370);
			rgb=(0.961681, 0.914672, 0.150520);
			rgb=(0.959276, 0.921407, 0.151566);
			rgb=(0.956808, 0.928152, 0.152409);
			rgb=(0.954287, 0.934908, 0.152921);
			rgb=(0.951726, 0.941671, 0.152925);
			rgb=(0.949151, 0.948435, 0.152178);
			rgb=(0.946602, 0.955190, 0.150328);
			rgb=(0.944152, 0.961916, 0.146861);
			rgb=(0.941896, 0.968590, 0.140956);
			rgb=(0.940015, 0.975158, 0.131326);
		}
	},
	colormap/magma/.style={colormap={magma}{
			rgb=(0.001462, 0.000466, 0.013866);
			rgb=(0.002258, 0.001295, 0.018331);
			rgb=(0.003279, 0.002305, 0.023708);
			rgb=(0.004512, 0.003490, 0.029965);
			rgb=(0.005950, 0.004843, 0.037130);
			rgb=(0.007588, 0.006356, 0.044973);
			rgb=(0.009426, 0.008022, 0.052844);
			rgb=(0.011465, 0.009828, 0.060750);
			rgb=(0.013708, 0.011771, 0.068667);
			rgb=(0.016156, 0.013840, 0.076603);
			rgb=(0.018815, 0.016026, 0.084584);
			rgb=(0.021692, 0.018320, 0.092610);
			rgb=(0.024792, 0.020715, 0.100676);
			rgb=(0.028123, 0.023201, 0.108787);
			rgb=(0.031696, 0.025765, 0.116965);
			rgb=(0.035520, 0.028397, 0.125209);
			rgb=(0.039608, 0.031090, 0.133515);
			rgb=(0.043830, 0.033830, 0.141886);
			rgb=(0.048062, 0.036607, 0.150327);
			rgb=(0.052320, 0.039407, 0.158841);
			rgb=(0.056615, 0.042160, 0.167446);
			rgb=(0.060949, 0.044794, 0.176129);
			rgb=(0.065330, 0.047318, 0.184892);
			rgb=(0.069764, 0.049726, 0.193735);
			rgb=(0.074257, 0.052017, 0.202660);
			rgb=(0.078815, 0.054184, 0.211667);
			rgb=(0.083446, 0.056225, 0.220755);
			rgb=(0.088155, 0.058133, 0.229922);
			rgb=(0.092949, 0.059904, 0.239164);
			rgb=(0.097833, 0.061531, 0.248477);
			rgb=(0.102815, 0.063010, 0.257854);
			rgb=(0.107899, 0.064335, 0.267289);
			rgb=(0.113094, 0.065492, 0.276784);
			rgb=(0.118405, 0.066479, 0.286321);
			rgb=(0.123833, 0.067295, 0.295879);
			rgb=(0.129380, 0.067935, 0.305443);
			rgb=(0.135053, 0.068391, 0.315000);
			rgb=(0.140858, 0.068654, 0.324538);
			rgb=(0.146785, 0.068738, 0.334011);
			rgb=(0.152839, 0.068637, 0.343404);
			rgb=(0.159018, 0.068354, 0.352688);
			rgb=(0.165308, 0.067911, 0.361816);
			rgb=(0.171713, 0.067305, 0.370771);
			rgb=(0.178212, 0.066576, 0.379497);
			rgb=(0.184801, 0.065732, 0.387973);
			rgb=(0.191460, 0.064818, 0.396152);
			rgb=(0.198177, 0.063862, 0.404009);
			rgb=(0.204935, 0.062907, 0.411514);
			rgb=(0.211718, 0.061992, 0.418647);
			rgb=(0.218512, 0.061158, 0.425392);
			rgb=(0.225302, 0.060445, 0.431742);
			rgb=(0.232077, 0.059889, 0.437695);
			rgb=(0.238826, 0.059517, 0.443256);
			rgb=(0.245543, 0.059352, 0.448436);
			rgb=(0.252220, 0.059415, 0.453248);
			rgb=(0.258857, 0.059706, 0.457710);
			rgb=(0.265447, 0.060237, 0.461840);
			rgb=(0.271994, 0.060994, 0.465660);
			rgb=(0.278493, 0.061978, 0.469190);
			rgb=(0.284951, 0.063168, 0.472451);
			rgb=(0.291366, 0.064553, 0.475462);
			rgb=(0.297740, 0.066117, 0.478243);
			rgb=(0.304081, 0.067835, 0.480812);
			rgb=(0.310382, 0.069702, 0.483186);
			rgb=(0.316654, 0.071690, 0.485380);
			rgb=(0.322899, 0.073782, 0.487408);
			rgb=(0.329114, 0.075972, 0.489287);
			rgb=(0.335308, 0.078236, 0.491024);
			rgb=(0.341482, 0.080564, 0.492631);
			rgb=(0.347636, 0.082946, 0.494121);
			rgb=(0.353773, 0.085373, 0.495501);
			rgb=(0.359898, 0.087831, 0.496778);
			rgb=(0.366012, 0.090314, 0.497960);
			rgb=(0.372116, 0.092816, 0.499053);
			rgb=(0.378211, 0.095332, 0.500067);
			rgb=(0.384299, 0.097855, 0.501002);
			rgb=(0.390384, 0.100379, 0.501864);
			rgb=(0.396467, 0.102902, 0.502658);
			rgb=(0.402548, 0.105420, 0.503386);
			rgb=(0.408629, 0.107930, 0.504052);
			rgb=(0.414709, 0.110431, 0.504662);
			rgb=(0.420791, 0.112920, 0.505215);
			rgb=(0.426877, 0.115395, 0.505714);
			rgb=(0.432967, 0.117855, 0.506160);
			rgb=(0.439062, 0.120298, 0.506555);
			rgb=(0.445163, 0.122724, 0.506901);
			rgb=(0.451271, 0.125132, 0.507198);
			rgb=(0.457386, 0.127522, 0.507448);
			rgb=(0.463508, 0.129893, 0.507652);
			rgb=(0.469640, 0.132245, 0.507809);
			rgb=(0.475780, 0.134577, 0.507921);
			rgb=(0.481929, 0.136891, 0.507989);
			rgb=(0.488088, 0.139186, 0.508011);
			rgb=(0.494258, 0.141462, 0.507988);
			rgb=(0.500438, 0.143719, 0.507920);
			rgb=(0.506629, 0.145958, 0.507806);
			rgb=(0.512831, 0.148179, 0.507648);
			rgb=(0.519045, 0.150383, 0.507443);
			rgb=(0.525270, 0.152569, 0.507192);
			rgb=(0.531507, 0.154739, 0.506895);
			rgb=(0.537755, 0.156894, 0.506551);
			rgb=(0.544015, 0.159033, 0.506159);
			rgb=(0.550287, 0.161158, 0.505719);
			rgb=(0.556571, 0.163269, 0.505230);
			rgb=(0.562866, 0.165368, 0.504692);
			rgb=(0.569172, 0.167454, 0.504105);
			rgb=(0.575490, 0.169530, 0.503466);
			rgb=(0.581819, 0.171596, 0.502777);
			rgb=(0.588158, 0.173652, 0.502035);
			rgb=(0.594508, 0.175701, 0.501241);
			rgb=(0.600868, 0.177743, 0.500394);
			rgb=(0.607238, 0.179779, 0.499492);
			rgb=(0.613617, 0.181811, 0.498536);
			rgb=(0.620005, 0.183840, 0.497524);
			rgb=(0.626401, 0.185867, 0.496456);
			rgb=(0.632805, 0.187893, 0.495332);
			rgb=(0.639216, 0.189921, 0.494150);
			rgb=(0.645633, 0.191952, 0.492910);
			rgb=(0.652056, 0.193986, 0.491611);
			rgb=(0.658483, 0.196027, 0.490253);
			rgb=(0.664915, 0.198075, 0.488836);
			rgb=(0.671349, 0.200133, 0.487358);
			rgb=(0.677786, 0.202203, 0.485819);
			rgb=(0.684224, 0.204286, 0.484219);
			rgb=(0.690661, 0.206384, 0.482558);
			rgb=(0.697098, 0.208501, 0.480835);
			rgb=(0.703532, 0.210638, 0.479049);
			rgb=(0.709962, 0.212797, 0.477201);
			rgb=(0.716387, 0.214982, 0.475290);
			rgb=(0.722805, 0.217194, 0.473316);
			rgb=(0.729216, 0.219437, 0.471279);
			rgb=(0.735616, 0.221713, 0.469180);
			rgb=(0.742004, 0.224025, 0.467018);
			rgb=(0.748378, 0.226377, 0.464794);
			rgb=(0.754737, 0.228772, 0.462509);
			rgb=(0.761077, 0.231214, 0.460162);
			rgb=(0.767398, 0.233705, 0.457755);
			rgb=(0.773695, 0.236249, 0.455289);
			rgb=(0.779968, 0.238851, 0.452765);
			rgb=(0.786212, 0.241514, 0.450184);
			rgb=(0.792427, 0.244242, 0.447543);
			rgb=(0.798608, 0.247040, 0.444848);
			rgb=(0.804752, 0.249911, 0.442102);
			rgb=(0.810855, 0.252861, 0.439305);
			rgb=(0.816914, 0.255895, 0.436461);
			rgb=(0.822926, 0.259016, 0.433573);
			rgb=(0.828886, 0.262229, 0.430644);
			rgb=(0.834791, 0.265540, 0.427671);
			rgb=(0.840636, 0.268953, 0.424666);
			rgb=(0.846416, 0.272473, 0.421631);
			rgb=(0.852126, 0.276106, 0.418573);
			rgb=(0.857763, 0.279857, 0.415496);
			rgb=(0.863320, 0.283729, 0.412403);
			rgb=(0.868793, 0.287728, 0.409303);
			rgb=(0.874176, 0.291859, 0.406205);
			rgb=(0.879464, 0.296125, 0.403118);
			rgb=(0.884651, 0.300530, 0.400047);
			rgb=(0.889731, 0.305079, 0.397002);
			rgb=(0.894700, 0.309773, 0.393995);
			rgb=(0.899552, 0.314616, 0.391037);
			rgb=(0.904281, 0.319610, 0.388137);
			rgb=(0.908884, 0.324755, 0.385308);
			rgb=(0.913354, 0.330052, 0.382563);
			rgb=(0.917689, 0.335500, 0.379915);
			rgb=(0.921884, 0.341098, 0.377376);
			rgb=(0.925937, 0.346844, 0.374959);
			rgb=(0.929845, 0.352734, 0.372677);
			rgb=(0.933606, 0.358764, 0.370541);
			rgb=(0.937221, 0.364929, 0.368567);
			rgb=(0.940687, 0.371224, 0.366762);
			rgb=(0.944006, 0.377643, 0.365136);
			rgb=(0.947180, 0.384178, 0.363701);
			rgb=(0.950210, 0.390820, 0.362468);
			rgb=(0.953099, 0.397563, 0.361438);
			rgb=(0.955849, 0.404400, 0.360619);
			rgb=(0.958464, 0.411324, 0.360014);
			rgb=(0.960949, 0.418323, 0.359630);
			rgb=(0.963310, 0.425390, 0.359469);
			rgb=(0.965549, 0.432519, 0.359529);
			rgb=(0.967671, 0.439703, 0.359810);
			rgb=(0.969680, 0.446936, 0.360311);
			rgb=(0.971582, 0.454210, 0.361030);
			rgb=(0.973381, 0.461520, 0.361965);
			rgb=(0.975082, 0.468861, 0.363111);
			rgb=(0.976690, 0.476226, 0.364466);
			rgb=(0.978210, 0.483612, 0.366025);
			rgb=(0.979645, 0.491014, 0.367783);
			rgb=(0.981000, 0.498428, 0.369734);
			rgb=(0.982279, 0.505851, 0.371874);
			rgb=(0.983485, 0.513280, 0.374198);
			rgb=(0.984622, 0.520713, 0.376698);
			rgb=(0.985693, 0.528148, 0.379371);
			rgb=(0.986700, 0.535582, 0.382210);
			rgb=(0.987646, 0.543015, 0.385210);
			rgb=(0.988533, 0.550446, 0.388365);
			rgb=(0.989363, 0.557873, 0.391671);
			rgb=(0.990138, 0.565296, 0.395122);
			rgb=(0.990871, 0.572706, 0.398714);
			rgb=(0.991558, 0.580107, 0.402441);
			rgb=(0.992196, 0.587502, 0.406299);
			rgb=(0.992785, 0.594891, 0.410283);
			rgb=(0.993326, 0.602275, 0.414390);
			rgb=(0.993834, 0.609644, 0.418613);
			rgb=(0.994309, 0.616999, 0.422950);
			rgb=(0.994738, 0.624350, 0.427397);
			rgb=(0.995122, 0.631696, 0.431951);
			rgb=(0.995480, 0.639027, 0.436607);
			rgb=(0.995810, 0.646344, 0.441361);
			rgb=(0.996096, 0.653659, 0.446213);
			rgb=(0.996341, 0.660969, 0.451160);
			rgb=(0.996580, 0.668256, 0.456192);
			rgb=(0.996775, 0.675541, 0.461314);
			rgb=(0.996925, 0.682828, 0.466526);
			rgb=(0.997077, 0.690088, 0.471811);
			rgb=(0.997186, 0.697349, 0.477182);
			rgb=(0.997254, 0.704611, 0.482635);
			rgb=(0.997325, 0.711848, 0.488154);
			rgb=(0.997351, 0.719089, 0.493755);
			rgb=(0.997351, 0.726324, 0.499428);
			rgb=(0.997341, 0.733545, 0.505167);
			rgb=(0.997285, 0.740772, 0.510983);
			rgb=(0.997228, 0.747981, 0.516859);
			rgb=(0.997138, 0.755190, 0.522806);
			rgb=(0.997019, 0.762398, 0.528821);
			rgb=(0.996898, 0.769591, 0.534892);
			rgb=(0.996727, 0.776795, 0.541039);
			rgb=(0.996571, 0.783977, 0.547233);
			rgb=(0.996369, 0.791167, 0.553499);
			rgb=(0.996162, 0.798348, 0.559820);
			rgb=(0.995932, 0.805527, 0.566202);
			rgb=(0.995680, 0.812706, 0.572645);
			rgb=(0.995424, 0.819875, 0.579140);
			rgb=(0.995131, 0.827052, 0.585701);
			rgb=(0.994851, 0.834213, 0.592307);
			rgb=(0.994524, 0.841387, 0.598983);
			rgb=(0.994222, 0.848540, 0.605696);
			rgb=(0.993866, 0.855711, 0.612482);
			rgb=(0.993545, 0.862859, 0.619299);
			rgb=(0.993170, 0.870024, 0.626189);
			rgb=(0.992831, 0.877168, 0.633109);
			rgb=(0.992440, 0.884330, 0.640099);
			rgb=(0.992089, 0.891470, 0.647116);
			rgb=(0.991688, 0.898627, 0.654202);
			rgb=(0.991332, 0.905763, 0.661309);
			rgb=(0.990930, 0.912915, 0.668481);
			rgb=(0.990570, 0.920049, 0.675675);
			rgb=(0.990175, 0.927196, 0.682926);
			rgb=(0.989815, 0.934329, 0.690198);
			rgb=(0.989434, 0.941470, 0.697519);
			rgb=(0.989077, 0.948604, 0.704863);
			rgb=(0.988717, 0.955742, 0.712242);
			rgb=(0.988367, 0.962878, 0.719649);
			rgb=(0.988033, 0.970012, 0.727077);
			rgb=(0.987691, 0.977154, 0.734536);
			rgb=(0.987387, 0.984288, 0.742002);
			rgb=(0.987053, 0.991438, 0.749504);
		}
	}
}

\journal{arXiv}

\begin{document}
\begin{frontmatter}

\title{Numerical Simulation of Phase Transition with the Hyperbolic Godunov-Peshkov-Romenski Model}

\author[]{Pascal Mossier \protect\affmark[a,$*$]}
\ead{pascal.mossier@iag.uni-stuttgart.de}
\author[]{Steven J\"ons \affmark[a]}
\author[]{Simone Chiocchetti \affmark[a,b]}
\author[]{Andrea D. Beck \affmark[a]}
\author[]{Claus-Dieter Munz \affmark[a]}

\affiliation[1]{organization={Institute of Aerodynamics and Gas Dynamics, University of Stuttgart},
            addressline={Pfaffenwaldring 21}, 
            city={Stuttgart},
            postcode={70569}, 
            country={Germany}}

\affiliation[2]{organization={Division of Mathematics, University of Cologne},
			addressline={Weyertal 86-90}, 
			city={Cologne},
			postcode={50931}, 
			country={Germany}}
\begin{abstract}
	In this paper, a thermodynamically consistent solution of the interfacial Riemann problem for the 
	first-order hyperbolic continuum model of Godunov, Peshkov and Romenski (GPR model) is presented. In the presence of phase transition, interfacial physics
	are governed by molecular interaction on a microscopic scale, beyond the scope of the macroscopic continuum model in the bulk phases. 
	The developed two-phase Riemann solvers tackle this multi-scale problem, by incorporating a local thermodynamic model to predict the interfacial 
	entropy production. Using phenomenological relations of non-equilibrium thermodynamics, interfacial mass and heat fluxes are derived from the entropy production
	and provide closure at the phase boundary. We employ the proposed Riemann solvers in an efficient sharp interface level-set Ghost-Fluid 
	framework to provide coupling conditions at phase interfaces under phase transition. As a single-phase benchmark, a Rayleigh-B\'{e}nard convection is studied to compare the hyperbolic thermal relaxation formulation of the GPR model 
	against the hyperbolic-parabolic Euler-Fourier system. The novel interfacial Riemann solvers are validated against
	molecular dynamics simulations of evaporating shock tubes with the Lennard-Jones shifted and truncated potential. On a macroscopic scale,
	evaporating shock tubes are computed for the material n-Dodecane and compared against Euler-Fourier results. Finally, the efficiency and robustness of 
	the scheme is demonstrated with shock-droplet interaction simulations that involve both phase transfer and surface tension, while featuring severe interface 
	deformations. 
\end{abstract}			


\end{frontmatter}

\section{Introduction}
\label{sec:Introduction}
The process of phase transition is a defining characteristic of interfacial flows and at the heart of fundamental environmental processes 
like the water cycle. In engineering, the understanding and reliable prediction of multiphase flows with phase transition is essential
e.g. in cooling circuits and combustion chambers of current aeronautical propulsion systems. Here, phase transition
often occurs under extreme ambient conditions close to the critical point and in the presence of strong thermodynamic non-equilibrium. 

Despite the indisputable relevance of interfacial flows with phase transition, high-fidelity simulations of such phenomena on a macroscopic scale 
remain a formidable challenge for current numerical methods due to the inherent multi-scale character of the problem. While the fluid in the bulk phases can be described
with continuum models, physical effects at the evaporating phase boundary are governed by molecular interaction on a microscopic length scale 
where the continuum assumption is lost. Further, in near critical conditions or under non-equilibrium, the fluids are strongly affected by compressible effects,
which leads to a tight coupling between hydrodynamics and thermodynamics. 

In literature, interfacial flows are studied with a variety of models depending on the scale of the problem of interest. 
Molecular dynamics (MD) simulations describe multi-phase flows through the interaction of molecules on a microscopic scale \cite{Dang1997,Kotsalis2004,Heinen2019}. 
Thereby, the material properties depend on intermolecular attractive and repulsive forces governed by a suitable potential. 
While MD simulations capture fluid behavior and interfacial physics intrinsically, they are restricted
to small-scale problems due to the immense number of molecules.

On a macroscopic continuum scale, two main strategies can be distinguished to model multi-phase flows: diffuse interface models and sharp interface models. 
While diffuse interface methods like the Navier-Stokes-Korteweg equations \cite{Anderson1998} and the Baer-Nunziato equations \cite{Bear1986,Kapila2001} model
the phase interface as a smooth transition layer of finite thickness, the sharp interface approach assumes a discontinuous transition of fluid properties
across a sharp interface of zero thickness. As a consequence, diffuse interface methods have to resolve the interface with sufficient accuracy to capture the
local physical behavior and are thus restricted to small-scale problems. Therefore, the present paper focuses on a sharp interface approach to study interfacial 
flows on a macroscopic scale, where the local physical behavior is included in the interfacial jump conditions. 

Multi-phase simulations with the sharp interface method rely on two essential building blocks: an interface tracking algorithm and a consistent coupling of the 
bulk phases across the phase interface. Following Sussman et al. \cite{Sussman1994}, we employ the level-set method to track the evolution
of the interface position and geometry. Coupling of the bulk phases at the phase boundary is implemented via the 
Ghost-Fluid idea of Fedkiw et al. \cite{Fedkiw1999}, which relies on the definition of ghost states at the interface. 

In this work, we follow the concept of Merkle and Rohde \cite{Merkle2007} to solve an interfacial Riemann problem to supply the required ghost states. 
Contrary to the single-phase case, where a multitude of Riemann solvers are available \cite{Toro1999}, the construction of thermodynamically consistent 
two-phase Riemann solvers in the presence of phase transition is still the subject of ongoing research. 

The main challenge arises due to a breakdown of the continuum
assumption across the phase interface. In the equation of state (EOS), this manifests as an unphysical spinodal region which is characterized by 
a non-convex behavior of the EOS and thus leads to imaginary eigenvalues. As Menikoff and Plohr pointed out \cite{Menikoff1989}, this results in anomalous 
wave structures and non-uniqueness of the solution of the Riemann problem. A possible approach to find a unique, admissible solution is the kinetic 
relation as proposed by Abeyaratne and Knowles \cite{Abeyaratne1991}, which enforces the correct amount of entropy production due to phase change.

The solution to the interfacial Riemann problem is further complicated in the presence of phase transition since  
the mechanism of heat transfer is essential to provide the latent heat of vaporization \cite{Hantke2019}. With heat transfer 
commonly modeled by the hyperbolic parabolic Euler-Fourier model, this results in a loss of self-similarity. In recent works by 
Hitz et al. \cite{Hitz2020,Hitz2021} and J\"ons et al. \cite{Joens2023,Joens2023a}, this issue was circumvented by solving the Riemann problem first without heat transfer
and then imposing the heat fluxes, obtained from an evaporation model, on the resulting interfacial fluxes. 

We choose a different approach and use the inviscid Godunov-Peshkov-Romenski (GPR) continuum model \cite{Dumbser2016}, which provides a 
first-order hyperbolic formulation for compressible, heat-conducting fluids, based on the work of Malyshev and Romenski \cite{Romenski1986}. 
A key advantage of the GPR approach is the treatment of irreversible, dissipative effects of heat conduction via algebraic source terms. This allows for the incorporation
of heat transfer effects in the solution of the Riemann solver. M\"uller et al. developed an approximate two-phase Riemann solver for the GPR system in \cite{Mueller2023,MuellerPhD}.
We refer the reader also to \cite{Thein2022}, where the authors formulate a Riemann solver for two-phase flow in the context of symmetric hyperbolic thermodynamically
compatible (SHTC) models, to which the GPR model pertains.

In the present paper, a novel approximate Riemann solver and a second simplified version are formulated for the interfacial Riemann problem under phase transition. 
Both solvers rely on a local thermodynamic model by Cipolla et al. \cite{Cipolla1974} that predicts the interfacial entropy production based on kinetic theory. 
Using phenomenological relations from non-equilibrium thermodynamics, the interfacial mass and heat flux are derived from the estimated entropy production. 
With these additional conditions for the thermodynamic interfacial fluxes, a closed non-linear equation system is obtained, which can be solved iteratively towards
the correct entropy solution. Crucially, the local phase transition model controls the entropy production and the dissipative effects of heat transition and thus 
allows for a thermodynamically consistent incorporation of the source terms of the GPR equation system in the solution process. 

The present work combines the novel interface solvers with an efficient semi-analytical source term integration scheme \cite{Chiocchetti2023},
which allows for accurate and robust source term treatment in the stiff regime. The resulting scheme is implemented in the hp-adaptive Discontinuous Galerkin (DG) 
multi-phase code \textit{FLEXI} \cite{fechter2015,Mueller2020,Joens2021,Appel2021,Mossier2023} to study interfacial flows with phase transition in multiple space dimensions. 

This paper is structured as follows: In Section \ref{sec:Equations} we recapitulate the GPR equation system for an inviscid, heat conducting continuum in the bulk phases and derive a 
kinetic relation that defines the entropy production across the interface. To close the system, expressions for the interfacial mass and heat 
flux are provided, based on kinetic theory and phenomenological force flux relations. In Section \ref{sec:Numerics}, the sharp interface framework 
is outlined briefly and two approximate two-phase Riemann solvers with phase transfer are formulated. Furthermore, the semi-analytical source 
term integration for the relaxation term of the hyperbolic heat transfer equation is addressed. Finally, in Section \ref{sec:Results} we apply the 
framework to a selection of numerical test cases starting with a single phase Rayleigh-B\'{e}nard convection to compare the GPR system against the 
Euler-Fourier equations. The Riemann solvers are applied for evaporating shock tube computations and validated against molecular dynamics data and 
Euler-Fourier computations. Finally, the robustness and efficiency of the scheme is demonstrated with two-dimensional shock-droplet interaction 
simulations, that exhibit severe interface deformations. The paper closes with a summary and conclusion in Section \ref{sec:Conclusion}. 
\section{Governing Equations}
\label{sec:Equations}
In the present work, we study compressible, two-phase flows with phase transition. We restrict our investigation to inviscid single-component fluids.%, thus only imiscible liquid and vapor phases of a single fluid are considered. 
The bulk phases are separated by a phase interface of zero thickness according to the sharp interface approach. Therefore, the computational domain $\Omega$ consists of a liquid subdomain $\Omega_l$ and a vapor subdomain $\Omega_v$, separated by a phase interface $\Gamma=\Omega_v\cap\Omega_l$. 
\subsection{Continuum Model of the Bulk Fluid}
As a continuum model for compressible, heat-conducting fluids, we chose the inviscid GPR equation system as found in \cite{Dumbser2016}, 
where the hyperbolic formulation of heat transfer developed by Malyshev and Romenski \cite{Romenski1986} is
found in conjunction with the model of mechanics given in \cite{Peshkov2014}. The system is defined by
\begin{subequations}
	\begin{align}
		\frac{\partial \rho}{\partial t}&+\nabla \cdot \left(\rho \uB \right) = 0,                                 \label{eq:HPR_equations_1} \\
		\frac{\partial \rho \uB}{\partial t}&+\nabla \cdot \left(\rho \uB \otimes \uB + p \mathbb{I} \right) = 0,  \label{eq:HPR_equations_2} \\
		\frac{\partial \rho \chi}{\partial t}&+\nabla \cdot \left(\rho \uB \chi  + \alpha^2 \jB \right) = -\frac{\rho}{\theta(\tau)T}\alpha^2 \jB\cdot\jB,    \label{eq:HPR_equations_3} \\
		\frac{\partial \rho \jB}{\partial t}&+\nabla \cdot \left(\rho \jB \otimes \uB + T \mathbb{I} \right)  = -\frac{\alpha^2 \rho \jB}{\theta(\tau)}, \label{eq:HPR_equations_4}
	\end{align}
	\label{eq:HPR_equations}
\end{subequations}
and satisfies an additional conservation equation for the total energy 
\begin{equation}
	\frac{\partial \rho e}{\partial t}+\nabla \cdot \left((\rho e + p) \uB + q\right) = 0.
	\label{eq:HPR_energy_equation}
\end{equation}
Even though the GPR system is derived for the entropy $\chi$ as an unknown and energy conservation as a consequence to obtain a thermodynamically compatible formulation \cite{Dumbser2016}, 
we solve the equation system with the energy instead of the entropy for practical numerical studies within this work. 
With the energy as an unknown, the GPR system can be expressed in matrix-vector notation as
\begin{equation}
\frac{\partial\QB}{\partial t} + \nabla_{\xB} \cdot \FB(\QB)=\SB(\QB)\quad \quad \text{in}  \quad \Omega \times [0,T],
\label{eq:HPR_matrix_vector}
\end{equation}
with the state vector $\QB$, the physical flux vector $\FB$ and the algebraic source term $\SB$, given in terms of the density $\rho$, velocity $\uB=(u_1,u_2,u_3)^T$, 
total energy per unit mass $e$, the pressure $p$, the temperature $T$ and the thermal impulse per unit mass $\jB=(j_1,j_2,j_3)^T$. 
The thermal heat flux $\qB$ is related to the temperature and the thermal impulse via the constitutive relation 
\begin{equation}
	\qB=\alpha^2Tj.
	\label{eq:heatflux}
\end{equation}
Thereby, the parameter $\alpha$ is connected to the propagation speed of the thermo-acoustic waves $c_h$ by 
\begin{equation}
	c_h=\frac{\alpha}{\rho}\sqrt{\frac{T}{c_v}},
	\label{eq:c_h} 
\end{equation}
with the specific heat capacity at constant volume $c_v$.
A key aspect of the GPR method is the definition of the specific total energy as a potential $e(\rho,\chi,\uB,\jB)$, which ensures 
a thermodynamically consistent formulation of the overdetermined equation system \eqref{eq:HPR_equations} and \eqref{eq:HPR_energy_equation}. 
The energy $e$ is obtained as the sum of the internal energy $u$, the kinetic energy $\frac{1}{2}\uB\cdot\uB$ and the mesoscopic 
non-equilibrium part of the energy $\frac{1}{2}\alpha^2 \jB\cdot\jB$ associated to thermal non-equilibrium:
\begin{equation}
	e = \epsilon+\frac{1}{2}\alpha\jB\cdot\jB+\frac{1}{2}\uB\cdot\uB.
	\label{eq:totalenergy} 
\end{equation}
For a detailed discussion of the energy potential, the reader is referred to Dumbser et al. \cite{Dumbser2016}.
The equation system \eqref{eq:HPR_equations} is closed with an EOS that relates the pressure with the density and specific internal energy 
\begin{equation}
	p=p(\rho,\epsilon).
	\label{eq:eos} 
\end{equation}
Our numerical framework provides algebraic EOS, like the ideal or stiffened gas EOS, as well as cubic EOS, like the Peng-Robinson EOS. 
Furthermore, multi-parameter EOS from the fluid library \textit{CoolProp} are available and can be accessed via the efficient 
tabulation approach of F\"oll et al. \cite{Foell2019}.

Finally, the source term $\SB$ contains a scalar function $\theta(\tau)$ that depends on the thermal impulse and relaxation time $\tau$. 
The remaining free parameter $\theta(\tau)$ can be defined as 
\begin{equation}
	\theta(\tau) = \tau\alpha^2\frac{\rho}{\rho_0}\frac{T_0}{T},
	\label{eq:Theta}
\end{equation}
with the reference density $\rho_0$ and reference $T_0$, set to to the initial conditions. Dumbser et at. \cite{Dumbser2016} 
demonstrated with an asymptotic analysis that this particular choice of $\theta(\tau)$ recovers the Fourier law in the stiff limit
\begin{equation}
	\qB=\alpha^2Tj=\tau\alpha^2\frac{T_0}{\rho_0}\nabla T:=-\lambda\nabla T.
	\label{eq:hpr_fourier_limit} 
\end{equation}
Thus, the relaxation time $\tau$ can be related to the thermal conductivity $\lambda$ through
\begin{equation}
	\lambda=\alpha^2\tau\frac{T_0}{\rho_0},
	\label{eq:lambda} 
\end{equation}
for vanishing relaxation times $\tau$. 

To simplify the evaluation of the eigenvalues of the GPR system, the parameter $\alpha$ is assumed to be constant within each phase
and determined once at the beginning of the computation by equation \eqref{eq:lambda}. Therefore, the relaxation time 
$\tau$ needs to be provided at the start of the computation. While the thermal relaxation time was set manually tuned to a stiff regime 
by Dumbser et al. in \cite{Dumbser2016,Dumbser2017a}, we follow the approach of M{\"u}ller et al. and choose between 
the kinetic theory-based model of Jordan \cite{Jordan2014} 
\begin{equation}
	\tau=\frac{3}{c_s^2}\frac{\lambda}{\rho c_v}\frac{T}{T_0}\frac{\rho_0}{\rho},
	\label{eq:tau_kin} 
	\end{equation}
and a thermomass theory-based model \cite{Chester1963,Kazemi2017} that predicts a relaxation time
\begin{equation}
	\tau=\frac{\lambda}{\rho c_v}\frac{1}{2 c_p T},
	\label{eq:tau_thermo} 
\end{equation}
with $c_s$ denoting the speed of sound. It is to be emphasized, that the validity of both models for real 
materials and a wide temperature range is highly questionable. However, in the absence of physically sound
models for $\tau$ in literature, the given models allow for a parameter-free closure and provide acceptable results
as long as a relaxation time in the stiff regime is ensured. 

\subsection{Thermodynamics of Phase Transition}
\label{subsec:evaporation_model}
At a sharp phase interface $\Gamma$ separating two immiscible phases, interfacial physics 
are governed by the exchange of mass, momentum and energy. From the Rankine-Hugoniot conditions,
a set of jump conditions can be defined for mass, momentum, energy and thermal impulse. 
For the GPR continuum model, they are obtained in the interface-normal direction as
\begin{subequations}
	\begin{align}
		\llbracket \dot{m} \rrbracket & = 0,                                                                     \label{eq:interface_jump_1}\\
		\dot{m}\llbracket u \rrbracket+\llbracket p \rrbracket&=\Delta p_{\sigma},								 \label{eq:interface_jump_2}\\
		\dot{m}\llbracket e \rrbracket+\llbracket up \rrbracket+\llbracket q \rrbracket&=\sB \Delta p_{\sigma},  \label{eq:interface_jump_3}\\
		\dot{m}\llbracket j \rrbracket+\llbracket T \rrbracket&=\Delta T, 								         \label{eq:interface_jump_4}
	\end{align}
	\label{eq:interface_jump}
\end{subequations}
with the definition of a jump operator $\llbracket z \rrbracket=z_v-z_l$ for an arbitrary quantity $z$. 
Here, the expression $\sB = \nB\cdot\sInt$ denotes the velocity of the phase boundary in an interface-normal reference space.
To account for surface tension, a pressure jump $\Delta p_{\sigma}$ is included in the momentum and energy equation with 
\begin{equation}
	\Delta p_{\sigma}=2\kappa\sigma,
	\label{eq:young_laplace_law}
\end{equation}
according to the Young-Laplace law. In the presence of phase transition, evaporation or condensation respectively
drives an interfacial mass flux denoted as $\dot{m}$. A further consequence of phase transition is an interfacial 
temperature jump $\Delta T$ in the thermal impulse equation. 

In the GPR system, dissipative effects of heat conduction are captured by an algebraic relaxation source term
\begin{equation}
	S_{\rho j}=-\frac{\alpha ^2 \rho j}{\Theta(\tau)}
	\label{eq:source_thermal}
\end{equation}
in the thermal impulse balance. Consequently, the term $S_{\rho j}$ reappears in the source term of the entropy balance equation  
\begin{equation}
	S_{\rho \chi}=-\frac{\rho}{\Theta(\tau) T} \alpha^4 \jB\jB = - S_{\rho j} \frac{\alpha^2 \jB}{T},
	\label{eq:source_entropy}
\end{equation}
where it describes the entropy production due to heat conduction. Since the GPR continuum model can not be expected to predict the
entropy production of phase transition, a local thermodynamic model is employed to obtain the correct entropy solution. 
In that sense, the temperature jump $\Delta T$ serves as an additional degree of freedom in the jump conditions, which
allows for an interfacial entropy production. Unknown a priori, it is determined through a surrogate phase transition model. 
It is to remark that the assumption of an interfacial temperature jump in the presence of phase change
is in line with experimental findings e.g. by Gatapova et al. \cite{Gatapova2017} or Kazemi et al. \cite{Kazemi2017}. 

To fulfill the second law of thermodynamics, an entropy jump condition can be derived from the entropy balance \eqref{eq:HPR_equations_3} as
\begin{equation}
\dot{m} \llbracket \chi\rrbracket + \llbracket \frac{\dot{q}}{T} \rrbracket=\chi^\Gamma \quad \text{, with} \quad \chi^\Gamma\geq 0,
\label{eq:entropy}
\end{equation}
where $\chi^\Gamma$ denotes the entropy production rate. This kinetic relation serves as a link between the jump conditions of the continuum model and 
a suitable phase transition model. Following \cite{MuellerPhD}, equation \eqref{eq:entropy} can be expressed in terms of the
Gibbs energy per unit mass $g$ and the enthalpy $h$ per unit mass as
\begin{equation}
	\dot{m} \llbracket \frac{-g+h_v}{T} \rrbracket + q_v \llbracket \frac{1}{T} \rrbracket = \chi^\Gamma,
	\label{eq:entropy_rewritten}
\end{equation}
with $h_v$ and $q_v$ denoting the enthalpy per unit mass and the heat flux evaluated at the vapor side of the phase boundary.
This relation can now be reformulated with the tools of non-equilibrium thermodynamics \cite{Lebon2008}. 
From equation \eqref{eq:entropy_rewritten} we can deduce the thermodynamic fluxes $\dot{m}$ and $q_v$ and the thermodynamic 
forces $\llbracket \frac{-g+h_v}{T} \rrbracket$ and $\llbracket \frac{1}{T} \rrbracket$. 
Assuming linear dependencies between forces and fluxes according to Onsager's theory, phenomenological relations for the 
thermodynamic fluxes can be derived as
\begin{align}
\dot{m}&=L_{mm} \llbracket \frac{-g+h_v}{T} \rrbracket + L_{em} \llbracket \frac{1}{T} \rrbracket, \\
q_v&=L_{me} \llbracket \frac{-g+h_v}{T} \rrbracket + L_{ee} \llbracket \frac{1}{T} \rrbracket,
\label{eq:onsager_relation}
\end{align}
with $L_{mm}$,$L_{me}$,$L_{em}$,$L_{ee}$ denoting the Onsager coefficients. 
With the Onsager reciprocal relation $L_{me}=L_{me}$, the number of the yet unknown Onsager coefficients can be reduced to three. 
The modeling of the non-equilibrium phase transition process is thus reduced to a closure problem for the three remaining Onsager coefficients. 
In this publication we follow the approach of J\"ons et al. \cite{Joens2023} and apply a model from kinetic theory derived by Cipolla et al. 
\cite{Cipolla1974}. Herein, the Onsager coefficients are reported as
\begin{align}
L_{mm}&=\frac{-k_2}{k_1 k_2-k_3^2}\rho_v\sqrt{\frac{2T_l}{R}},  \\
L_{me}=L_{em}&=\frac{-k_3}{k_1 k_2-k_3^2}\rho_v T_l\sqrt{2T_l}, \\
L_{ee}&=\frac{-k_1}{k_1 k_2- k_3^2}p_s(T_l)T_l\sqrt{2T_l},
\label{eq:cipolla}
\end{align} 
with the ideal gas constant $R$ and the saturation pressure at the liquid temperature $p_s(T_l)$. 
The coefficients $k_1$, $k_2$ and $k_3$ are defined as
\begin{align}
k_1 &= \frac{9}{8}\sqrt{\pi}\left(\frac{1}{2}+\frac{16}{9\pi}\right)-\sqrt{\pi}\frac{1-\sigma_c}{\sigma_c}, \\
k_2&=\frac{1}{2}\sqrt{\pi}\left(\frac{1}{2}+\frac{52}{25\pi}\right), \\
k_3&=\frac{1}{4}\sqrt{\pi}\left(\frac{1}{2}+\frac{8}{5\pi}\right),
\end{align} 
with the condensation coefficient $\sigma_c$. We determine the condensation coefficient with a model of Nagayama et al. \cite{Nagayama2003} as
\begin{equation}
\sigma_c=\left(1-\sqrt[3]{\frac{\nu_l}{\nu_v}}\right)\exp \left(-\frac{1}{2}\frac{1}{\sqrt[3]{\nu_l/\nu_v}-1}\right),
\label{eq:condCoeff}
\end{equation}
with 
\begin{align}
\nu_l &= \frac{1}{\rho_l}-\frac{1}{3\rho_{c}}, \\
\nu_v &= \frac{1}{\rho_v}-\frac{1}{3\rho_{c}}, \\
\end{align}
and the critical density $\rho_{c}$. We want to emphasize that a different choice for a microscopic model for the Onsager coefficients 
or even a different choice for an evaporation model, like the Hertz Knudsen model \cite{Hertz1882,Knudsen1915} that predicts the thermodynamic 
fluxes directly is possible within the present framework and does not affect the construction of the two-phase Riemann solvers. 

\subsection{Level-set Interface Tracking}
The interface tracking algorithm is an essential building block of a sharp interface framework. We deduce the interface position and geometry 
following Sussman et al. \cite{Sussman1994} from a level-set function $\Phi$ that is advected with the velocity $\sInt=(s_1,s_2,s_3)$ according to
\begin{equation}
\frac{\partial\Phi}{\partial t}+\sInt\cdot\nabla_{\xB}\Phi=0.
\label{eq:levelset}
\end{equation}
The velocity $\sInt$ is determined by solving a two-phase Riemann problem at the phase interface and subsequently extrapolated to the volume 
by solving a Hamilton-Jacobi type equation, according to Peng et al. \cite{Peng1998}. The signed distance property of the level-set function 
is maintained with the level-set reinitialization procedure of Sussman et al. \cite{Sussman1994}. Geometric properties like the local 
normal vector $\nB$ and curvature $\kappa$ of the phase interface can obtained in terms of derivatives of the level-set field. For a detailed 
description of the interface tracking algorithm, the reader is referred to \cite{FechterPhd,Joens2021}.

\section{Numerical Method}
\label{sec:Numerics}
The sharp interface framework used in this work requires three major building blocks: the bulk fluid solver, the level-set interface 
tracking algorithm and a Ghost-Fluid coupling at the interface via a two-phase Riemann solver. In this Section, we provide a brief 
overview of the bulk fluid solver and the interface tracking algorithm. For a thorough description of these operators, the reader is 
referred to \cite{FechterPhd,MuellerPhD,Zeifang2020}. 
Since the bulk flow is modeled with the GPR equation system, additional care is required when treating the stiff thermal relaxation 
source term. In the present framework, we implemented the semi-analytical source term integration, recently introduced in \cite{Chiocchetti2023}. 
The main focus of this Section is dedicated to the construction of two novel approximate interfacial Riemann 
solvers for the GPR equation system. 

\subsection{The Bulk Fluid Solver}
\label{subsec:BulkSolver}
The GPR equation system \eqref{eq:HPR_equations}, governing the bulk flow, is discretized in space with the high-order discontinuous Galerkin 
spectral element method (DGSEM)\cite{kopriva2009,flexi}. Hereby the computational domain $\Omega\subset  \mathbb{R}^3$ is divided in $K\in \mathbb{N}$ 
non-overlapping hexaedral elements $\Omega^e$. Within each element, the solution is represented by a Lagrange polynomial of degree $N$, 
defined on Gauss-Legendre interpolation nodes. This results in a piecewise polynomial solution representation in the computational domain $\Omega$. 
Adjacent elements are coupled at their faces by classical single-phase approximate Riemann solvers like the HLLC solver \cite{Toro1994}. 
Even though the DGSEM scheme provides efficient and high-order accurate results in smooth flow regions, it suffers from Gibbs oscillations in the 
presence of shocks, phase interfaces or severe under-resolution. Here, a robust Finite Volume (FV) scheme is applied on an h-refined sub-cell 
grid of $N_{FV}^d$ FV sub-cells per DG element to provide reliable shock capturing and precise interface localization. This involves a conservative 
transformation of the element-local solution between a polynomial and a piecewise constant solution representation. 
The accuracy of the FV method is enhanced by a second-order TVD reconstruction scheme. 

In this work, we use the hp-adaptive extension for the DGSEM by Mossier et al. \cite{Mossier2022,Mossier2023} that allows for a p-adaptive 
discretization with variable element-local degree $N$ and a variable FV sub-cell resolution $N_{FV}$. The adaptive scheme allows for high local 
accuracy at the phase interface while maintaining a coarser resolution in most regions of the computational domain. Both p-adaptivity and FV sub-cell 
shock capturing are controlled by an indicator, which infers the smoothness and an error estimate from the analysis of the modal spectrum 
of the element local solution polynomials \cite{Mavriplis1990,Mossier2022}. 
The resulting scheme is advanced in time with an explicit fourth-order low-storage Runge-Kutta (RK) method \cite{Kennedy2003}.
 
\subsection{The Level-Set Ghost-Fluid Method}
\label{subsec:InterfaceTracking}
The level-set Ghost-Fluid method relies on the sign of the level-set function to distinguish the bulk phases and infers the interface position from
the zero iso-contour of the level-set field. In the present framework, the level-set transport equation \eqref{eq:levelset} is discretized with a 
path-conservative DGSEM \cite{castro2006high,Dumbser2016a} operator with FV sub-cell limiting \cite{Joens2021}. Analogously to Section \ref{subsec:BulkSolver}, 
the solution is therefore represented either by piece-wise DG polynomials or a piecewise constant FV discretization on an h-refined sub-cell grid. 
In smooth regions, the DGSEM operator is applied for its high-order accuracy and efficiency. If changes in the topology of the phase interface 
occur, e.g. during the merging of droplets, the level-set may exhibit discontinuities and the robust FV sub-cell scheme is applied. 
In the present work, we employ an hp-adaptive extension \cite{Mossier2023}, that allows for variable polynomial degrees $N$ and variable FV sub-cell 
resolutions $N_{FV}$ within each element. The level-set transport equation 
is advanced in time with an explicit RK scheme, analogously to the bulk fluid operator. The remaining building blocks for the level-set 
method, the velocity extrapolation and level-set reinitialization are implemented by solving Hamilton-Jacobi-type equations with a fifth-order 
WENO scheme \cite{Peng1998,FechterPhd}. To reduce computational costs, the interface tracking is restricted to a narrow band of two to three 
elements around the interface \cite{Adalsteinsson1995}. 

With the interface tracking established, the remaining step is the coupling of the bulk fluid phases. The present scheme employs the 
Ghost-Fluid method of Fedkiw et al. \cite{Fedkiw1999}. It relies on the definition of ghost states at the phase interface 
to provide boundary conditions and achieve a sound coupling between the bulk phases. Following Merkle et al. \cite{Merkle2007}, 
we solve an interfacial Riemann problem and assign the intermediate Riemann states as ghost states. The interfacial Riemann problem further
provides the local advection velocity of the phase boundary. Solution strategies for the two-Phase Riemann problem are covered in 
depth in Section \ref{subsec:TwoPhaseRP}. 

When a cell changes its phase during the convection of the phase boundary, the cell needs to be initialized with a physically sound state.
In the present Ghost-Fluid scheme, cells that underwent a phase change are populated with the intermediate states of the corresponding
interfacial Riemann problem. For an analysis of the effect of this non-conservative procedure, the reader is referred to the publication 
of J\"ons et al. \cite{Joens2023}.

\subsection{Semi-Analytical Source Term Integration}
\label{subsec:SourceTermTreatment}
Solving the GPR equation system \eqref{eq:HPR_equations} for inviscid, heat-conducting fluids 
requires the discretization of possibly stiff algebraic source terms in the thermal impulse equations. Therefore, a simple 
explicit source term integration is impracticable due to a prohibitive time step limitation in the stiff relaxation regime.  
A common strategy to solve an inhomogeneous partial differential equation (PDE) is the splitting approach. 
This involves the solution of the homogeneous PDE in the first step 
and a subsequent correction for the contribution of the source term by solving an ordinary differential equation (ODE).
In the present case, this requires advancing the PDE \eqref{eq:HPR_matrix_vector} from a time level $t^n$ towards $t^{n+1}$ by first solving 
the homogeneous system 
\begin{equation}
\frac{\partial\QB}{\partial t} + \nabla_{\xB} \cdot \FB(\QB)=0 \quad \text{in}  \quad \Omega\times[t^n,t^{n+1}] \quad \text{with}  \quad \QB(\xB,t^n)=\QB^n.
\label{eq:HPR_PDE_homogen}
\end{equation}
The solution of the homogeneous system at $t^{n+1}$ is denoted $\QB^*$.
According to the splitting approach, the solution at $\QB^{n+1}$ is then obtained by solving the ODE
\begin{equation}
\frac{\partial\QB}{\partial t} = S(\QB) \quad \text{in}  \quad \Omega\times[t^n,t^{n+1}] \quad \text{, with}  \quad \QB(\xB,t^n)=\QB^*.
\label{eq:HPR_ODE}
\end{equation}
However, in recent studies \cite{Boscheri2022,Chiocchetti2023} it was demonstrated that a simple splitting
approach fails to recover the non-trivial equilibrium of the thermal relaxation formulation in the stiff relaxation limit. 
The authors therefore highlighted the requirement of an asymptotic preserving scheme, that ensures convergence towards 
the Euler-Fourier system in the stiff regime. They therefore suggested solving a modified 
ODE for the thermal impulse vector $\JB=\rho \jB$, that accounts for the discrete update of the left-hand side of equation 
\eqref{eq:HPR_matrix_vector}
\begin{equation}
\frac{\partial\JB}{\partial t} = \PB^* + S(\JB) \quad \text{in}  \quad \Omega\times[t^n,t^{n+1}] \quad \text{with}  \quad \JB(\xB,t^n)=\JB^*,
\label{eq:HPR_ODE_Mod}
\end{equation}
with the additional term $\PB^*=(\JB*-\JB^n)/\Delta t$. The modified ODE is shown to recover the Fourier law in the stiff limit and is thus 
called an asymptotic preserving scheme. Following \cite{Chiocchetti2023}, an exact solution to the ODE \eqref{eq:HPR_ODE_Mod} can be found as
\begin{equation}
\JB^{n+1}=(\JB^n-\tau_H \PB^*)\exp(-\Delta t/\tau_H)+\tau_H \PB^* \quad \text{, with} \quad \tau_H=\frac{\Theta(\tau)}{\alpha^2}=\frac{\rho\lambda}{T \alpha^2},
\label{eq:HPR_ODE_exac}
\end{equation}
allowing for a computationally efficient evaluation. This semi-analytical scheme is used throughout the present paper for its robustness 
in the stiff relaxation limit and its low computational cost. It is applied for the bulk discretization in combination with a fourth-order 
Runge-Kutta method, where the semi-analytical solver simply replaces each Euler step found in the standard RK scheme. 
Thus, the ODE \eqref{eq:HPR_ODE_Mod} is solved between subsequent RK time stages instead of the time interval $t\in[t^n,t^{n+1}]$.

Note that such a semi-analytical approach, while specifically designed for the relaxation systems, namely heat impulse and strain relaxation found in the GPR model, 
is based on a rather general principle and its early origins trace back to the solution of Baer-Nunziato
relaxation sources \cite{Chiocchetti2020}, meaning the coupled system of mechanical friction,
pressure relaxation, and compaction dynamics.

In particular, it aims solely at being an accurate and efficient integrator and does not take advantage of
the symmetric hyperbolic thermodynamically compatible (SHTC) structure of the governing equations.
See e.g. \cite{Abgrall2023, Boscheri2023} for recent numerical schemes specifically designed to mimic
the SHTC structure at the discrete level, in a general GPR context and \cite{Thomann2023} in multiphase
flows.
 
\subsection{Approximate Two-Phase Riemann Solvers}
\label{subsec:TwoPhaseRP}
At a sharp phase interface $\Gamma$, separating a liquid and vapor bulk phase, a two-phase Riemann problem can be defined. 
Using the rotational invariance of system \eqref{eq:HPR_equations}, we consider the initial value problem 
\begin{equation}
\QB(\xB,t=0)=
\begin{cases*}
\QB_l \quad \text{for} \; x\le x_0 \; \text{(liq)},\\
\QB_v \quad \text{for} \; x > x_0 \; \text{(vap)},
\end{cases*}
\end{equation}
with a liquid state $\QB_l$ on the left and a vapor state $\QB_v$ on the right in a reference space, normal to the phase interface $\Gamma$. 
A key challenge when solving the two-phase Riemann problem with phase transition is the dissipative nature of heat conduction that causes a loss 
of self-similarity in the solution. Previous publications saw a variety of strategies to circumvent this issue. In \cite{Hantke2013,Rohde2018} 
the isothermal Euler equations were considered, thus assuming an instantaneous heat transfer. Hitz et al. \cite{Hitz2020,Hitz2021} and J\"ons et al. 
\cite{Joens2023} formulated interfacial Riemann solvers for the Euler-Fourier system that neglect heat conduction across all waves except the phase 
boundary, where a constant heat flux from a local phase transition model is imposed. In a recent study of  M\"uller et al. \cite{Mueller2023},
an approximate two-phase Riemann solver was proposed for the GPR system. It is formulated for the homogeneous part of the PDE \eqref{eq:HPR_equations},
thus neglecting the dissipative effects of heat transfer. While their approach employs a phase transition model to predict the interfacial mass flux,
it failed to enforce an interfacial heat flux and thus the interfacial entropy production.

The goal of this paper is the formulation of a thermodynamically consistent interfacial Riemann solver for the GPR system, that
finds an entropy solution in agreement with a local phase transition model. In the following, the construction of two approximate 
Riemann solvers denoted as $HLLP_{mq}$ and $HLLP_{m}$ is motivated. 

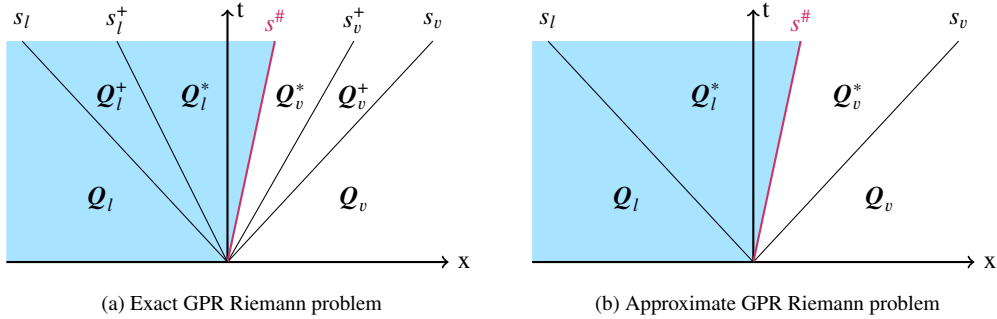
\begin{figure}[htb]
	\begin{subfigure}[t]{0.49\textwidth}
		\centering
		\begin{tikzpicture}
		[
		scale=0.21,
		grid/.style={black},
		every node/.style={scale=0.8}
		]
		
		\def \x{0}
		\def \y{0}					
		
		\coordinate (A) at (-14,0);
		\coordinate (B) at (0,0);
		\coordinate (C) at (3,14);
		\coordinate (D) at (-14,14);
		\fill[fill=lightblue!50] (A) -- (B) -- (C) -- (D) -- cycle;

		\draw[thick,->] (\x-14,\y) -- (\x+14,\y) node[anchor=west] {x};
		\draw[thick,->] (\x,\y)    -- (\x,\y+16) node[anchor=west] {t};
		
		\draw[] (\x,\y) -- (\x+13,\y+14);
		\draw[] (\x,\y) -- (\x+8,\y+14);
		
		\draw[purple!80,thick] (\x,\y) -- (\x+3,\y+14);
		
		\draw[] (\x,\y) -- (\x-13,\y+14);
		\draw[] (\x,\y) -- (\x-7,\y+14);

		\draw[]  (\x+8,\y+4) node [shape=rectangle,rounded corners,text centered] %
		{$\QB_{v}$};
		
		\draw[]  (\x-8,\y+4) node [shape=rectangle,rounded corners,text centered] %
		{$\QB_{l}$};

		\draw[]  (\x+4,\y+10.5) node [shape=rectangle,rounded corners,text centered] %
		{$\QB_{v}^*$};
		\draw[]  (\x+8,\y+10.5) node [shape=rectangle,rounded corners,text centered] %
		{$\QB_{v}^+$};
		
		\draw[]  (\x-2,\y+10.5) node [shape=rectangle,rounded corners,text centered] %
		{$\QB_{l}^*$};
		\draw[]  (\x-7.3,\y+10.5) node [shape=rectangle,rounded corners,text centered] %
		{$\QB_{l}^+$};	
		
		\draw[] (\x-13,\y+15.3) node [shape=rectangle,rounded corners,text centered] %
		{$s_{l}$};
		\draw[] (\x-7,\y+15.3) node [shape=rectangle,rounded corners,text centered] %
		{$s_{l}^+$};

		\draw[purple!80] (\x+3,\y+15.3) node [shape=rectangle,rounded corners,text centered] %
		{$\sB$};
		
		\draw[] (\x+13,\y+15.3) node [shape=rectangle,rounded corners,text centered] %
		{$s_{v}$};
		\draw[] (\x+8,\y+15.3) node [shape=rectangle,rounded corners,text centered] %
		{$s_{v}^+$};
		
		\end{tikzpicture}
		\caption{Exact GPR Riemann problem}
		\label{fig:hpr_exact}
	\end{subfigure}
	\hfill
	\begin{subfigure}[t]{0.49\textwidth}
		\centering
		\begin{tikzpicture}
		[
		scale=0.21,
		grid/.style={black},
		every node/.style={scale=0.8}
		]
		
		\def \x{0}
		\def \y{0}					
		
		\coordinate (A) at (-14,0);
		\coordinate (B) at (0,0);
		\coordinate (C) at (3,14);
		\coordinate (D) at (-14,14);
		\fill[fill=lightblue!50] (A) -- (B) -- (C) -- (D) -- cycle;
		
				\draw[thick,->] (\x-14,\y) -- (\x+14,\y) node[anchor=west] {x};
		\draw[thick,->] (\x,\y)    -- (\x,\y+16) node[anchor=west] {t};
		
		\draw[] (\x,\y) -- (\x+13,\y+14);

		\draw[purple!80,thick] (\x,\y) -- (\x+3,\y+14);
		
		\draw[] (\x,\y) -- (\x-13,\y+14);
		\draw[]  (\x+8,\y+4) node [shape=rectangle,rounded corners,text centered] %
		{$\QB_{v}$};
		
		\draw[]  (\x-8,\y+4) node [shape=rectangle,rounded corners,text centered] %
		{$\QB_{l}$};
		
		\draw[]  (\x+6,\y+10.5) node [shape=rectangle,rounded corners,text centered] %
		{$\QB_{v}^*$};
		
		\draw[]  (\x-3,\y+10.5) node [shape=rectangle,rounded corners,text centered] %
		{$\QB_{l}^*$};

		\draw[] (\x-13,\y+15.3) node [shape=rectangle,rounded corners,text centered] %
		{$s_{l}$};

		\draw[purple!80] (\x+3,\y+15.3) node [shape=rectangle,rounded corners,text centered] %
		{$\sB$};

		\draw[] (\x+13,\y+15.3) node [shape=rectangle,rounded corners,text centered] %
		{$s_{v}$};
		
		\end{tikzpicture}
		\caption{Approximate GPR Riemann problem}
		\label{fig:hpr_approx}
	\end{subfigure}
	\caption{Exact and approximate wave pattern for a GPR two-phase Riemann problem with phase transition. 
	The non-classical wave of the phase interface, associated with the velocity of the phase boundary $\sB$, 
	is highlighted in pink. The remaining waves are shock or rarefaction waves. }
	\label{fig:riemann_wave_pattern}
\end{figure}

Consider the wave pattern associated with the exact solution of the homogeneous GPR equation system in figure \ref{fig:hpr_exact}. It consists
of four shock or rarefaction-like waves that are related to the eigenvalues of the system. The material interface at the phase boundary is 
represented by an additional non-classical, undercompressive shock wave. As discussed by Menikoff and Plohr \cite{Menikoff1989}, it is not 
associated with an eigenvalue of the system in the presence of phase transition. This is explained by the non-convex spinodal region that 
separates the phases and causes a local breakdown of hyperbolicity due to imaginary eigenvalues. As outlined in Section \ref{subsec:evaporation_model},
we evaluate a kinetic relation \cite{Abeyaratne1991} to control the entropy production at the phase interface to unitize the undercompressive shock. The interfacial 
entropy production in turn is predicted by a phase transition model. While exact two-phase Riemann solvers with phase transition 
were previously reported for the Euler-Fourier system \cite{Hitz2021,Joens2023}, the authors found approximate solution strategies to be
of comparable accuracy, while being significantly more robust and computationally efficient. 
Therefore, we focus on the construction of an approximate two-phase Riemann solver for the GPR system, that assumes a simplified
HLLC-like wave pattern, as depicted in figure \ref{fig:hpr_approx}. 
The reduced wave fan consists of two outer classical waves and a central undercompressive wave, which represents the phase boundary. 
By solving the approximate two-phase Riemann problem, the inner states $\QB^*_l$ and $\QB^*_v$ and the interface velocity $\sB$ are determined. 

\subsubsection{The $\text{HLLP}_{mq}$ Two-Phase Riemann Solver} 
\label{subsec:hllp_mq}
The $\text{HLLP}_{mq}$ solver consists of two main building blocks: a non-linear equation system $\smacro$ 
for the intermediate states, defined by jump conditions across the three waves and a thermodynamic closure
relation $\smicro$, based on a local phase transition model.
The equation system $\smacro$ is constructed from the interfacial jump conditions \eqref{eq:interface_jump}
and Rankine-Hugoniot conditions for the two outer waves. It relates the inner states $\QB^*_l$ and $\QB^*_v$ and the 
velocity of the phase boundary $\sB$ to the initial states $\QB_l$ and $\QB_v$: 
\begin{equation}
\left(\QB^*_l,\QB^*_v,\sB \right)^T := 
\smacro \left(\QB_l,\QB_v,\dot{m}^*,q^*_v \right).
\end{equation}
When phase transition is considered, estimates for the interfacial mass $\dot{m}^*$ and heat flux $q^*_v$ have to be supplied through
a closure relation to find a unique and thermodynamically consistent entropy solution. Throughout this work, 
we use the evaporation model of Cipolla et al. \cite{Cipolla1974} as outlined in Section \ref{subsec:evaporation_model}.
Given the inner states $\QB^*_l$ and $\QB^*_v$, the phase transition model $\smicro$ provides estimates for $\dot{m}^*$ and $q^*_v$ based on
the phenomenological force flux relations \eqref{eq:onsager_relation}:
\begin{equation}
	\left(\dot{m}^*,q^*_v \right)^T := 
	\smicro \left(\QB^*_l,\QB^*_v,\sB \right).
\end{equation}
In conjunction, the integral jump relations $\smacro$ and the thermodynamic closure relation $\smicro$ form a closed 
non-linear system of equations. Starting with an initial guess $(\dot{m}^*,q^*_v)^T=(0,0)^T$, the system can be solved iteratively in $(\dot{m}^*,q^*_v)^T$
with a Newton algorithm:
\begin{equation}
\left(\chi^{\Gamma},\dot{m}^*,q^*_v \right)^T -  
\smicro \left(\smacro \left(\QB_l,\QB_v,\dot{m}^*,q^*_v \right)\right)\overset{!}{=}0.
\end{equation}
In the following, we derive the non-linear equation system $\smacro$ for the GPR model. First,
the wave speeds of the outer waves $s_l$ and $s_v$ are approximated by two-phase adapted HLL estimates of 
Davis et al. \cite{Davis1988} as     
\begin{align}
	s_l&=u_l-c_{s,l}, \\
    s_v&=u_v+c_{s,v},
\end{align}
with $c_{s,l}$ and $a_{s,v}$ denoting the sound speed in the liquid and vapor respectively. 
Next, Ranine-Hugoniot conditions are defined for the two outer waves $s_l$ and $s_v$
\begin{subequations}
	\begin{align}
	\llbracket \dot{m} \rrbracket_{l,v} &= 0,\\
	\dot{m}\llbracket u \rrbracket_{l,v}+\llbracket p \rrbracket_{l,v}&=0,\\
	\dot{m}\llbracket e \rrbracket_{l,v}+\llbracket up \rrbracket_{l,v}+\llbracket q \rrbracket_{l,v}&=0,\\
	\dot{m}\llbracket j \rrbracket_{l,v}+\llbracket T \rrbracket_{l,v}&=0,
	\end{align}
	\label{eq:iouter_jump}
\end{subequations}
with the notation $\llbracket z \rrbracket_{l}=z_l^*-z_l$ and $\llbracket z \rrbracket_{v}=z_v^*-z_v$.
Given the jump relations across the phase boundary \eqref{eq:iouter_jump}, previously defined in \ref{subsec:evaporation_model},
a set of balance equations for the mass 
\begin{align}
\dot{m}_l=\rho_l(u_l-s_l)&=\rho_l^*(u_l^*-s_l), \\
\rho_l^*(u_l^*-\sB)&=\rho_v^*(u_v^*-\sB),       \\
\dot{m}_v=\rho_v(u_v-s_v)&=\rho_v^*(u_v^*-s_v), 
\end{align}
and impulse  
\begin{align}
\dot{m}_l u_l+p_l&=\dot{m}_l u_l^*+p_l^*, \\
\dot{m}^*u_v^*+p_v^*&=\dot{m}^* u_l^*+p_l^*+\Delta p_{\sigma}, \\
\dot{m}_v u_v+p_v&=\dot{m}_v u_v^*+p_v^*
\end{align}
can be defined across the three waves. In total, we obtain six equations for the density, 
velocity and pressure $\rho^*_l,\rho^*_v,u^*_l,u^*_v,p^*_l,p^*_v$ of the inner state vectors $\QB^*_l$ and $\QB^*_v$. 
Further, the unknown velocity of the phase boundary $\sB$ can be related to the interfacial mass flux estimate $\dot{m}^*$ by 
\begin{equation}
\dot{m}^*= \rho_l^*(u_l^*-\sB)=\rho_v^*(u_v^*-\sB).
\end{equation}
The resulting linear equation system is called the \textit{mechanical system} and can be solved algebraically. 

To find expressions for the remaining conservative quantities $(\rho e)^*_l,(\rho e)^*_v$ and $(\rho j)^*_l,(\rho j)^*_v$, 
we consider the jump relations for the energy   
\begin{align}
\dot{m}_l e_l+u_l p_l+q_l&=\dot{m}_l e_l^*+u_l^* p_l^*+q_l^*, \label{eq:e_l} \\
\dot{m}_v e_v^*+u_v^* p_v^*+q_v^*&=\dot{m}_l e_l^*+u_l^* p_l^*+q_l^*+\sB \Delta p_{\sigma}, \label{eq:e_lv} \\
\dot{m}_v e_v+u_v p_v+q_v&=\dot{m}_v e_v^*+u_v^* p_v^*+q_v^*, \label{eq:e_v}
\end{align}
and the thermal impulse
\begin{align}
\dot{m}_l j_l+T_l&=\dot{m}_l j_l^*+T_l^*, \label{eq:j_l}\\
\dot{m}_v j_v^*+T_v^*&=\dot{m}_l j_l^*+T_l^*+\Delta T, \label{eq:j_lv} \\
\dot{m}_v j_v+T_v&=\dot{m}_v j_v^*+T_v^*. \label{eq:j_v}
\end{align}
With constitutive relations, linking the heat flux to the temperature and thermal impulse
\begin{align}
q_l^*&=\alpha^2T_l^*j_l^*, \label{eq:q_l}\\
q_v^*&=\alpha^2T_v^*j_v^*, \label{eq:q_v}
\end{align}
we obtain a total of eight equations, that we call the \textit{thermodynamic system}.
In the present work, we solve the \textit{thermodynamic system} in two steps. First, we use the jump condition \eqref{eq:j_l}-\eqref{eq:j_v}
for the energy together with an estimate for the heat flux $q_v^*$ from the phase transition model $\smicro$ to compute the unknown energies 
$e^*_l$ and $e^*_v$ and the heat flux $q_l^*$. The resulting linear equation system can be easily solved algebraically. 
To determine the remaining quantities $j^*_l$, $j^*_v$, $T_l^*$ and $T_l^*$, we insert equations \eqref{eq:q_l} and \eqref{eq:q_v} in \eqref{eq:j_l} 
and \eqref{eq:j_v}:
\begin{align}
	\dot{m}_l j_l+T_l=\dot{m}_l \frac{q_l^*}{\alpha^2 T_l^*}+T_l^*, \\
	\dot{m}_v j_l+T_v=\dot{m}_v \frac{q_l^*}{\alpha^2 T_l^*}+T_v^*.
\end{align}
This produced two quadratic equations in $T^*_l$ and $T^*_v$ that can be solved as follows:
\begin{align}
T_{l_{1,2}}^*&=\frac{1}{2}\left(\dot{m}_l j_l+T_l\pm\sqrt{(\dot{m}_l j_l+T_l)^2-4  \frac{\dot{m}_l q_l^*}{\alpha^2 T_l^*}}\right), \\
T_{v_{1,2}}^*&=\frac{1}{2}\left(\dot{m}_v j_v+T_v\pm\sqrt{(\dot{m}_v j_v+T_v)^2-4  \frac{\dot{m}_v q_v^*}{\alpha^2 T_v^*}}\right).
\end{align}
Even though the quadratic equations yield two possible solutions, we found from numerical experiments, 
that only $T_{l_{1}}^*$ and $T_{v_{1}}^*$ are physically meaningful solutions. 

\begin{figure}[t]
	\centering
	\includegraphics[width=0.70\textwidth]{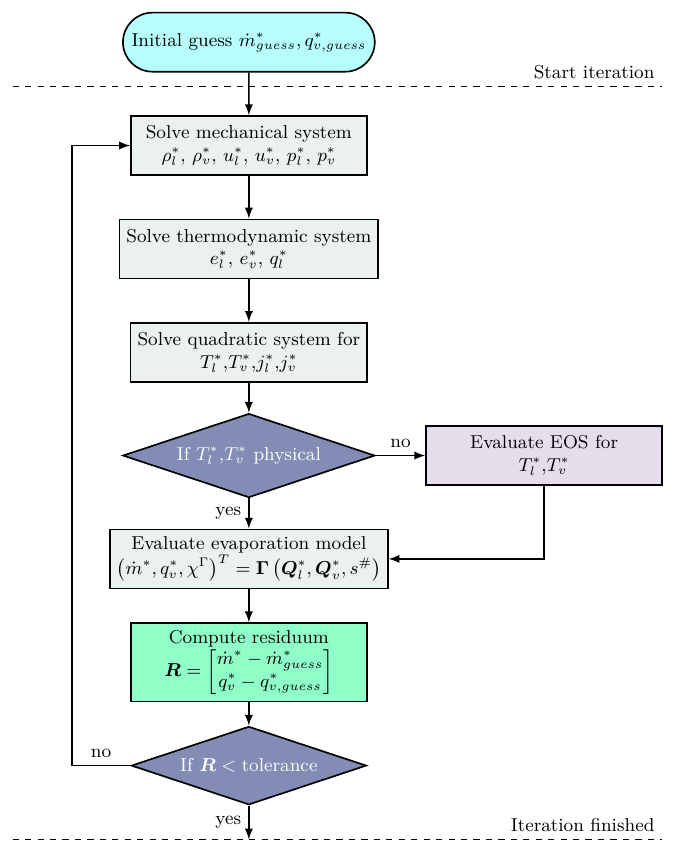}
	\caption{Flowchart, illustrating the iterative solution procedure of the $\text{HLLP}_{mq}$ Riemann solver.}
	\label{fig:hllp_mq}
\end{figure}

An inherent issue of the quadratic
equation system is the possibility of a negative discriminant during the iterative solution procedure
when given inaccurate initial guesses for $\dot{m}^*$ and $\dot{q}^*$. Such behavior was encountered 
when computing more complex multi-dimensional setups involving severe phase interface deformations, surface tension and 
strong thermodynamic non-equilibrium. In the $\text{HLLP}_{mq}$ Riemann solver, we circumvent this issue with 
a fallback to an EOS evaluation for the temperatures, in case the discriminant exhibits a negative sign. The non-equilibrium contribution 
$\frac{1}{2}\alpha j^2$ to the total energy \eqref{eq:totalenergy} is neglected during the EOS evaluation, since $j_l^*$ and $j_v^*$ are 
still unknown. We want to emphasize, that this fallback is only encountered during the initial iteration steps and not allowed in the final step. 
With the temperatures $T_l^*$ and $T_v^*$ known, we can finally obtain the thermal impulses from equations \eqref{eq:j_l} 
and \eqref{eq:j_v}. The main steps during the iterative solution process of the $\text{HLLP}_{mq}$ solver are visualized in figure \ref{fig:hllp_mq}.
The iterative solution procedure is illustrated in figure \ref{fig:hllp_mq}.

Notice, that equation \eqref{eq:j_lv} was not used in the solution procedure so far. To explain this curious decision, we 
want to recapitulate the effects of the relaxation source terms in the GPR model. In the balance equation
for the thermal impulse \eqref{eq:HPR_equations_4}, the thermal relaxation source term controls the dissipation and therefore
the entropy production due to heat conduction. Since the $\text{HLLP}_{mq}$ two-phase Riemann solver is formulated for the homogeneous GPR system,
it would appear that we enforced a negligible entropy production across the interface. However, even if the thermal relaxation source terms had 
been included, the GPR model could not be expected to predict the correct entropy production in the presence of phase transition 
due to the breakdown of the continuum assumption across the phase boundary.

By incorporating a temperature jump $\Delta T$ in the jump relation \eqref{eq:interface_jump_4}, 
we allow the phase transition model to impose an interfacial entropy production.
In that sense the additional degree of freedom $\Delta T$ in the thermal impulse jump relation acknowledges the existence of an 
unknown thermal relaxation source term, that is determined through the local thermodynamic model. The proposed strategy thus 
allows a consistent coupling between the GPR model and a local thermodynamic phase transition model, that guides the 
iterative solution procedure to the correct entropy solution. 

\subsubsection{The $\text{HLLP}_{m}$ Two-Phase Riemann Solver} 
\label{subsec:hllp_m}
While the proposed $\text{HLLP}_{mq}$ interface Riemann solver combines the GPR model consistently with a thermodynamic closure relation, 
the iteration in the two variables poses a possible source for instabilities and increased computational costs. This motivates,
the construction of a further simplified two-phase Riemann solver, called $\text{HLLP}_{m}$, which requires only an iteration in
the mass flux $\dot{m}^*$. Since the linear \textit{mechanical system} is evaluated independent of the heat flux $q_v^*$, the differences between 
the $\text{HLLP}_{mq}$ and $\text{HLLP}_{m}$ solvers are restricted to the construction of the \textit{thermodynamic system}. 

\begin{figure}[t]
	\centering
	\includegraphics[width=0.7\textwidth]{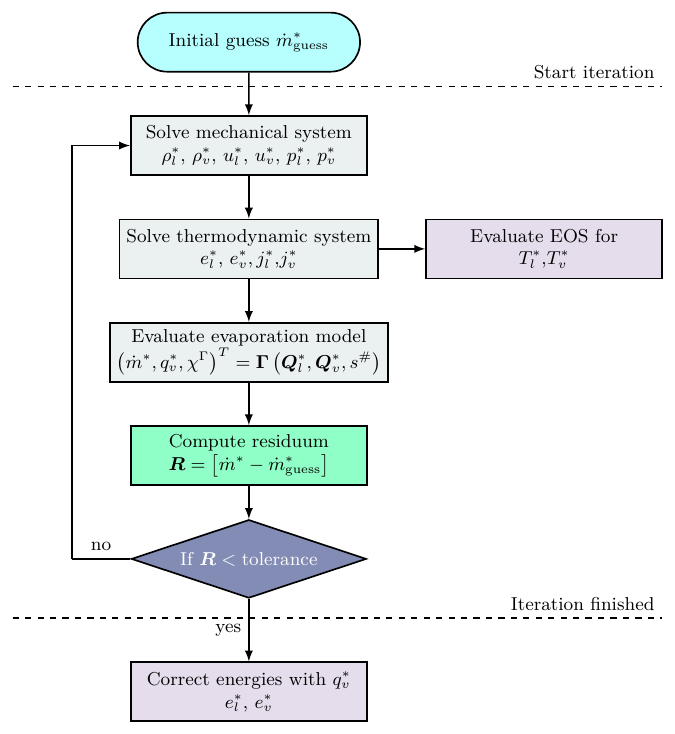}
	\caption{Flowchart, illustrating the iterative solution procedures of the $\text{HLLP}_{m}$ Riemann solver.}
	\label{fig:hllp_m}
\end{figure}

The iteration in $q_v^*$ is necessary since the jump relations for the energy across the three waves provide only three equations 
for the four unknown energies and heat fluxes $e_l^*$, $e_v^*$, $q_l^*$  and $q_v^*$. 
A possible simplification to obtain closure without a prediction of $q_v^*$  is to neglect the heat fluxes in the energy jump 
conditions across the outer waves. Equations \eqref{eq:e_l} and \eqref{eq:e_v} thus reduce to
\begin{align}
\dot{m}_l e_l+u_l p_l&=\dot{m}_l e_l^*+u_l^* p_l^* \label{eq:e_l_no_q}, \\
\dot{m}_v e_v+u_v p_v&=\dot{m}_v e_v^*+u_v^* p_v^* \label{eq:e_v_no_q}. 
\end{align}
Consequently, they can be solved for the unknown energies $e_l^*$ and $e_v^*$ independently of the phase transition model. 
Subsequently, the temperatures $T_l^*$ and $T_v^*$ are determined by the EOS, while the non-equilibrium contributions $\frac{1}{2}\alpha j^2$ 
to the total energy \eqref{eq:totalenergy} are neglected. Finally, with the temperatures known, the thermal impulses $j_l^*$ and $j_v^*$ are 
computed with equations \eqref{eq:j_l} and \eqref{eq:j_v}. As depicted in figure \ref{fig:hllp_m}, the simplified \textit{thermodynamic system} 
allows to solve the two-phase Riemann problem with an iteration over the mass flux $\dot{m}^*$ only.

To account for the heat flux $q_v^*$, predicted by the evaporation model, a correction step is performed after the iteration in $\dot{m}^*$ 
has converged. Using the states $\QB^*_l$ and $\QB^*_v$ from the final iteration and the matching heat flux prediction $q_v^*$ from the 
evaporation model, the heat flux $q_l^*$ is computed via equation \eqref{eq:e_lv}:
\begin{align}
q_l^* = \dot{m}_v e_v^*+u_v^* p_v^*+q_v^*-(\dot{m}_l e_l^*+u_l^* p_l^*+\sB \Delta p_{\sigma}).\label{eq:q_star_l}
\end{align}
When the heat fluxes $q_v^*$ and $q_l^*$ are available, the energies $e_l^*$, $e_v^*$ can finally be computed with \eqref{eq:e_l} and \eqref{eq:e_v}, 
without neglecting the heat conduction across the outer waves, as done in \eqref{eq:e_l_no_q} and \eqref{eq:e_v_no_q}. This correction step 
for the energies is crucial to achieve a solution, that is in agreement with the predicted thermodynamic flux $q_v^*$ of the evaporation model.
The final solution procedure for the $\text{HLLP}_{mq}$ solver is outlined in figure \ref{fig:hllp_m}.

\section{Numerical Results}
\label{sec:Results}
In this Section, we apply the proposed numerical method to a range of representative test cases. First, the thermal relaxation formulation
of the GPR continuum model is compared against the Euler-Fourier system. Therefore, a one-dimensional heat conduction problem and the
well-known Rayleigh-B\'{e}nard convection are studied. Next, the novel $\text{HLLP}_{mq}$ and $\text{HLLP}_{m}$ Riemann solvers
are validated with evaporating shock tube computations against MD reference data and solutions obtained with the Euler-Fourier approach. 
The investigation focuses on a qualitative analysis, convergence studies and a comparison of the computational costs.  
Finally, we apply the method to shock-droplet interactions that involve phase transition, surface tension and 
complex deformations of the phase interface. 

\subsection{One-Dimensional Heat Conduction}
\label{subsec:1D_heat_conduction}
\begin{wrapfigure}{htb!}{0.45\textwidth}
	\vspace{-5pt}
	\begin{center}
		\begin{tikzpicture}[scale=1.0,every node/.style={scale=1}]
			\draw[thick,fill=lightgrey] (0,0) rectangle (4,1);
			\draw[-latex] (-0.2,0.3) -- node[left] {$q^+$} (0.2,0.7);
			\draw[-latex] (3.8,0.3)  -- node[right]{$q^-$} (4.2,0.7);
			\draw[text=black,thick] (-0.6,1.2) node {$T_B^+=3$};
			\draw[text=black,thick] (4.5 ,1.2) node {$T_B^-=1$};
			\draw[-latex] (-0.6,-0.3) -- node[below]{$x$} (-0.2,-0.3);
			\end{tikzpicture}	
	\end{center}
	\vspace{-15pt}
	\caption{Computational setup for the one-dimensional heat conduction test case. Heat fluxes are imposed on the left and right boundaries.}
	\vspace{0pt}
	\label{fig:1D_setup}
\end{wrapfigure}
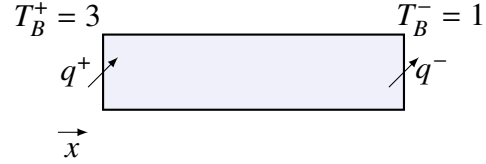

\begin{figure}[htb]
	\begin{subfigure}[t]{0.325\textwidth}
		\centering
		\includegraphics[width=1.\textwidth]{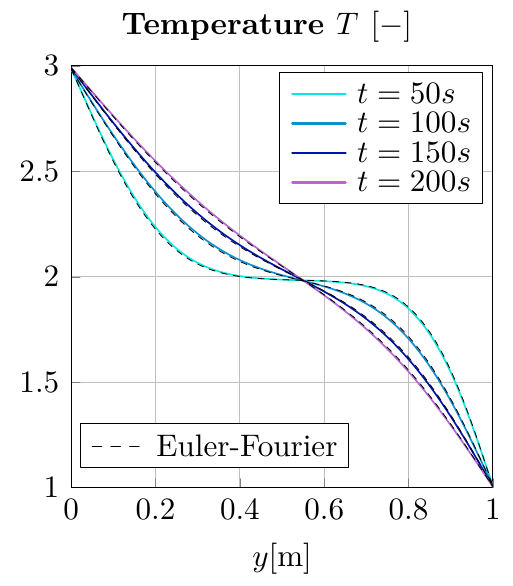}
		\caption{$\lambda=1\cdot 10^{-3}$}
	\end{subfigure}
	\hfill
	\begin{subfigure}[t]{0.325\textwidth}
		\centering
		\includegraphics[width=1.\textwidth]{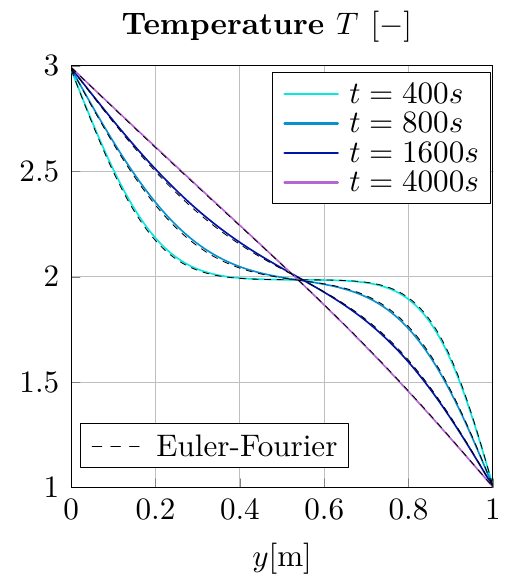}
		\caption{$\lambda=1\cdot 10^{-4}$}
	\end{subfigure}
	\hfill
	\begin{subfigure}[t]{0.325\textwidth}
		\centering
		\includegraphics[width=1.\textwidth]{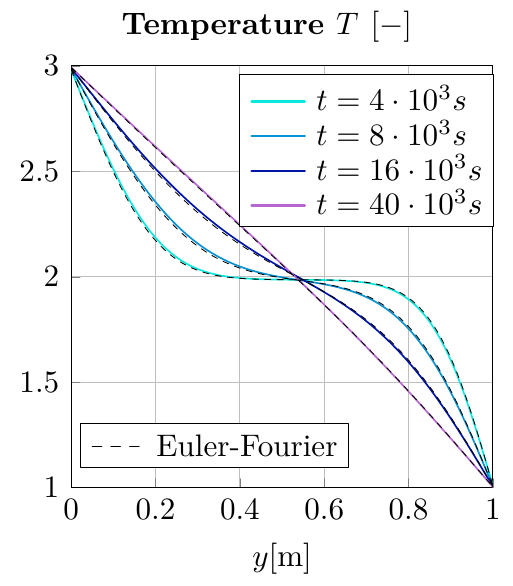}
		\caption{$\lambda=1\cdot 10^{-5}$}
	\end{subfigure}
	\caption{Temperature profile of the one-dimensional heat conduction test case at different time instances. 
	The heat conductivity is increased from left to right and the GPR model is compared against a Euler-Fourier reference solution.}
	\label{fig:1D_Heat}
\end{figure}

In this paragraph, we investigate a one-dimensional heat conductivity dominated singe-phase flow to validate the thermal relaxation model
of the GPR system against the Euler-Fourier model. We consider a computational domain $\Omega=[0,1]$ that contains a resting perfect 
gas with a constant temperature $T=2.0$ and pressure $p=2.5$ at the initial time $t=0.0$. Heat capacities at constant volume and pressure 
are chosen as $c_v=0.718$ and $c_p=1.005$ respectively and the relaxation time $\tau$ of the GPR model is defined according to \eqref{eq:tau_kin}. 
The setup is visualized in figure \ref{fig:1D_setup} with periodic boundary conditions in y-direction and heat fluxes $q^+$ and $q^-$ imposed 
at the lower and upper boundaries in x-direction. With the heat transmission coefficient $\alpha=100\lambda$, the heat fluxes are defined as  
\begin{align*}
	q^+=\rho\alpha(T-T_B^+), \\
	q^-=\rho\alpha(T-T_B^-), 
\end{align*}
and depend on the density and temperature at the wall. The left wall is heated to a temperature $T_B^+=3$ and the right wall is cooled 
to a temperature $T_B^-=1$. Wherever not indicated, standard SI units can be assumed. 
The domain $\Omega$ is discretized with $[256]\times[1]$ DG elements of degree $N=2$. We perform three computations  
with different heat conductivities $\lambda=1\cdot 10^{-3}$, $\lambda=1\cdot 10^{-4}$ and $\lambda=1\cdot 10^{-5}$ until the final computation times $t=200$, 
$t=4000$ and $t=40000$ respectively. 
Figure \ref{fig:1D_Heat} depicts the temperature profiles in x-direction at different time instances, 
computed with the GPR model and the Euler-Fourier system as a reference solution. A near-perfect agreement of both solutions can be observed 
for the considered range of thermal conductivities.

\subsection{Rayleigh-B\'{e}nard Convection}
\label{subsec:Rayleigh_Bennard}
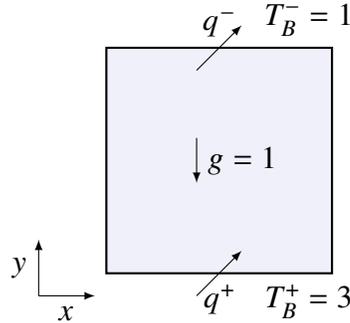
\begin{wrapfigure}{r}{0.45\textwidth}
	\vspace{-15pt}
	\begin{center}
		\begin{tikzpicture}[scale=1.5,every node/.style={scale=1}]
			\draw[thick,fill=lightgrey] (0,0) rectangle (2,2);
			\draw[-latex] (0.80,-0.2) -- node[below]{$q^+$} (1.2,0.2);
			\draw[-latex] (0.80,1.8) -- node[above]{$q^-$} (1.2,2.2);
			\draw[text=black,thick] (1.8,2.25) node {$T_B^-=1$};
			\draw[text=black,thick] (1.8,-0.25) node {$T_B^+=3$};
			\draw[-latex] (0.80,1.2) -- node[right]{$g=1$} (0.8,0.8);
			\draw[-latex] (-0.6,-0.2) -- node[below]{$x$} (-0.1,-0.2);
			\draw[-latex] (-0.6,-0.2) -- node[left] {$y$} (-0.6,0.3);
		\end{tikzpicture}	\end{center}
	\vspace{-15pt}
	\caption{Computational setup for the Rayleigh-B\'{e}nard convection test case. }
	\vspace{0pt}
	\label{fig:RB_setup}
\end{wrapfigure}

\begin{figure}[htb]
	\centering
	\begin{subfigure}[t]{0.325\textwidth}
		\includegraphics[width=1\textwidth]{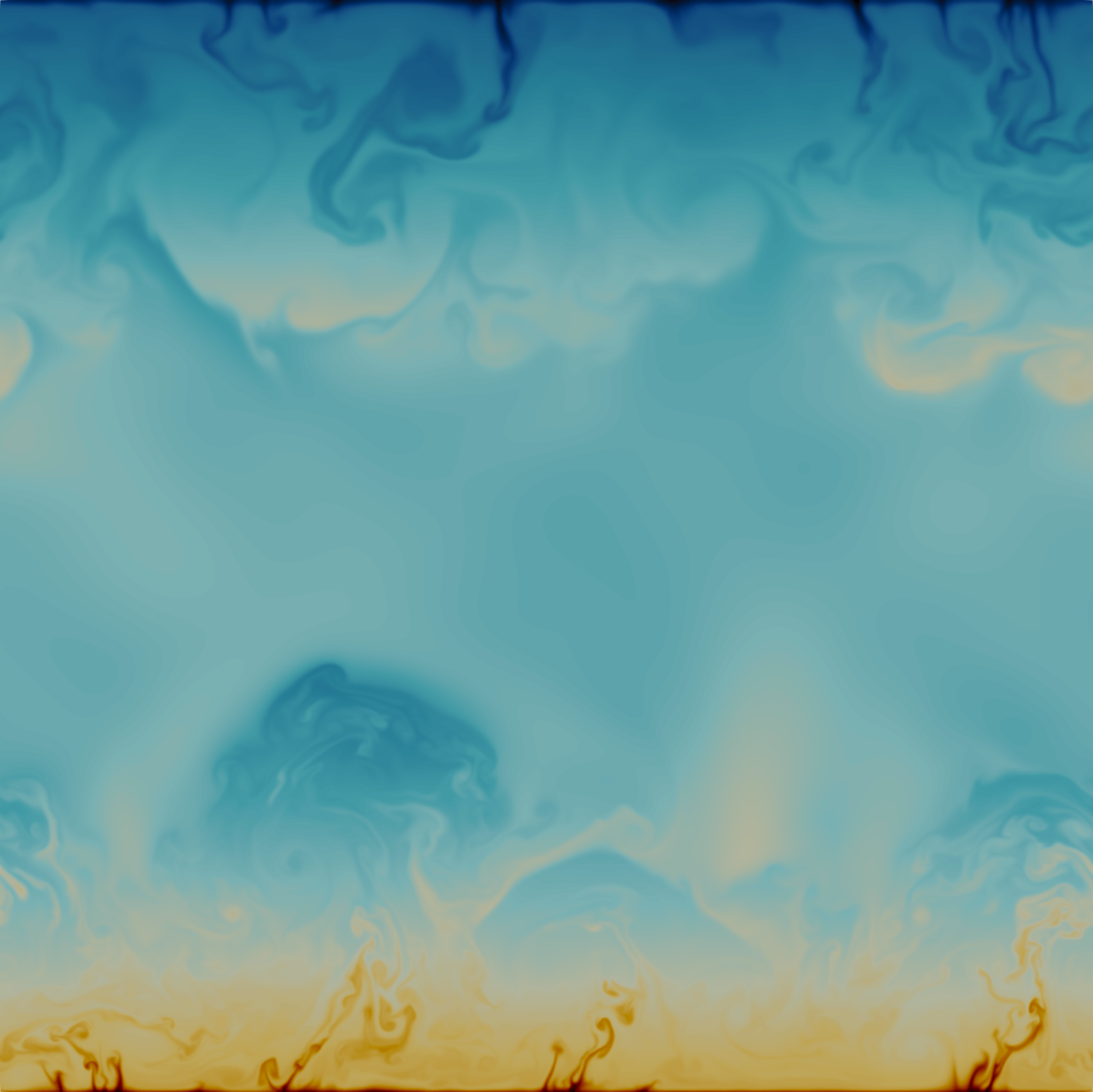}
	\end{subfigure}
	\hfill
	\begin{subfigure}[t]{0.325\textwidth}
		\includegraphics[width=1\textwidth]{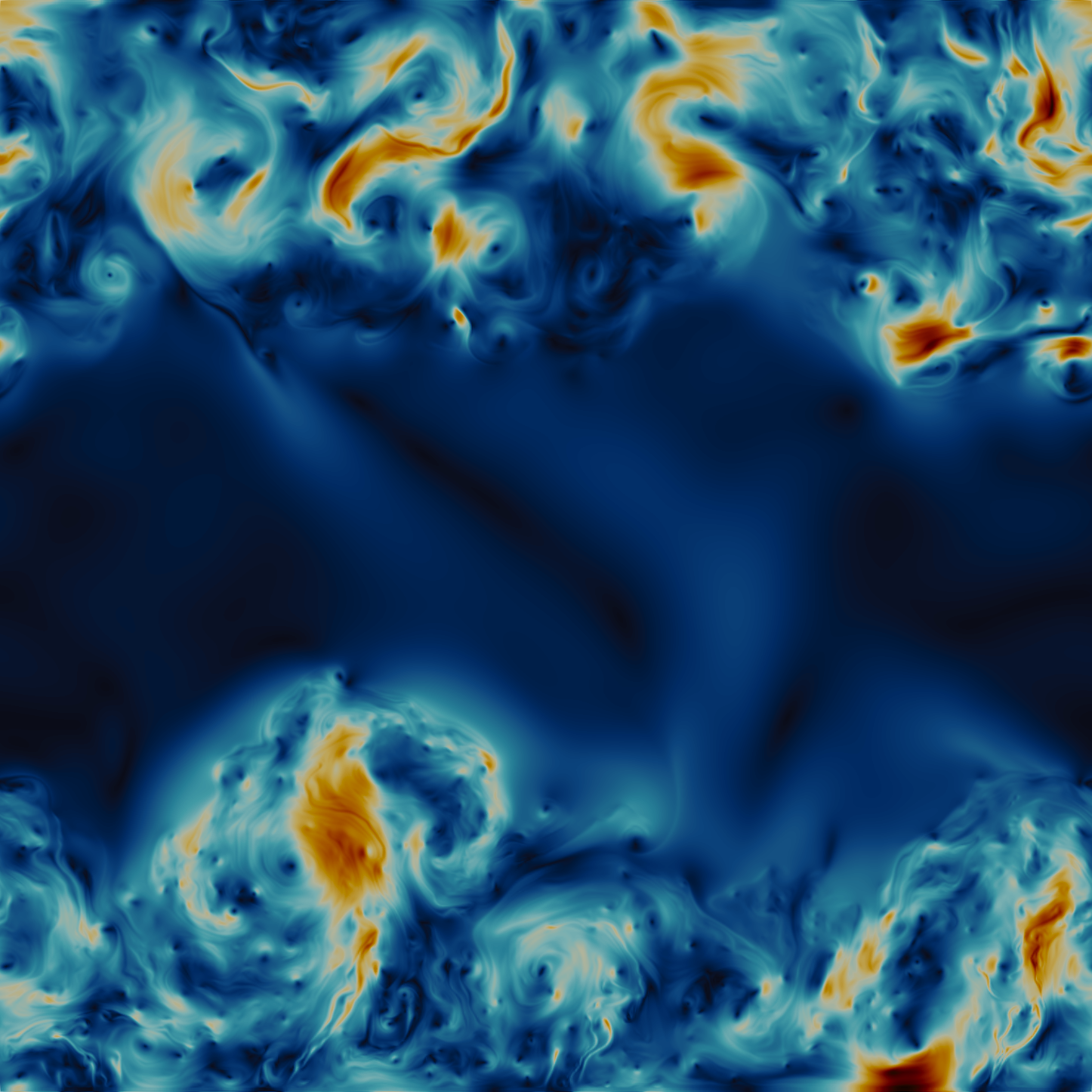}
	\end{subfigure}
	\hfill
	\begin{subfigure}[t]{0.325\textwidth}
		\includegraphics[width=1\textwidth]{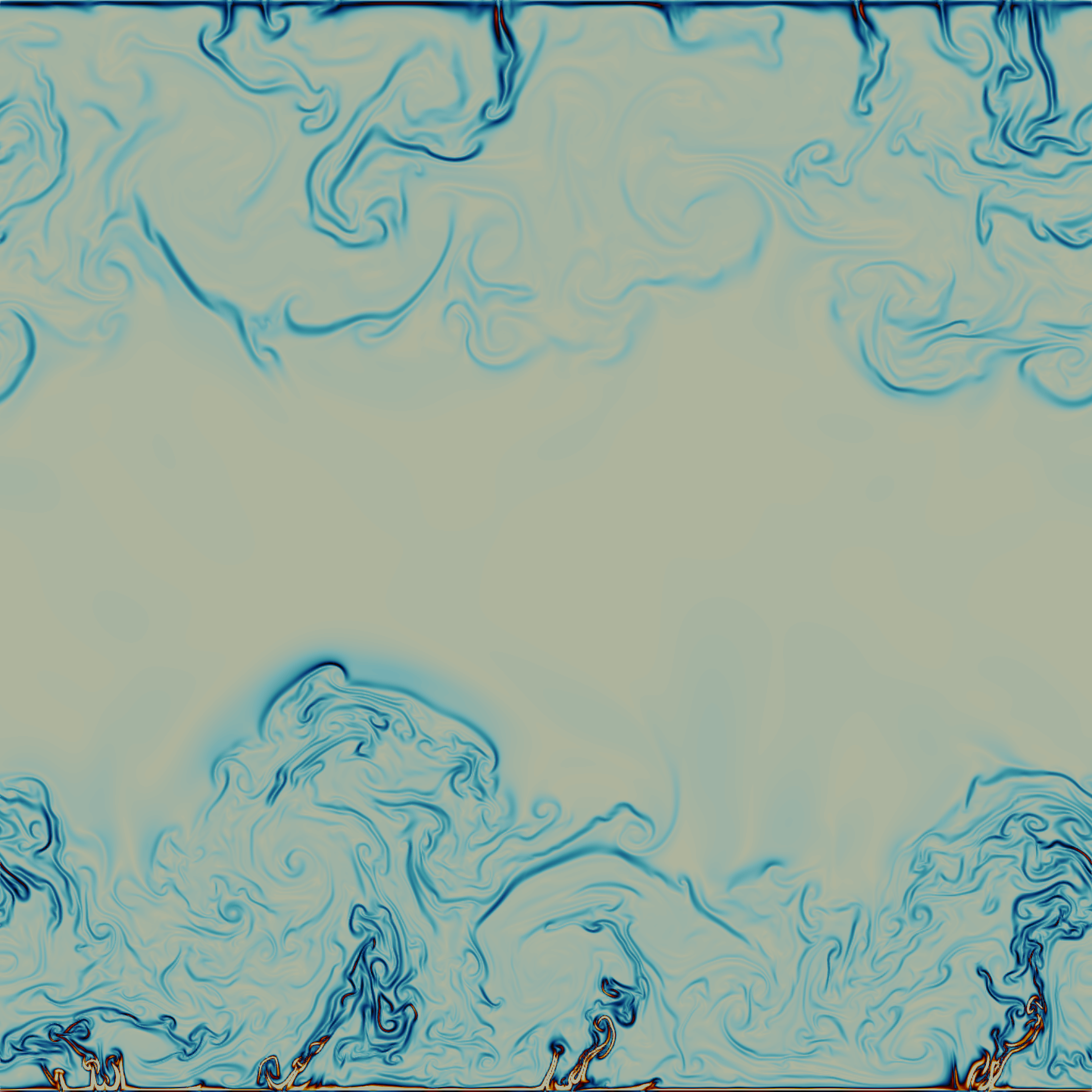}
	\end{subfigure}
	\vfill
	\begin{subfigure}[t]{0.32\textwidth}
		\centering
		\begin{tikzpicture}
		\pgfplotscolorbardrawstandalone[ 
		colormap/erdciceFireH,
		axis on top,
		colorbar horizontal,
		enlargelimits=false,
		point meta min=1.490,
		point meta max=2.695,
		colorbar,
		colorbar style={
			font=\footnotesize,
			width=3.0cm,
			height=0.2cm,
			title={T},
			ylabel near ticks,
			title style={at={(0.50,0.45)},anchor=south},
			xtick={1.49,2.0925,2.695},
			scaled ticks=false,
			ticklabel style={
				/pgf/number format/.cd,
				fixed,
				precision=2,
				fixed zerofill,
				/tikz/.cd
			},
			at={(0.1,0.)},
			anchor=south,
		}]
		\end{tikzpicture}
		\caption{Temperature}
	\end{subfigure}
	\hfill
	\begin{subfigure}[t]{0.32\textwidth}
		\centering
		\begin{tikzpicture}
		\pgfplotscolorbardrawstandalone[ 
		colormap/erdciceFireH,
		axis on top,
		colorbar horizontal,
		enlargelimits=false,
		point meta min=-0.005,
		point meta max=0.196,
		colorbar,
		colorbar style={
			font=\footnotesize,
			width=3.0cm,
			height=0.2cm,
			title={$|\vB|$},
			ylabel near ticks,
			title style={at={(0.50,0.45)},anchor=south},
			xtick={-0.005,0.0955,0.196},
			scaled ticks=false,
			ticklabel style={
				/pgf/number format/.cd,
				fixed,
				precision=3,
				fixed zerofill,
				/tikz/.cd
			},
			at={(0.1,0.)},
			anchor=south,
		}]
		\end{tikzpicture}
		\caption{Velocity magnitude}
	\end{subfigure}
	\hfill
	\begin{subfigure}[t]{0.32\textwidth}
		\centering
		\begin{tikzpicture}
		\pgfplotscolorbardrawstandalone[ 
		colormap/erdciceFireL,
		axis on top,
		colorbar horizontal,
		enlargelimits=false,
		point meta min=0.0,
		point meta max=0.002,
		colorbar,
		colorbar style={
			font=\footnotesize,
			width=3.0cm,
			height=0.2cm,
			title={$|\jB|$},
			ylabel near ticks,
			title style={at={(0.50,0.45)},anchor=south},
			xtick={0.0,0.001,0.002},
			scaled ticks=false,
			ticklabel style={
				/pgf/number format/.cd,
				fixed,
				precision=3,
				fixed zerofill,
				/tikz/.cd
			},
			at={(0.1,0.)},
			anchor=south,
		}]
		\end{tikzpicture}
		\caption{Thermal impulse magnitude}
	\end{subfigure}
	\caption{Temperature (left), velocity (center) and thermal impulse (left) field of a Rayleigh-B\'{e}nard convection at $t=40.0$ 
	with $\lambda=1\cdot 10^{-5}$, computed with the GPR continuum model.}
	\label{fig:RayleighBennard_hpr}
\end{figure}

\begin{figure}[htb]
	\centering
	\begin{subfigure}[t]{0.325\textwidth}
		\includegraphics[width=1.05\textwidth]{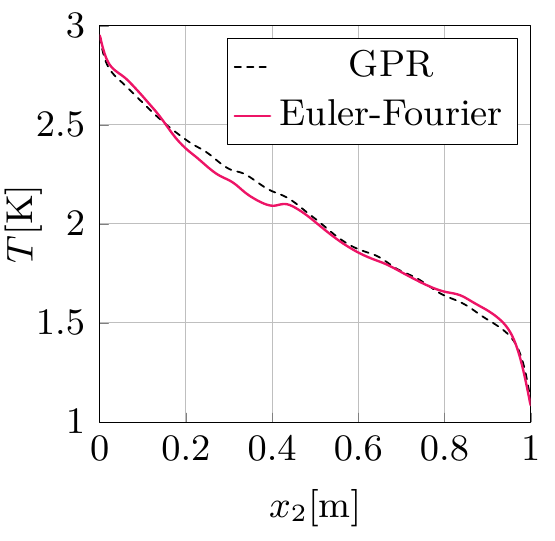}
		\caption{$\lambda=1\cdot 10^{-3}$}
	\end{subfigure}
	\hfill
	\begin{subfigure}[t]{0.325\textwidth}
		\includegraphics[width=1.05\textwidth]{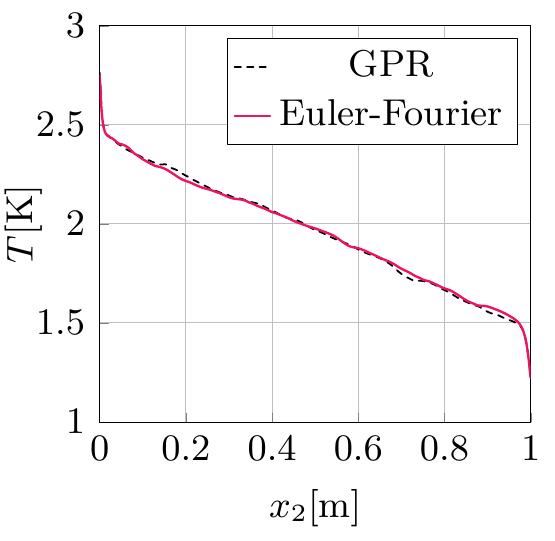}
		\caption{$\lambda=1\cdot 10^{-4}$}
	\end{subfigure}
	\hfill
	\begin{subfigure}[t]{0.325\textwidth}
		\includegraphics[width=1.05\textwidth]{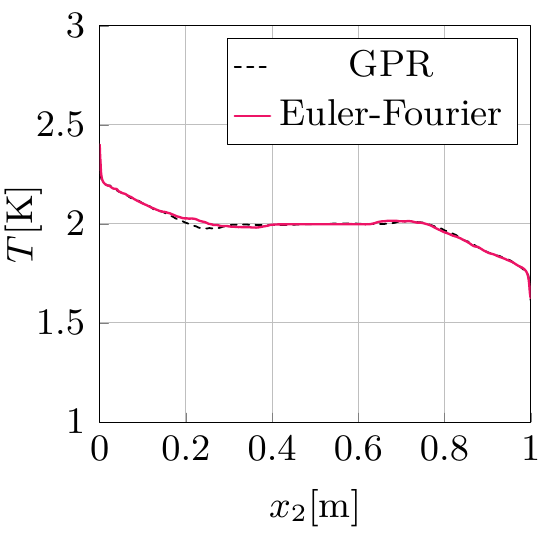}
		\caption{$\lambda=1\cdot 10^{-5}$}
	\end{subfigure}
	\caption{Temperature profiles in y-direction, averaged over the x-direction for different thermal conductivities $\lambda$. 
	The GPR computation is compared against a reference solution, obtained with the Euler-Fourier system.}
	\label{fig:RayleighBennard_statistic}
\end{figure}

The Rayleigh-B\'{e}nard convection is a well-known benchmark, based on the interaction of heat conductivity and gravitational 
forces. In this paragraph, it is used as a two-dimensional test case to assess the computational performance of the present GPR implementation. 
As an initial setup, we consider a two-dimensional domain $\Omega=[0,1]\times[0,1]$ that contains a resting perfect gas with a constant 
temperature $T=2.0$. A gravitational force with gravitation constant $g=1.0$ is imposed in negative y-direction. For the compressible fluid, 
this results in a non-linear hydrostatic pressure profile
\begin{equation}
	p(y)=p_0 \exp{\left( \frac{gy}{T_0 k} \right)}  \quad \text{, with} \quad T_0=2.0, p_0=2.5, k=\frac{c_p}{c_v},
\end{equation}  
with $T_0=2.0$, $p_0=2.5$, $R=\frac{c_p}{c_v}$ and $c_v=0.718$ and $c_p=1.005$. Periodic boundary conditions are imposed in x-direction, 
while heat fluxes $q^+$ and $q^-$ are prescribed at the lower and upper boundaries in y-direction. The heat fluxes are defined similarly 
to Section \ref{subsec:1D_heat_conduction} with the addition of a sine-shaped perturbation in x-direction
\begin{align*}
q^+=\rho h(T-T_B^+) (1+0.1 \text{cos}(8\pi x)), \\
q^-=\rho h(T-T_B^-) (1+0.1 \text{cos}(8\pi x)), \\
\end{align*}
and the transmission coefficient $h=100\lambda$. We consider three different thermal conductivities $\lambda=1\cdot 10^{-3}$, $\lambda=1\cdot 10^{-4}$ 
and $\lambda=1\cdot 10^{-5}$ and chose the thermal relaxation time $\tau$ according to equation \eqref{eq:tau_kin}. Since we study the inviscid GPR 
model in the present work, viscosity is neglected. The final setup is depicted in figure \ref{fig:RB_setup}
 and discretized with $256\times256$ DG elements of degree $N=3$ in a time interval 
$t\in[0,40]$. 
We compare the results obtained with the GPR continuum model with Euler-Fourier computations. Since we neglect physical viscosity, 
the gravitational forces are solely damped by numerical viscosity. Therefore, the characteristic convection process of the Rayleigh-B\'{e}nard 
instability develops, as visualized in figure \ref{fig:RayleighBennard_hpr}. 
For a quantitative comparison between the GPR and Euler-Fourier model, temperature profiles in y-direction are averaged along the x-axis 
at the final computation time $t=40.0$. The resulting temperature statistics are visualized in figure \ref{fig:RayleighBennard_statistic} and 
show a good agreement.    

To quantify the computational efficiency of both schemes, we compare the wall time and the total number of time steps
in table \ref{tab:wall_time_RB}. 
\begin{table}[htb]
	\centering
	\begin{tabular}{c|c|c|c}
		&  $\lambda$ & Time steps  & Wall time [CPU h]\\
		\hline
		GPR           & \multirow{2}{*}{$1\cdot 10^{-3}$}  & $1.76\cdot 10^{5}$  & 743.8  \\
		Euler-Fourier &  & $5.13\cdot 10^{5}$  & 1910.6 \\
		\hline
		GPR           & \multirow{2}{*}{$1\cdot 10^{-4}$}& $1.69\cdot 10^{5}$  & 623.0   \\
		Euler-Fourier &  & $1.69\cdot 10^{5}$  & 468.5  \\
		\hline		
		GPR           & \multirow{2}{*}{$1\cdot 10^{-5}$} & $1.58\cdot 10^{5}$  & 618.8   \\
		Euler-Fourier &  & $1.56\cdot 10^{5}$  & 432.7  \\
		\hline	                                                   
	\end{tabular}
	\caption{Comparison of the number of times steps and the wall time for the Rayleigh-B\'{e}nard convection benchmark with the GPR model 
	and the Euler-Fourier system.}
	\label{tab:wall_time_RB}	
\end{table}
For the lower heat conductivities $\lambda=1\cdot 10^{-4}$ and $\lambda=1\cdot 10^{-5}$, the number of time steps is nearly identical for both methods. 
A different trend is observed for the highest thermal conductivity $\lambda=1\cdot 10^{-3}$. Here, the Euler-Fourier computation requires about three
times more time steps than the GPR computation. This behavior is a result of the parabolic time step constraint of the Euler-Fourier model, 
which is linked to the thermal diffusivity $d_{\alpha}=\frac{\lambda}{\rho c_p}$. Since the hyperbolic GPR model avoids this constraint,
it outperforms the Euler-Fourier model in the presence of high thermal conductivities. In contrast, when the time step is restricted by the
convective process, the GPR model requires about $30\%$ more wall time due to the additional variables for the thermal impulse.

\subsection{Evaporating LJTS Shock-Tube}
\label{subsec:shocktube_LJTS}
To validate the $\text{HLLP}_{mq}$ and $\text{HLLP}_{m}$ two-phase Riemann solvers, we study an evaporating shock tube setup, introduced by 
Hitz et al. \cite{Hitz2021}. It is defined as a Riemann problem for the LJTS fluid with piece-wise constant initial conditions
\begin{equation*}
	\QB(\xB,0) = 
	\begin{cases}
	(\rho,u,T)^T =(0.6635,0.0,0.9)^T   & \quad \text{, for} \; x < 0 \; \text{(liq)},\\
    (\rho,u,T)^T =(0.013844,0.0,0.8)^T & \quad \text{, for} \; x > 0 \; \text{(vap)},
	\end{cases}  
\end{equation*}
within a computational domain $\Omega=[-200,1000]$. The left state is in a saturated liquid state, while the right state consists of superheated 
vapor. The initial conditions are provided in non-dimensionalized form with the reference length $\sigma_{\text{ref}}= 1 \text{\AA}$, the reference 
energy $\frac{\epsilon_{\text{ref}}}{k_B}=1K$, the reference mass $m_{\text{ref}}=1u$ and the reference time 
$t_{\text{ref}}=\sigma_{\text{ref}}\sqrt{m_{\text{ref}}/\epsilon_{\text{ref}}}$ following Merker et al. \cite{Merker2012}. We apply the EOS of 
Heier et al. \cite{Heier2018} for the LJTS fluid and model the heat conductivities in the liquid and vapor with the models of Lautenschl\"ager 
and Hasse \cite{Lautenschlaeger2019} and Lemmon and Jacobsen \cite{Lemmon2004}, respectively. The thermal relaxation time $\tau$ is defined by
equation \eqref{eq:tau_thermo} according to a thermomass theory model. We discretize the computational domain with $N_{\text{elems}}=300$ DG elements of degree 
$N=3$. In the presence of shocks or the phase boundary, the FV sub-cell scheme with a sub-cell resolution of $N_{FV}=4$ is applied.
The setup is computed until a final time $t=600$. 

\begin{figure}[t]
	\begin{centering}
		\includegraphics[width=1\textwidth]{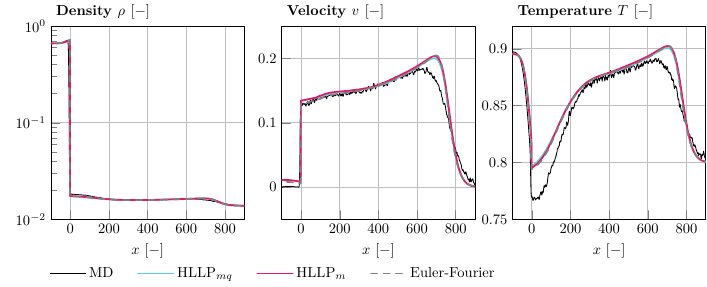}
		\caption{Density, velocity and temperature of the evaporating LJTS shock tube at time $t=600$. 
		The $\text{HLLP}_{mq}$ and $\text{HLLP}_{m}$ Riemann solvers are compared against the Euler-Fourier 
		solution of J\"ons et al. \cite{Joens2023} and a reference solution from molecular dynamic data of Hitz et al. \cite{Hitz2021}.}
		\label{fig:LJTS_ef_hpr_LJTS}
	\end{centering}
\end{figure}

In figure \ref{fig:LJTS_ef_hpr_LJTS}, we compare the solutions of the introduced $\text{HLLP}_{mq}$ 
and $\text{HLLP}_{m}$ Riemann solvers to Euler-Fourier computations by J\"ons et al. \cite{Joens2023} and a molecular dynamics simulation of 
Hitz et al. \cite{Hitz2021}. 
An excellent agreement between the solutions of the $\text{HLLP}_{mq}$ and $\text{HLLP}_{m}$ Riemann solvers can be observed. 
Further, the GPR computations match the solution of the Euler-Fourier system almost perfectly. 

\begin{figure}[!b]
	\begin{centering}
		\includegraphics[width=0.80\textwidth]{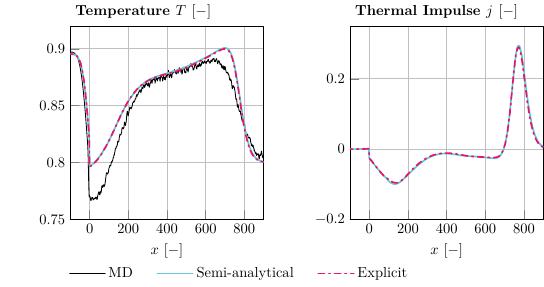}
		\caption{Comparison of the temperature and thermal impulse distribution at $t=600$ computed with a semi-analytical
		and explicit source term integration scheme.}
		\label{fig:LJTS_Time_Integrator}
	\end{centering}
\end{figure}

When compared against the molecular dynamics data, the $\text{HLLP}_{mq}$ and $\text{HLLP}_{m}$ solvers show a good qualitative agreement. 
Deviations are most prominent at the phase boundary, where the GPR computations slightly overpredict the temperature.
A possible cause for the discrepancy could be an underprediction of the evaporation by the thermodynamic closure model. 
Deviations between the GPR results and the MD data away from the interface can be explained by the fact that in this work 
we neglected viscous effects. 

Next, the semi-analytical source term integration scheme of Section \ref{subsec:SourceTermTreatment} is validated. 
Therefore, the shock tube computation is repeated with an explicit source term integration scheme. Due to the
high thermal diffusivity $d_{\alpha}=\frac{\lambda}{\rho c_p}$ of the LJTS fluid, the thermal relaxation time obtained with
equation \eqref{eq:tau_thermo} leads to a mildly stiff behavior of the source and an explicit reference computation 
is affordable. Both methods are compared in figure \ref{fig:LJTS_Time_Integrator} and demonstrate a near-perfect agreement in 
the temperature and thermal impulse profiles.  

\begin{figure}[t]
	\begin{centering}
		\includegraphics[width=1\textwidth]{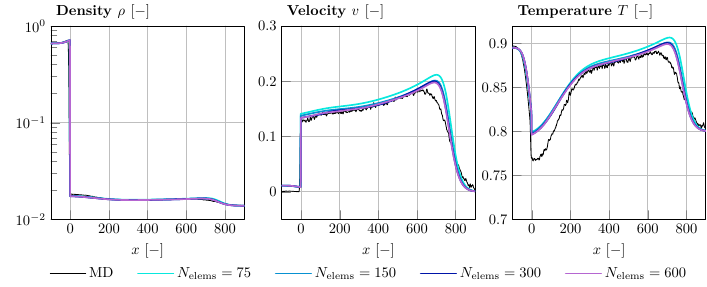}
		\caption{Mesh convergence study with the $\text{HLLP}_{mq}$ two-phase Riemann solver for an evaporating LJTS shock tube. }
		\label{fig:LJTS_convergence}
	\end{centering}
\end{figure}

Further, a mesh convergence study is performed for the evaporating shock tube in figure \ref{fig:LJTS_convergence}. 
Therefore, the computation is repeated for a range of mesh resolutions $N_{\text{elems}}\in[75,150,300,600]$. With an increased resolution, 
a noticeable decrease in the vapor temperatures and velocities is observable. Beyond a resolution of $N_{\text{elems}}=300$, 
a further increase in the element number has only a very minor effect on the solution. 

Finally, the computational efficiency of the GPR method is assessed in comparison to the Euler-Fourier methodology. 
Since the LJTS fluid exhibits a high thermal conductivity $\lambda$ and consequently high thermal diffusivity $d_{\alpha}$, 
the Euler-Fourier scheme suffers from a parabolic time step restriction. 
As indicated in table \ref{tab:wall_time_LJTS}, this results in significantly larger computation times due to an increased number of 
time steps when compared to the GPR model.

In summary, the introduced $\text{HLLP}_{mq}$ and $\text{HLLP}_{m}$ Riemann solvers provide results in good agreement 
with MD reference data. A comparison to a sharp interface study with the Euler-Fourier method revealed a significant performance
advantage of the GPR method due to a lack of a parabolic time step constraint. Further, a mesh convergence study
achieves convergence at a reasonable resolution. Finally, a perfect agreement between the semi-analytical source term integration 
and an explicit reference solution was demonstrated. 

\begin{table}[h!]
	\centering
	\begin{tabular}{c|c|c}
		&  Time steps & Wall time [CPU h]\\
		\hline
		GPR $\text{HLLP}_{mq}$       & $3.76\cdot 10^{3}$  & 29.7  \\
		Euler-Fourier & $3.50\cdot 10^{4}$  & 197.9 \\
        \hline                                   
	\end{tabular}
	\caption{Comparison of the number of times steps and the wall time for the evaporating LJTS shock tube with the GPR model and the Euler-Fourier system.}
	\label{tab:wall_time_LJTS}	
\end{table}

\subsection{Evaporating n-Dodecane Shock-Tube}
\label{subsec:shocktube_n-dodecane}
In the previous paragraph, we investigated an evaporating shock tube for an artificial model fluid, derived from the LJTS potential. 
This facilitated a validation against molecular dynamics simulations to establish the proposed interfacial Riemann solvers. 
With this Section, we extend the study to a shock tube problem for the material n-Dodecane. A Riemann problem is derived 
from an evaporating n-Dodecane shock-droplet setup, reported by J\"ons et al. \cite{Joens2023}, with piecewise constant 
initial liquid and vapor states 
\begin{equation*}
\QB(\xB,0) = 
\begin{cases}
(\rho,u,p)^T =(539.94[\frac{Kg}{m^3}],0.0 [\frac{m}{s}],0.13 [\text{MPa}])^T &  \text{, for} \; x < 0 \; \text{(liq)},\\
(\rho,u,p)^T =(4.3830 [\frac{Kg}{m^3}],0.0 [\frac{m}{s}],0.10 [\text{MPa}])^T & \text{, for} \; x > 0 \; \text{(vap)},
\end{cases}  
\end{equation*}
separated by a phase boundary at $x=0$. N-dodecane is modeled with the Peng-Robinson EOS, based on the material parameter 
listed in table \ref{tab:eos_parameter_preos}. The thermal conductivity is determined with the model of Chung et al. \cite{Chung1988}.
We chose the thermal relaxation time $\tau$ according to equation \eqref{eq:tau_kin}. 
\begin{table}[htb]
	\centering
	\begin{tabular}{c|c|c|c|c}
		$\rho_{c}[\frac{Kg}{m^3}]$ & $p_c[\text{MPa}]$ & $T_c[K]$ & $M[\frac{Kg}{mol}]$ & $\omega [-]$ \\
		\hline
		226.55 & 18.17 & 658.1 & 0.1703 & 0.576 \\
		\hline                                   
	\end{tabular}
	\caption{Material parameter of the Peng-Robinson EOS for n-Dodecane}
	\label{tab:eos_parameter_preos}	
\end{table}

The setup considers a domain $\Omega=[0,0.001]\mathrm{m}$ that is discretized with $600$ DG elements of degree $N=3$. 
At shocks and the interface, a local refinement is applied based on an FV sub-cell scheme with a resolution of $N_{FV}=4$ sub-cells 
per DG element. The setup is advanced until a final computation time $t=2\cdot 10^{-6}\mathrm{s}$ with the novel $\text{HLLP}_{mq}$ and $\text{HLLP}_{m}$ interface 
solvers.

\begin{figure}[t]
	\begin{centering}
		\includegraphics[width=1\textwidth]{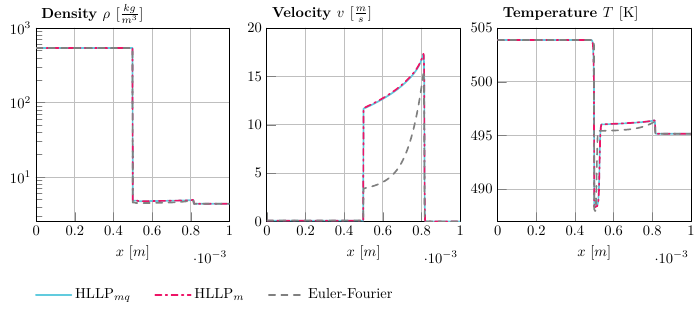}
		\caption{Density, velocity and temperature profile of an evaporating shock tube computation with the material n-Dodecane at $t=2\cdot 10^{-6}\mathrm{s}$. 
		Sharp interface computations with the $\text{HLLP}_{mq}$ and $\text{HLLP}_{m}$ solvers are compared against an Euler-Fourier solution.}
		\label{fig:ndodecane_ef_hpr}
	\end{centering}
\end{figure}

\begin{figure}[!b]
	\begin{centering}
		\includegraphics[width=0.65\textwidth]{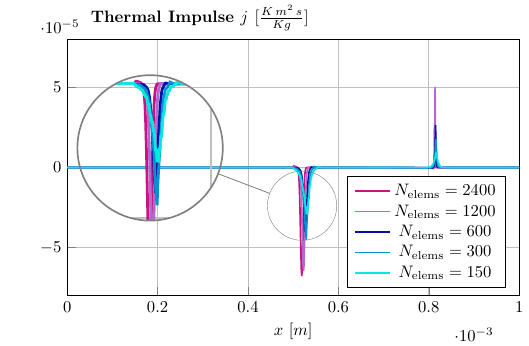}
		\caption{Convergence study with the $\text{HLLP}_{mq}$ solver for the thermal impulse $j$. }
		\label{fig:ndodecane_convergence_j}
	\end{centering}
\end{figure}

Figure \ref{fig:ndodecane_ef_hpr} provides the density, velocity and temperature profiles at the final time $t=2\cdot 10^{-6}\mathrm{s}$. 
Results, obtained with the GPR model are compared against an Euler-Fourier reference solution, computed with the framework of 
J\"ons et al. \cite{Joens2023}. While an excellent agreement can be observed between the novel
$\text{HLLP}_{mq}$ and $\text{HLLP}_{m}$ two-phase Riemann solvers, the Euler-Fourier solution predicts a significantly lower 
velocity in the vapor and a sharper temperature jump at the interface. The wider temperature profile indicates a more dissipative behavior of the GPR model, compared to the Euler-Fourier system. 

\begin{figure}[t]
	\begin{centering}
		\includegraphics[width=1\textwidth]{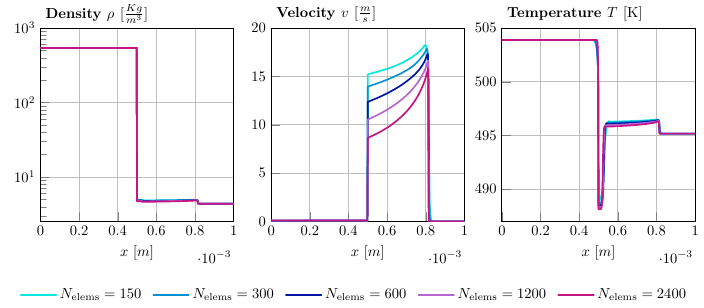}
		\caption{Mesh convergence study with the $\text{HLLP}_{mq}$ two-phase Riemann solver for an evaporating n-Dodecane shock tube. }
		\label{fig:ndodecane_convergence}
	\end{centering}
\end{figure}

This observation can be explained by the thermal impulse distribution, visualized in figure \ref{fig:ndodecane_convergence_j}. 
Due to low significantly lower thermal conductivity of n-Dodecane compared to the LJTS fluid, steep gradients in the temperature cause 
a delta-pulse-like distribution of the thermal impulse. Capturing these sharp solution features is demanding for a discretization scheme
and incurs a particularly high resolution requirement. This argument is further supported by a mesh convergence study in figure 
\ref{fig:ndodecane_convergence}. Similar observations were reported by Peshkov et al in \cite{Peshkov2021} with regard to the simulation
of viscous flows in the GPR framework.

While mesh convergence is not reached for the GPR model at $N_{\text{elems}}=2400$, the investigation suggests that the solution 
approaches the Euler-Fourier reference with increasing mesh resolution. This trend is particularly pronounced in the velocity of 
the freshly evaporated vapor. Furthermore, the temperature dip at the phase boundary appears sharper with increasing mesh resolutions. 

\begin{figure}[htb]	
	\centering
	\includegraphics[width=0.75\textwidth]{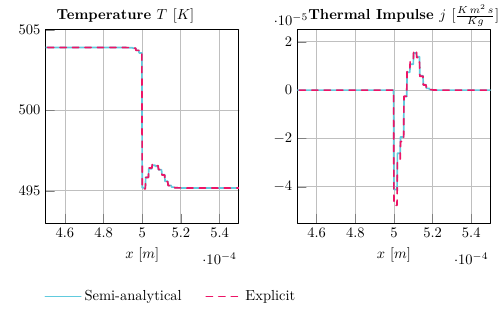}
	\caption{Comparison of the temperature profile and thermal impulse distribution for the semi-analytical and explicit integration of 
	the thermal relaxation source term.}
	\label{fig:source_term_integration_dodecane}
\end{figure}

Next, we analyze the accuracy of the semi-analytical integration scheme for the thermal relaxation source term. 
Again, a reference solution is computed with an explicit source term integration scheme. Due to the prohibitive time step restriction
of the explicit scheme, the comparison is performed for a coarse resolution of $N_{\text{elems}}=150$ elements and evaluated at $t=7\cdot 10^{-8}\mathrm{s}$.
Figure \ref{fig:source_term_integration_dodecane} depicts the temperature and thermal impulse profiles for both simulations. 
Both schemes produce near identical results with a slightly damped thermal impulse profile obtained with the semi-analytical scheme.

Finally, the computational costs of the GPR computation and the Euler-Fourier reference are evaluated in table \ref{tab:wall_time_dodecanee}.
For the given thermal conductivity of n-Dodecane, the time step is dominated by the convection process and the parabolic time step 
constraint of the Euler-Fourier system has no effect. Thus both schemes require roughly the same number of times steps. 
The increased wall time of the GPR implementation is the result of the additional variables for the thermal impulse and the source term integration.

\begin{table}[htb]
	\centering
	\begin{tabular}{c|c|c}
		&  Time steps &  Wall time [CPU h]  \\
		\hline
		GPR $\text{HLLP}_{mq}$        & $2.89\cdot 10^{3}$   & 49.5  \\
		Euler-Fourier & $2.80\cdot 10^{3}$   & 16.0  \\
		\hline                                   
	\end{tabular}
	\caption{Comparison of the number of times steps and the wall time for the evaporating n-Dodecane shock tube with the GPR model and the Euler-Fourier system.}
	\label{tab:wall_time_dodecanee}	
\end{table}

In conclusion, the low thermal conductivity of n-Dodecane compared to the LJTS fluid results in low thermal diffusivity.
As a consequence, the solution exhibits sharp temperature gradients, that cause delta-pulse-like solution features in the thermal impulse.
Therefore, a particularly high-resolution requirement is observed for GPR computations. 
Since the Euler-Fourier system does not suffer from a parabolic time step restriction in this case, it proves to be the more 
efficient method for the given setup. 
 
\subsection{Shock-Droplet Interaction}
\label{Shock_Droplet}
\begin{figure}[h]
	\centering
	\begin{subfigure}{0.49\textwidth}
		\begin{tikzpicture}[scale=1.0,every node/.style={scale=0.8}]
		\draw[thick,fill=grey] (0,0) rectangle (0.6,4);
		\draw[thick,fill=lightgrey] (0.6,0) rectangle (4,4);
		\draw[draw=black,fill=lightblue] (1.3,2) circle (0.48);
		\draw[text=black,thick] (0.3,3.2) node {$\Omega_{v_s}$};
		\draw[text=black,thick] (3.5,3.2) node {$\Omega_{v}$};
		\draw[text=black,thick] (1.9,2.73) node {$\Omega_{l}$};
		\draw[draw=black] (1.8,2.6) -- (1.4,2.3);
		\draw[-latex] (-0.6,-0.2) -- node[below]{$x$} (-0.1,-0.2);
		\draw[-latex] (-0.6,-0.2) -- node[left] {$y$} (-0.6,0.3);
		\draw[draw=black,dashed] (0,2) -- (4,2);
		\end{tikzpicture}
		\caption{Setup 1}
		\label{fig:shock_droplet}	
	\end{subfigure}
	\hfill
	\begin{subfigure}{0.49\textwidth}
		\begin{tikzpicture}[scale=1.0,every node/.style={scale=0.8}]
		\draw[thick,fill=grey] (0,0) rectangle (0.6,4);
		\draw[thick,fill=lightgrey] (0.6,0) rectangle (4,4);
		\draw[draw=black,fill=lightblue] (1.3,2) circle (0.48);
		\draw[draw=black,fill=lightgrey] (1.3,2) circle (0.24);
		\draw[text=black,thick] (0.3,3.2) node {$\Omega_{v_s}$};
		\draw[text=black,thick] (3.5,3.2) node {$\Omega_{v}$};
		\draw[text=black,thick] (2.4,2.4) node {$\Omega_{v_c}$};
		\draw[text=black,thick] (1.9,2.73) node {$\Omega_{l}$};
		\draw[draw=black] (1.8,2.6) -- (1.4,2.3);
		\draw[draw=black] (1.3,2.1) -- (2.1,2.4);
		\draw[-latex] (-0.6,-0.2) -- node[below]{$x$} (-0.1,-0.2);
		\draw[-latex] (-0.6,-0.2) -- node[left] {$y$} (-0.6,0.3);
		\draw[draw=black,dashed] (0,2) -- (4,2);
		\end{tikzpicture}
		\caption{Setup 2}
		\label{fig:shock_droplet_cavity}
	\end{subfigure}
	\caption{Initial setup for a 2D shock-droplet interaction. Setup 1 considers an initially resting n-Dodecane droplet. Setup 2 is an extension 
	of setup 1 with the droplet containing a vapor-filled cavity.}
	\label{fig:2D_shock_droplet_setup}
\end{figure}
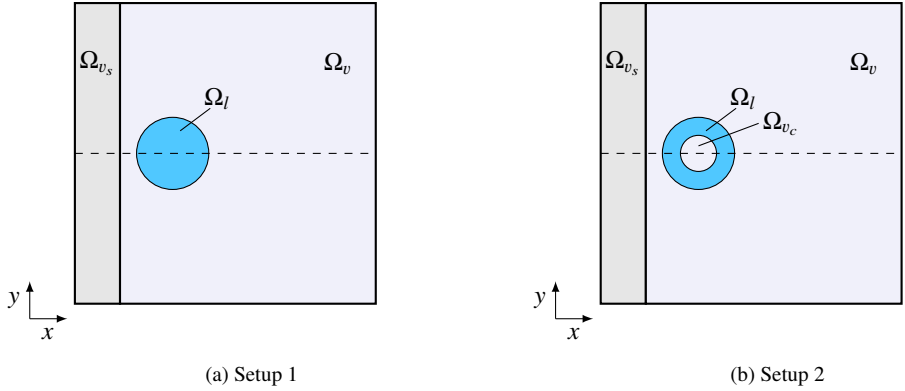

Finally, we consider a two-dimensional shock-droplet interaction with phase transition for n-Dodecane. 
The test case involves surface tension and severe phase boundary deformations and is chosen in order to demonstrate the robustness 
and efficiency of the proposed two-phase Riemann solver for complex two-phase simulations. We adopt a setup, proposed by 
Fechter et al. \cite{Fechter2018} and recently studied by J\"ons et al. \cite{Joens2023} as a test case for their Euler-Fourier two-phase Riemann 
solver. The computational domain $\Omega=[-2.5,7.5]\mathrm{mm}\times[-5.0,5.0]\mathrm{mm}$ contains a planar incident shock at $x=-1.5\mathrm{mm}$ 
and an initially resting droplet of radius $r=1.\mathrm{mm}$ at $x=0.\mathrm{mm}$ under evaporating conditions. The surface tension is set to $\sigma_c=0.009\mathrm{N}\, \mathrm{m}^{-1}$ resulting in a 
Weber number of $We=25549$. The setup is illustrated in figure \ref{fig:shock_droplet} and initial conditions are provided in table 
\ref{table:shock_droplet_ics}. As an extension to the setup of Fechter et al. \cite{Fechter2018}, we consider a second test case with a vapor-filled cavity of radius 
$r=0.5\mathrm{mm}$ inside the droplet, sketched in figure \ref{fig:shock_droplet_cavity}. 

\begin{table}[htb]
	\centering
	\begin{tabular}{c|c|c|c|c}
		& $\rho[\text{kg}\,\text{m}^{-3}]$ & $u_1[\mathrm{m}\, \mathrm{s}^{-1}]$ & $p[\mathrm{MPa}]$\\
		\hline
		Vapor (pre-shock) $\Omega_{v}$   & 4.383  & 0.0    & 0.10  \\
		Vapor (post-shock) $\Omega_{v_s}$ & 9.696  & 108.87 & 0.227 \\
		\hline
		Vapor (cavity) $\Omega_{v_c}$ & 4.383  & 0.0    & 0.10  \\	
		Liquid $\Omega_{l}$   & 539.94 & 0.0    & 0.13  \\
		\hline                                                   
	\end{tabular}
	\caption{Initial conditions of the evaporating n-Dodecane shock-droplet interaction.}
	\label{table:shock_droplet_ics}	
\end{table}

The computational domain $\Omega$ is discretized by $240\times120$ DG elements with a polynomial degree in a range of $N=[2,4]$ and a FV sub-cell 
resolution of $N_{FV}=9$. The sharp phase interface is always discretized by FV sub-cells, leading to an effective resolution of $432$ DOFs per 
bubble diameter.
Due to the symmetric setup, we only compute half of the domain $\Omega$ and impose symmetry boundary conditions along the x-axis. 
For the remaining boundaries, we impose non-reflecting boundary conditions at the right and top and an inflow boundary condition at the left. 
The numerical flux is computed by an HLLC Riemann solver in the bulk and by the $\text{HLLP}_{mq}$ Riemann solver at the phase interface. 
\begin{figure}[htb!]	
	\centering
	\begin{subfigure}[t]{0.490\textwidth}
		\includegraphics[width=1.0\textwidth]{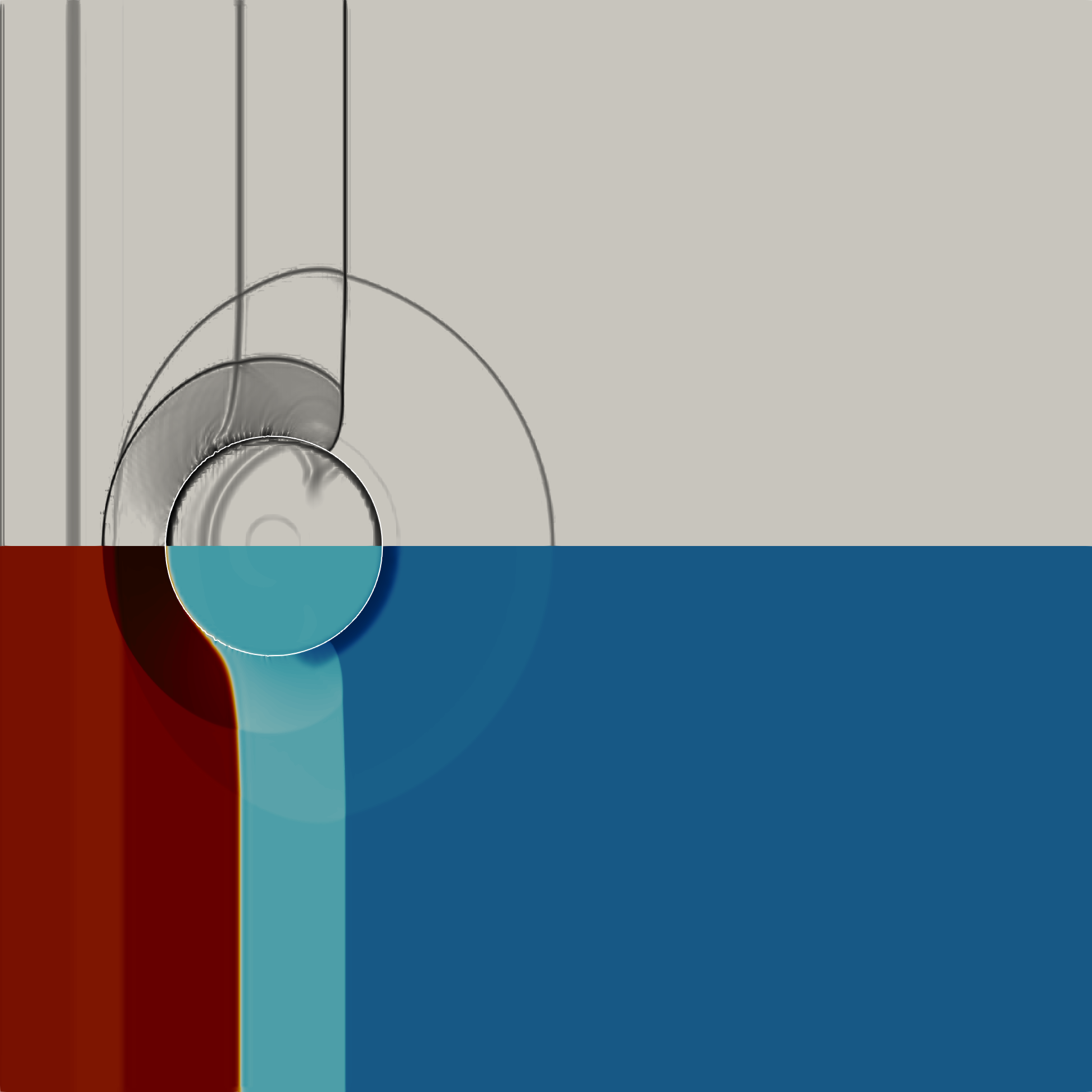}
		\caption{$t=10\mu \mathrm{s}$}
		\label{fig:shock_drop_T_10}
	\end{subfigure}
	\hfill
	\begin{subfigure}[t]{0.490\textwidth}
		\includegraphics[width=1.0\textwidth]{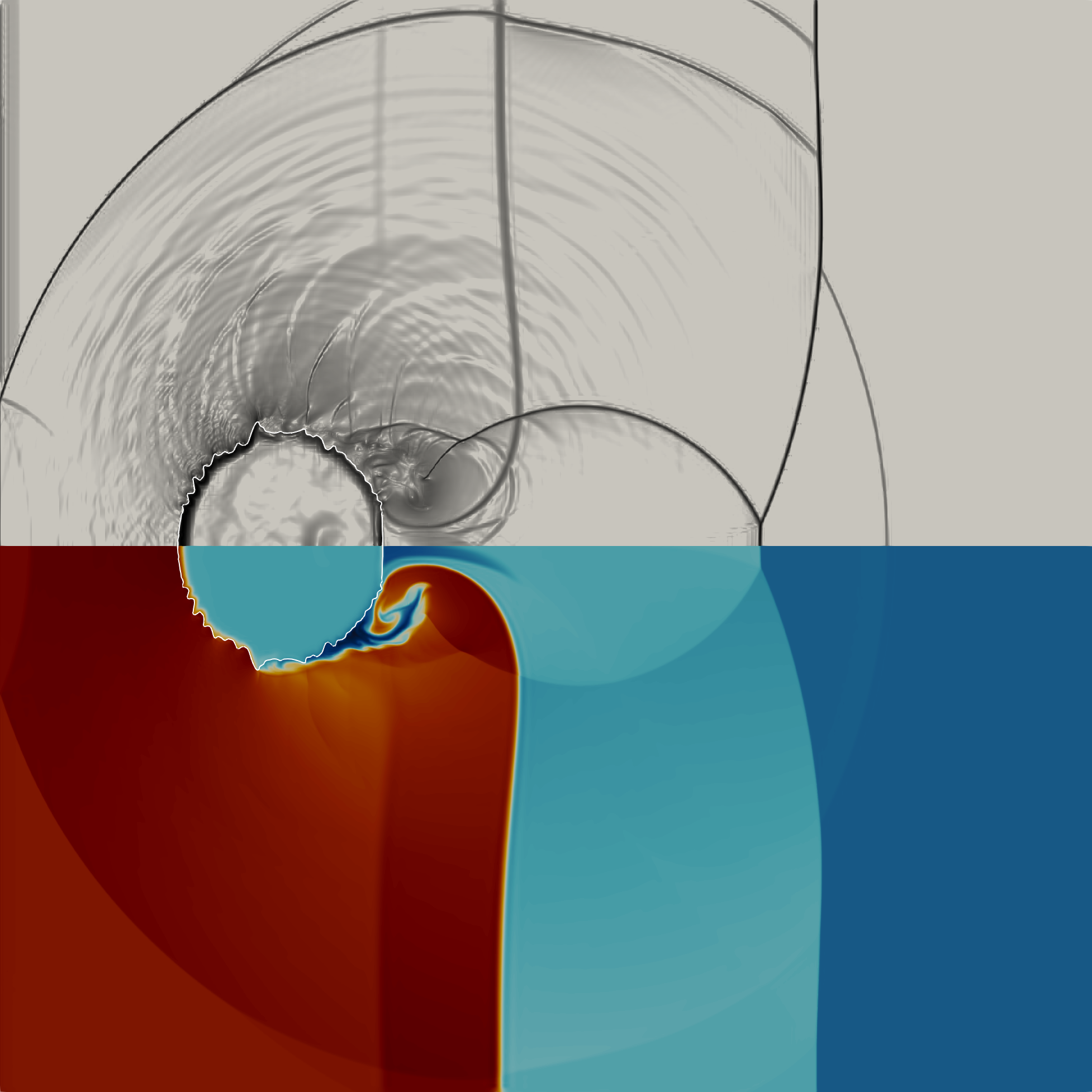}
		\caption{$t=30\mu \mathrm{s}$}
		\label{fig:shock_drop_T_30}
	\end{subfigure}
	\vfill
	\begin{subfigure}[t]{0.490\textwidth}
		\includegraphics[width=1.0\textwidth]{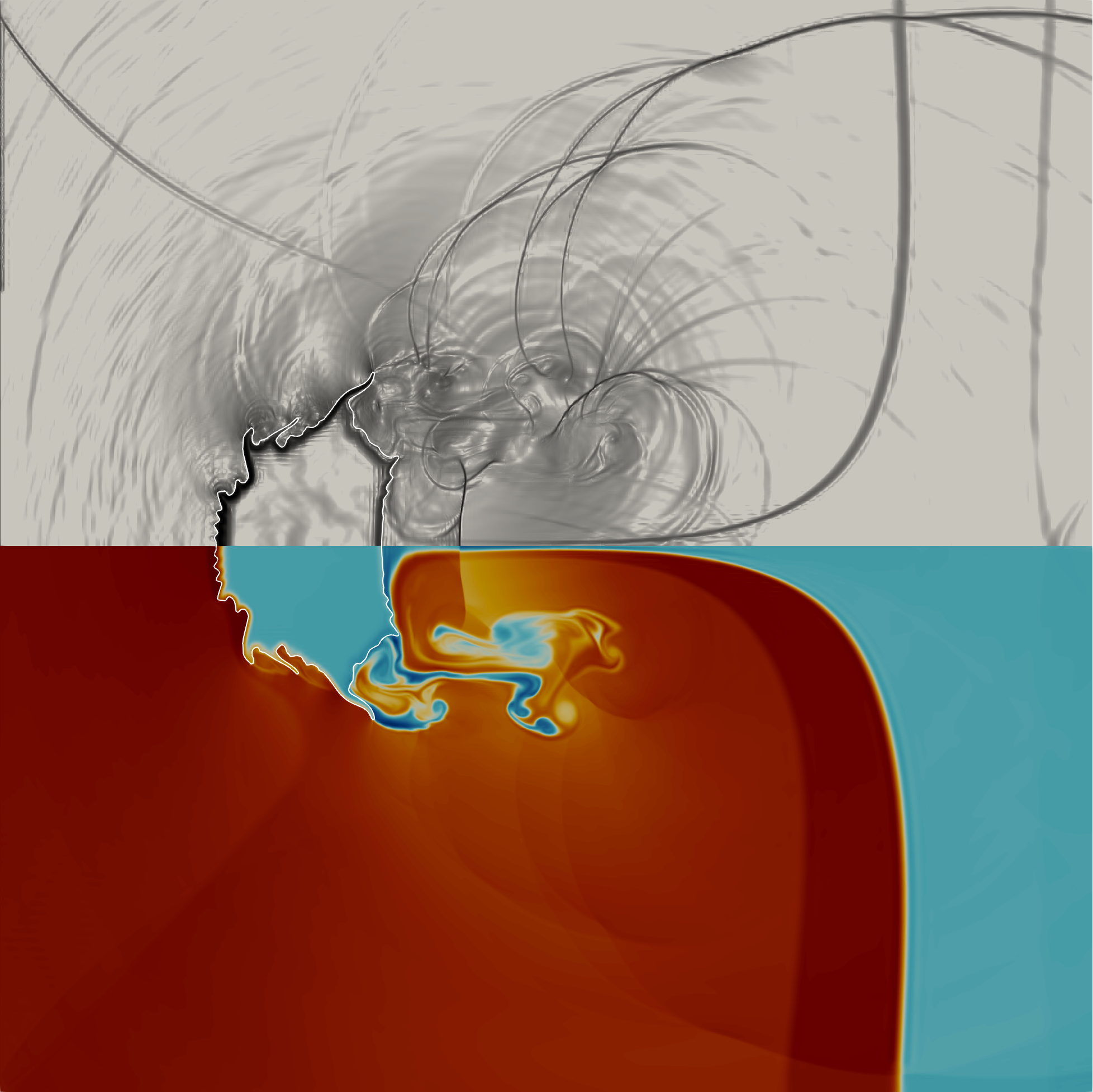}
		\caption{$t=60\mu \mathrm{s}$}
		\label{fig:shock_drop_T_60}
	\end{subfigure}
	\hfill
	\begin{subfigure}[t]{0.490\textwidth}
		\includegraphics[width=1.0\textwidth]{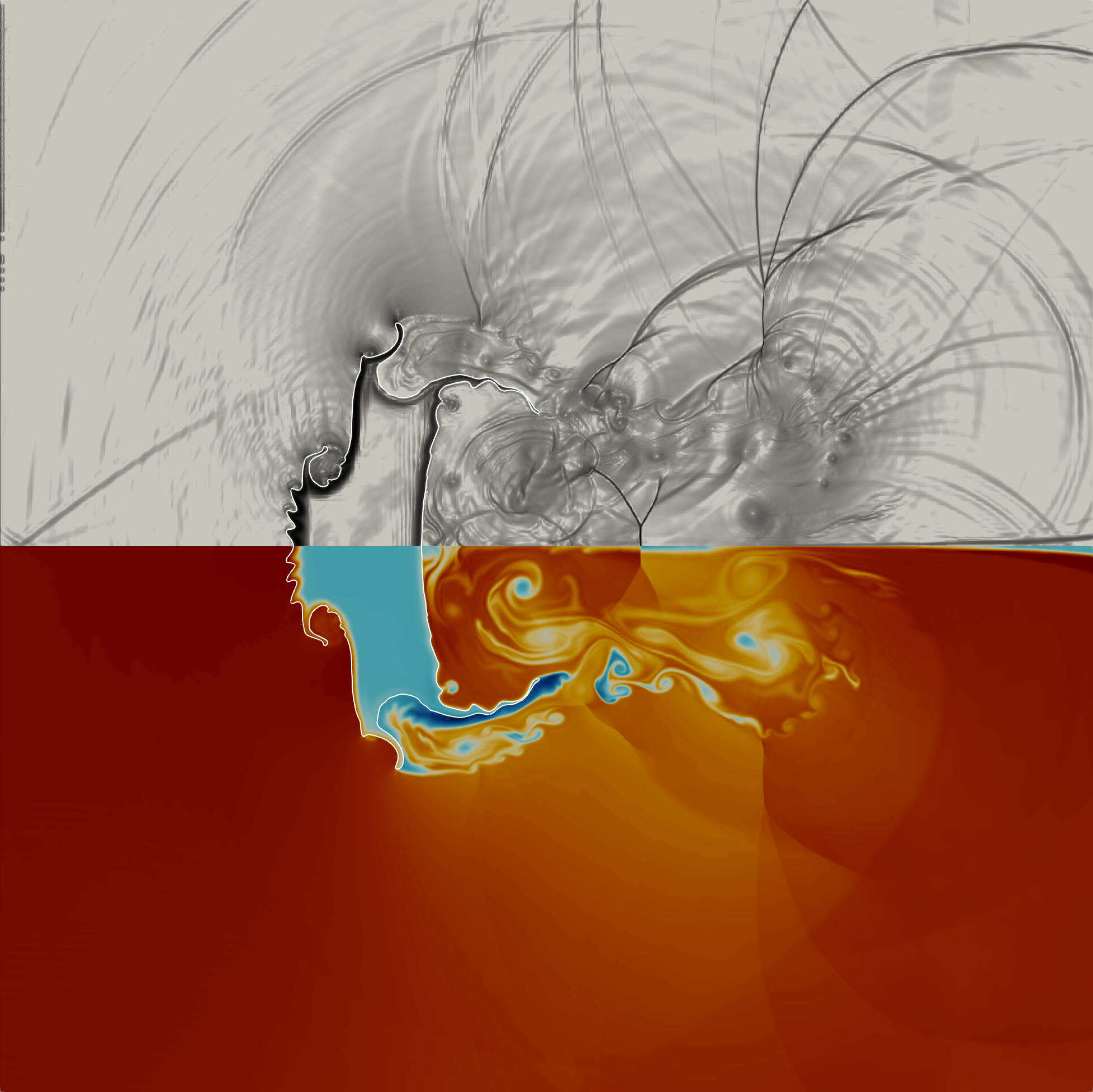}
		\caption{$t=100\mu \mathrm{s}$}
		\label{fig:shock_drop_T_100}
	\end{subfigure}
	\vfill
	\begin{subfigure}[t]{0.495\textwidth}
		\centering
		\begin{tikzpicture}
		\pgfplotscolorbardrawstandalone[ 
		colormap/erdciceFireH,
		axis on top,
		colorbar horizontal,
		enlargelimits=false,
		point meta min=481.0,
		point meta max=541.7,
		colorbar,
		colorbar style={
			font=\footnotesize,
			width=4cm,
			height=0.2cm,
			title={Temperature$[K]$},
			ylabel near ticks,
			title style={at={(0.50,0.45)},anchor=south},
			xtick={481.0,500,520,541.7},
			scaled ticks=false,
			ticklabel style={
				/pgf/number format/.cd,
				fixed,
				precision=1,
				fixed zerofill,
				/tikz/.cd
			},
			at={(0.1,0.)},
			anchor=south,
		}]
		\end{tikzpicture}
	\end{subfigure}
	\hfill
	\begin{subfigure}[t]{0.495\textwidth}
		\centering
		\begin{tikzpicture}
		\pgfplotscolorbardrawstandalone[ 
		colormap/XRay,
		axis on top,
		colorbar horizontal,
		enlargelimits=false,
		point meta min=3.0,
		point meta max=5.7,
		colorbar,
		colorbar style={
			font=\footnotesize,
			width=4cm,
			height=0.2cm,
			title={Schlieren},
			ylabel near ticks,
			title style={at={(0.50,0.45)},anchor=south},
			xtick={3.0,3.9,4.8,5.7},
			scaled ticks=false,
			ticklabel style={
				/pgf/number format/.cd,
				fixed,
				precision=1,
				fixed zerofill,
				/tikz/.cd
			},
			at={(0.1,0.)},
			anchor=south,
		}]
		\end{tikzpicture}
	\end{subfigure}
	\caption{Numerical schlieren images analyzed temperature fields of an n-Dodecane shock-droplet interaction with phase transition using 
	the $\text{HLLP}_{mq}$ two-phase Riemann solver.  }	
	\label{fig:shock_droplet_T}
\end{figure}
We advance the setup until a final time $t=120\mu \mathrm{s}$. The computation is performed on $1024$ processor units and load imbalances caused by 
the adaptive discretization, the interface tracking and the two-phase Riemann solver are balanced with the dynamic load balancing (DLB) scheme, 
introduced in \cite{Appel2021,Mossier2023}. 

Figure \ref{fig:shock_droplet_T} provides numerical schlieren images and the temperature fields 
of the shock-droplet interaction at times $t=10\mu \mathrm{s}$, $t=30\mu \mathrm{s}$, $t=60\mu \mathrm{s}$ and $t=100\mu \mathrm{s}$. 
\begin{figure}[t]	
	\centering
	\begin{subfigure}[t]{0.34\textwidth}
		\includegraphics[width=1.0\textwidth]{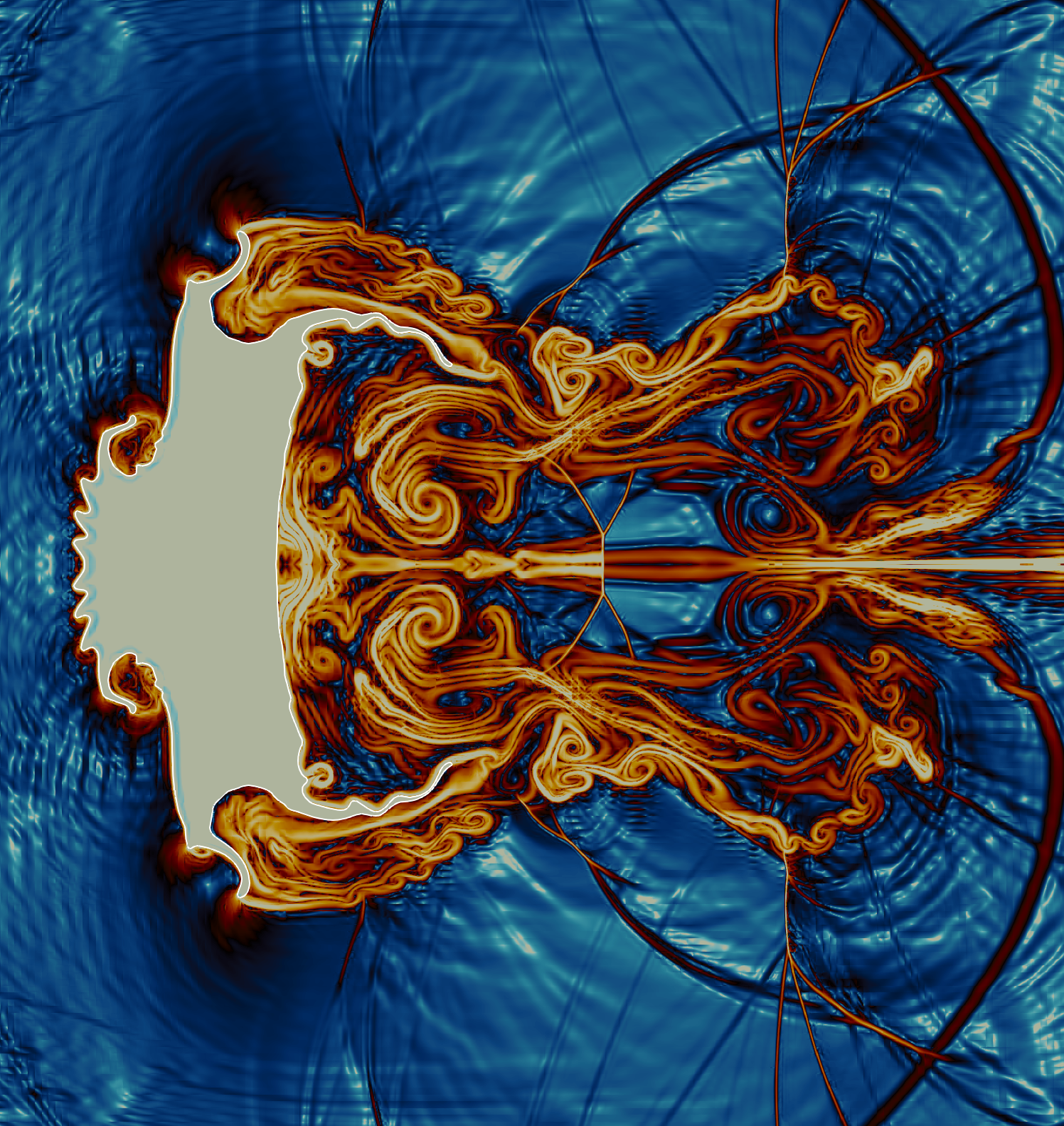}
		\caption{Thermal impulse $|\rho\jB|$}
		\label{fig:thermal_impulse_t100}
	\end{subfigure}
	\hfill
	\begin{subfigure}[t]{0.12\textwidth}
		\centering
		\vspace{-4.0cm}
		\begin{tikzpicture}
		\pgfplotscolorbardrawstandalone[ 
		colormap/erdciceFireL,
		axis on top,
		enlargelimits=false,
		point meta min=1,
		point meta max=4,
		colorbar,
		colorbar right,
		colorbar style={
			font=\tiny,
			width=0.15cm,
			height=2.8cm,
			title={Thermal Impulse},
			ylabel near ticks,
			title style={at={(-1.0,0.40)},anchor=south,rotate=90},
			ytick={1,2,3,4},
			yticklabels={$1.0\cdot 10^{-8}$,$1.0\cdot 10^{-7}$,$1.0\cdot 10^{-6}$,$1.0\cdot 10^{-5}$},
			ticklabel style={
				/pgf/number format/.cd,
				fixed,
				precision=1,
				fixed zerofill,
				/tikz/.cd
			},
			at={(0.0,0.0)},
			anchor=north,
		}]
		\end{tikzpicture}
	\end{subfigure}
	\hfill
	\begin{subfigure}[t]{0.36\textwidth}
		\includegraphics[width=1.0\textwidth]{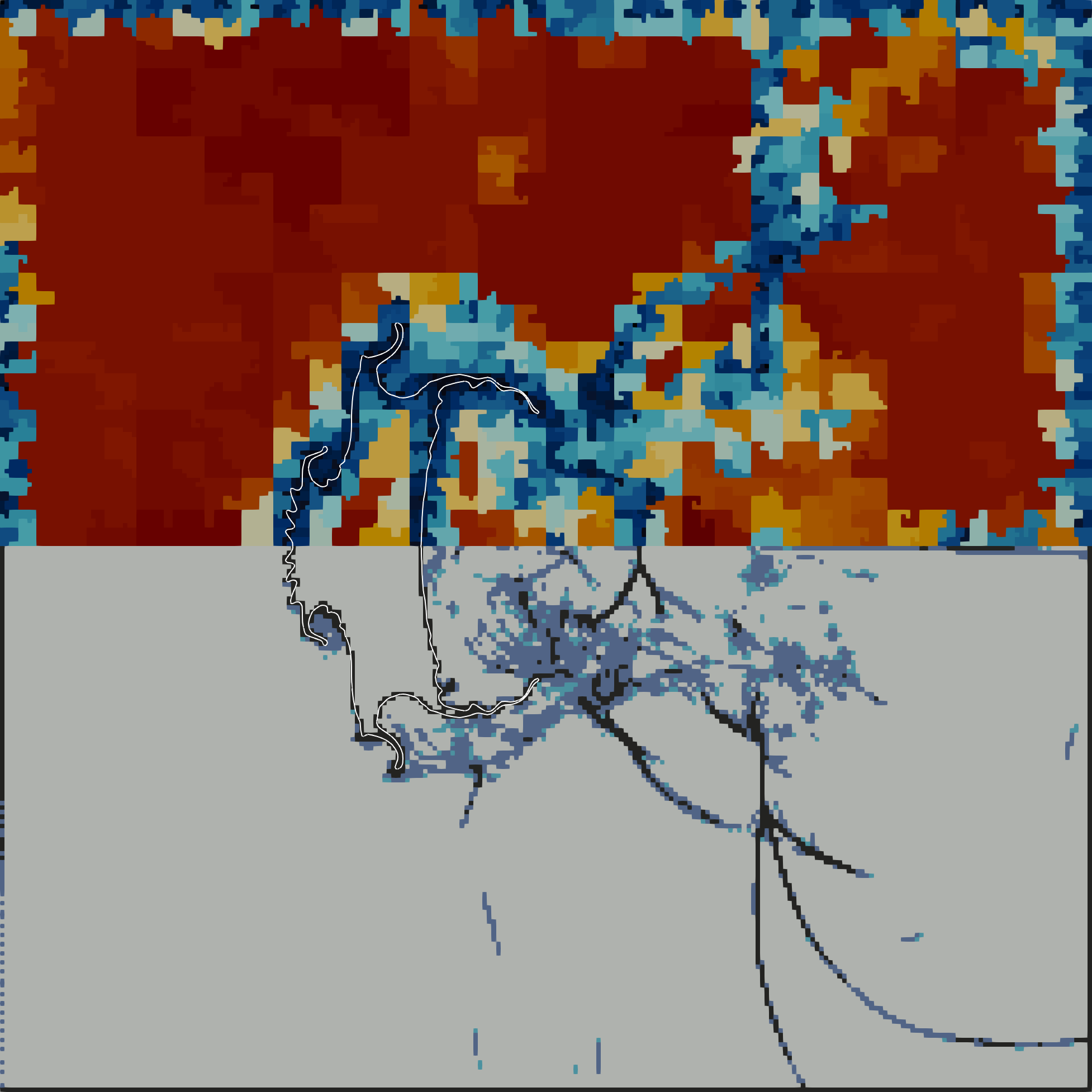}
		\caption{Adaptive discretization}
		\label{fig:shock_drop_hp_t100}
	\end{subfigure}
	\hfill
	\begin{subfigure}[t]{0.12\textwidth}
		\centering
		\vspace{-4.6cm}
		\begin{tikzpicture}
			\pgfplotscolorbardrawstandalone[ 
			colormap/erdciceFireH,
			axis on top,
			enlargelimits=false,
			point meta min=1,
			point meta max=64,
			colorbar,
			colorbar right,
			colorbar style={
				font=\tiny,
				width=0.1cm,
				height=1.8cm,
				title={Elements per proc},
				ylabel near ticks,
				title style={at={(-1.0,0.40)},anchor=south,rotate=90},
				ytick={1,22.3,43.7,64},
				ticklabel style={
					/pgf/number format/.cd,
					fixed,
					precision=0,
					fixed zerofill,
					/tikz/.cd
				},
				at={(0.0,0.0)},
				anchor=north,
			}]
			\pgfplotscolorbardrawstandalone[ 
			colormap/hp_1,
			axis on top,
			enlargelimits=false,
			point meta min=1,
			point meta max=4,
			colorbar,
			colorbar sampled,
			colorbar right,
			colorbar style={
				font=\tiny,
				colormap access=piecewise const,
				samples=5,
				width=0.1cm,
				height=1.8cm,
				title={Resolution},
				ylabel near ticks,
				title style={at={(-1.0,0.40)},anchor=south,rotate=90},
				ytick={1.0,2.0,3.0,4.0},
				yticklabels={DG{,}N=2,DG{,}N=3,DG{,}N=4,FV},
				ticklabel style={
					font=\tiny,
					/pgf/number format/.cd,
					fixed,
					precision=0,
					fixed zerofill,
					/tikz/.cd
				},
				at={(0.0,-2.15cm)},
				anchor=north,
			}]
			\end{tikzpicture}
	\end{subfigure}
	\caption{The left plot shows the absolute thermal impulse $|\rho\jB|$ at time $t=100\mu \mathrm{s}$. The bottom half of the right figure 
	shows the hp-adaptive hybrid DG/FV discretization. The top half depicts the current partition of the domain and indicates the number 
	of elements per processor.}
	\label{fig:shock_drop_J_hp}
\end{figure}
Results are in good qualitative agreement with the reported computations of J\"ons et al. \cite{Joens2023}. At $t=10\mu \mathrm{s}$, 
the incident shock wave passed through most of the droplet and has been reflected at the droplet surface. Inside the droplet, the transmitted wave 
is reflected at the back and a weak retransmitted wave behind the droplet is observable. Due to the initial pressure difference between the droplet 
and the surrounding vapor, a circular shock wave has formed around the droplet and interacts with the incident shock. At the droplet surface, still 
unaffected by the incident shock, a lower temperature in the vapor is visible due to the latent heat of the evaporating droplet. 

At the front of the droplet, a slight increase in the droplet temperature is visible. We follow the explanation of J\"ons et al. and Fechter et al. 
\cite{Fechter2018,Joens2023} and assume that this phenomenon is caused by the condensation of hot vapor impinging on the cool droplet surface. 
At $t=30\mu \mathrm{s}$, the onset of instabilities at the droplet surface due to the high Weber number becomes apparent. They develop into filaments 
at $t=60\mu \mathrm{s}$ and grow until they almost detach from the main liquid body at $t=100\mu \mathrm{s}$. During later stages of the simulation, condensation 
at the front of the droplet becomes stronger, indicated by a heating of the surface. In the back of the droplet, freshly evaporated and cooled 
vapor detaches and mixes with vortical structures in the wake of the droplet. Since the GPR continuum model is used for the bulk fluid, heat conduction is modeled through the thermal impulse. As visualized in figure 
\ref{fig:thermal_impulse_t100}, the absolute value of the thermal impulse is an excellent choice for the visualization of temperature gradients 
in the fluid. 
\begin{figure}[htb!]	
	\centering
	\begin{subfigure}[t]{0.49\textwidth}
		\includegraphics[width=1.0\textwidth]{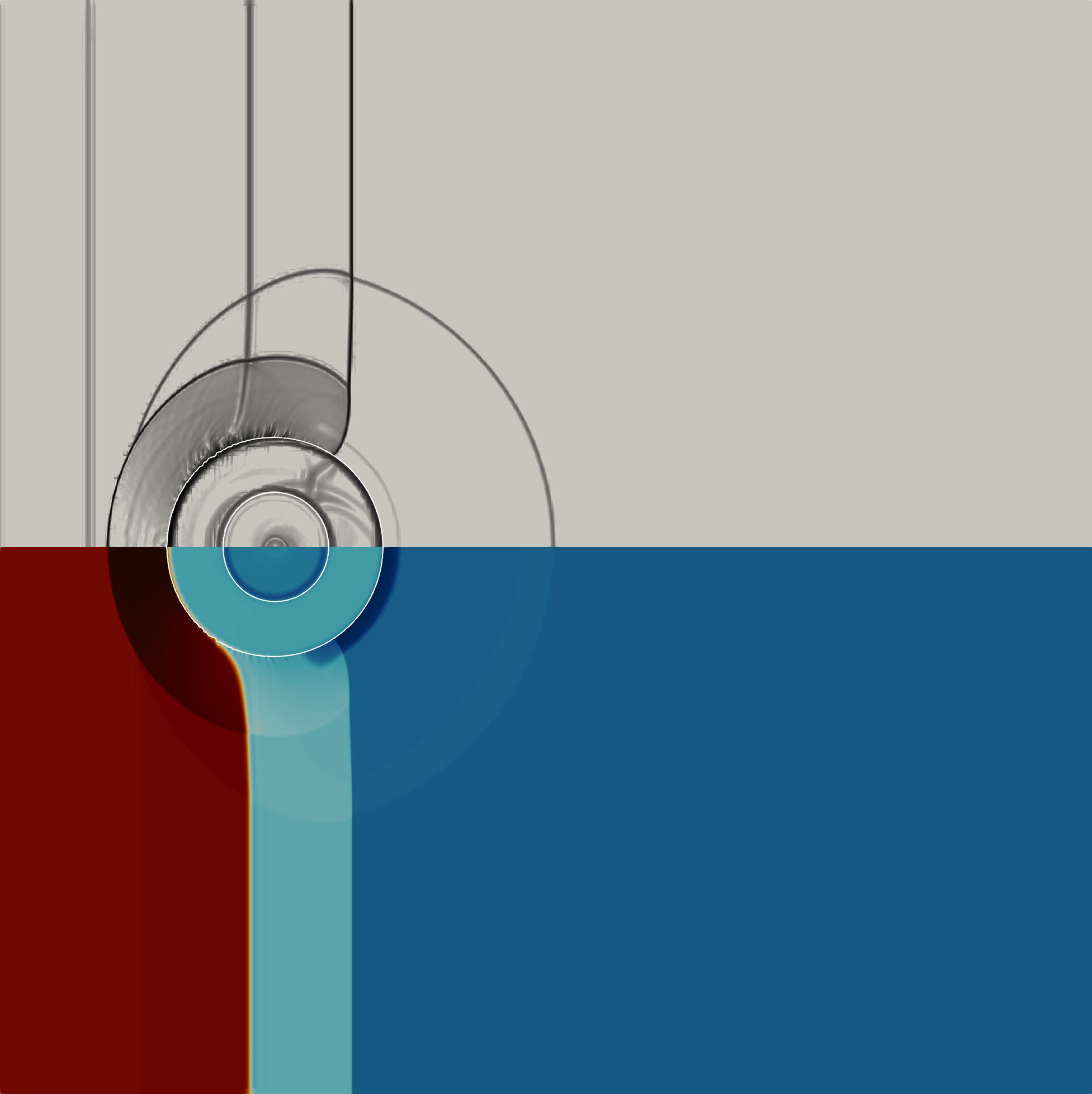}
		\caption{$t=10\mu \mathrm{s}$}
		\label{fig:shock_drop_cavity_T_t10}
	\end{subfigure}
	\hfill
	\begin{subfigure}[t]{0.49\textwidth}
		\includegraphics[width=1.0\textwidth]{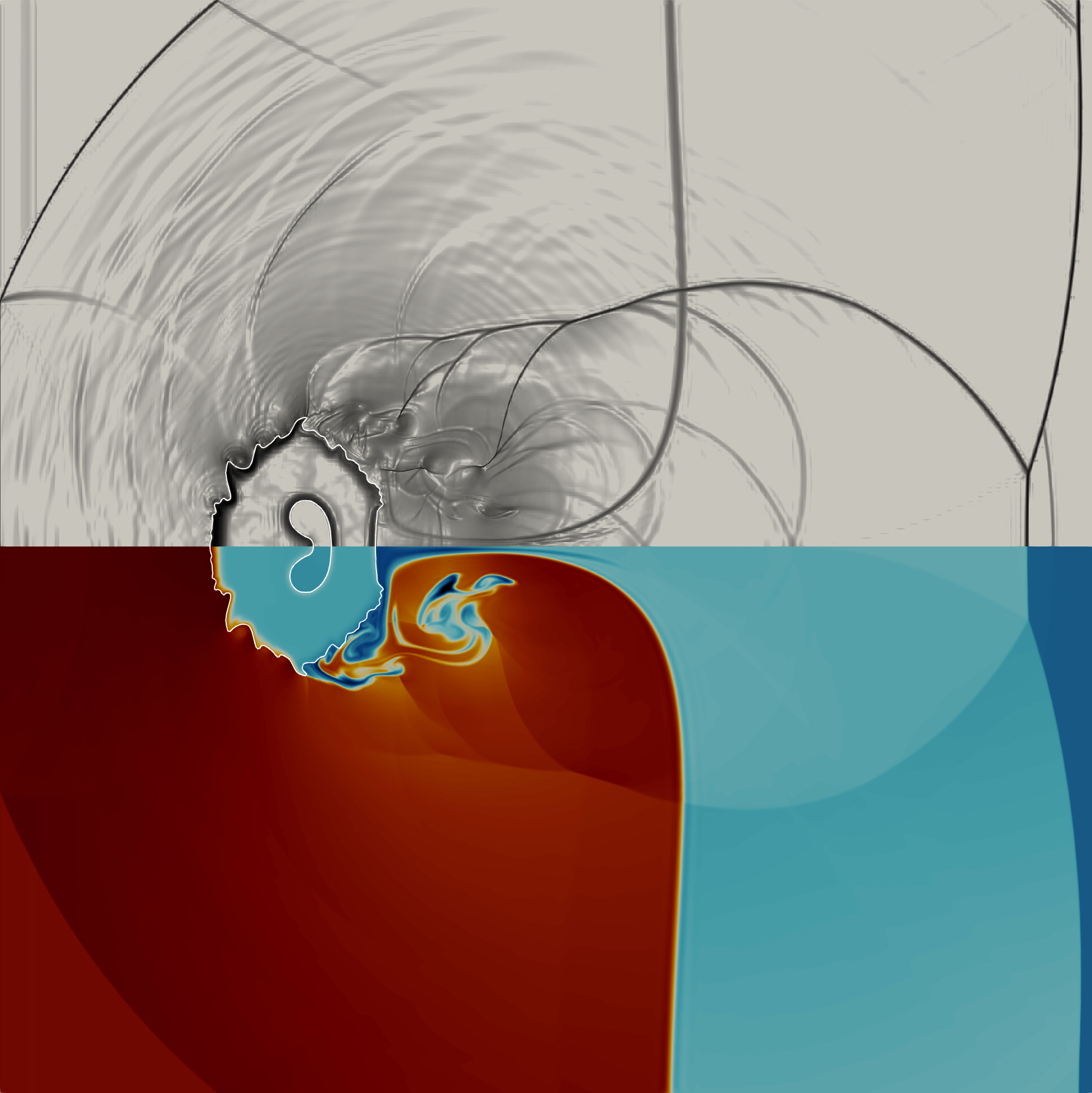}
		\caption{$t=40\mu \mathrm{s}$}
		\label{fig:shock_drop_cavity_T_t40}
	\end{subfigure}
	\vfill
	\begin{subfigure}[t]{0.49\textwidth}
		\includegraphics[width=1.0\textwidth]{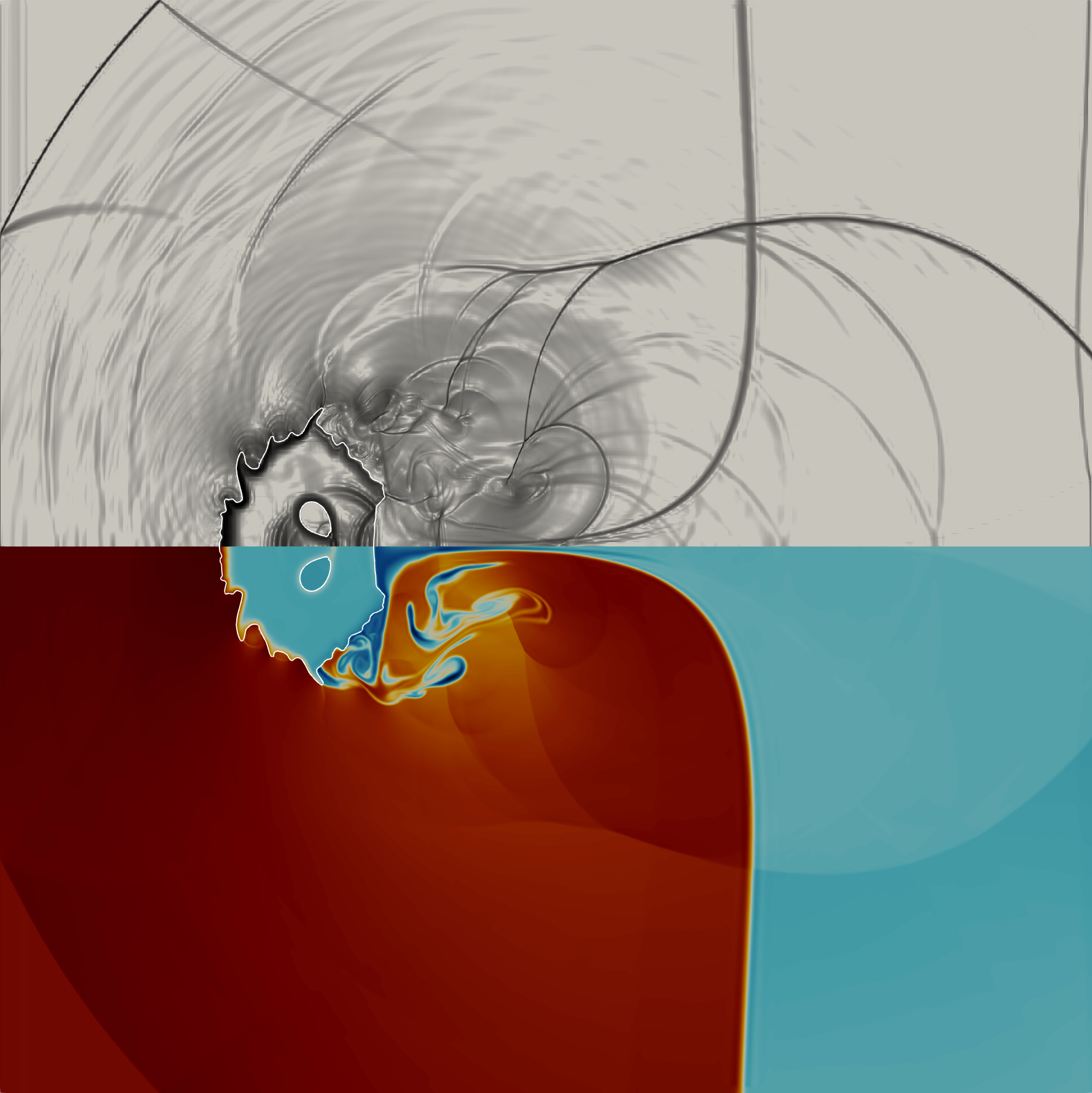}
		\caption{$t=45\mu \mathrm{s}$}
		\label{fig:shock_drop_cavity_T_t45}
	\end{subfigure}
	\hfill
	\begin{subfigure}[t]{0.49\textwidth}
		\includegraphics[width=1.0\textwidth]{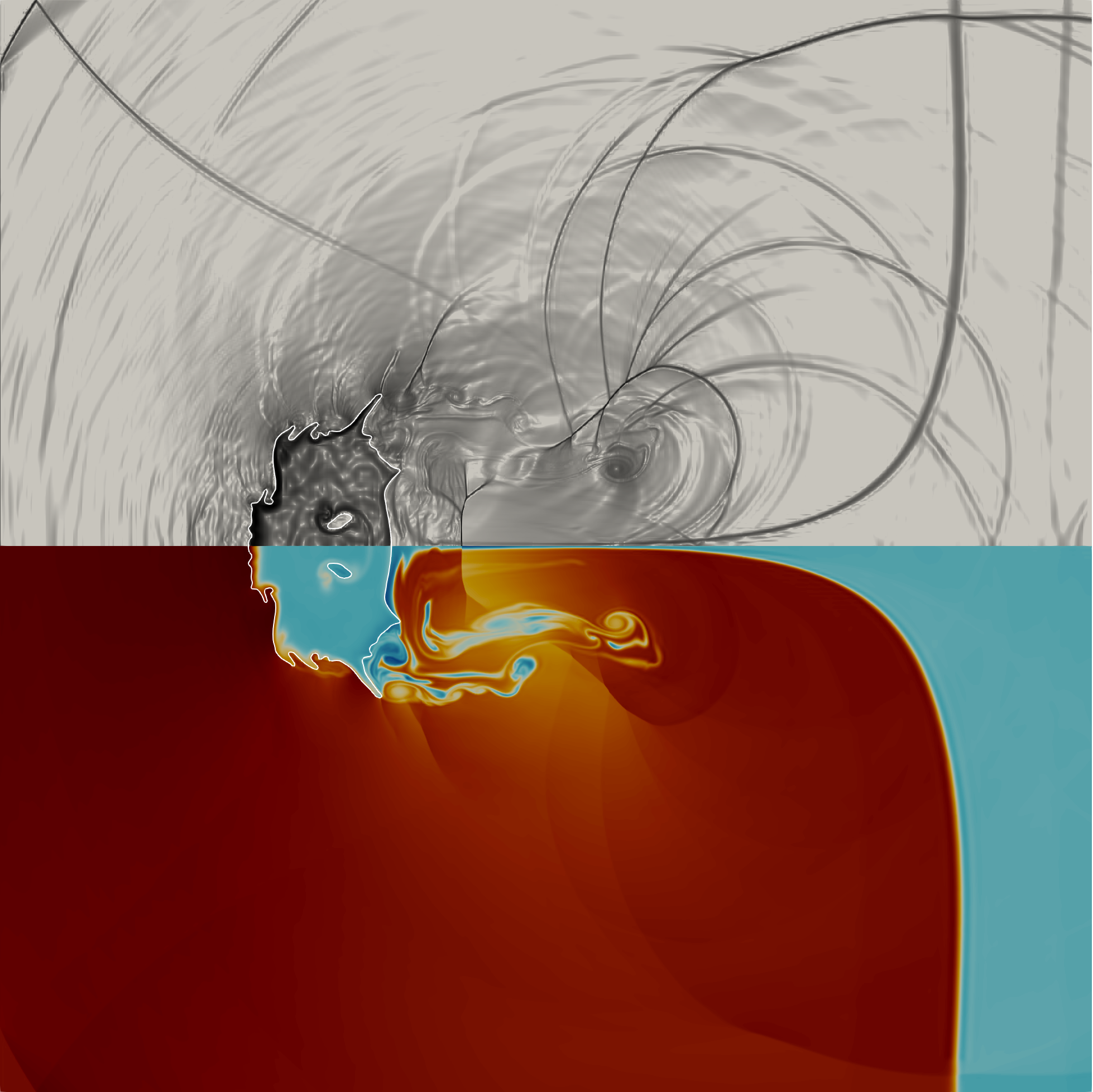}
		\caption{$t=60\mu \mathrm{s}$}
		\label{fig:shock_drop_cavity_T_t60}
	\end{subfigure}
	\begin{subfigure}[t]{0.495\textwidth}
	\centering
	\begin{tikzpicture}
		\pgfplotscolorbardrawstandalone[ 
		colormap/erdciceFireH,
		axis on top,
		colorbar horizontal,
		enlargelimits=false,
		point meta min=481.0,
		point meta max=541.7,
		colorbar,
		colorbar style={
			font=\footnotesize,
			width=4cm,
			height=0.2cm,
			title={Temperature$[K]$},
			ylabel near ticks,
			title style={at={(0.50,0.45)},anchor=south},
			xtick={481.0,500,520,541.7},
			scaled ticks=false,
			ticklabel style={
				/pgf/number format/.cd,
				fixed,
				precision=1,
				fixed zerofill,
				/tikz/.cd
			},
			at={(0.1,0.)},
			anchor=south,
		}]
		\end{tikzpicture}
	\end{subfigure}
	\hfill
	\begin{subfigure}[t]{0.495\textwidth}
		\centering
		\begin{tikzpicture}
		\pgfplotscolorbardrawstandalone[ 
		colormap/XRay,
		axis on top,
		colorbar horizontal,
		enlargelimits=false,
		point meta min=3.0,
		point meta max=5.7,
		colorbar,
		colorbar style={
			font=\footnotesize,
			width=4cm,
			height=0.2cm,
			title={Schlieren},
			ylabel near ticks,
			title style={at={(0.50,0.45)},anchor=south},
			xtick={3.0,3.9,4.8,5.7},
			scaled ticks=false,
			ticklabel style={
				/pgf/number format/.cd,
				fixed,
				precision=1,
				fixed zerofill,
				/tikz/.cd
			},
			at={(0.1,0.)},
			anchor=south,
		}]
		\end{tikzpicture}
	\end{subfigure}
	\caption{Interaction of an incident shock wave with an evaporating n-Dodecane droplet that contains a vapor-filled cavity. 
	Numerical results, obtained with the $\text{HLLP}_{mq}$ two-phase Riemann solver show a collapse of the cavity. }
	\label{fig:shock_drop_cavity_T}
\end{figure}
Finally, insight into the adaptive discretization and DLB is provided by figure \ref{fig:shock_drop_hp_t100}.
The lower half of the plot indicates regions where FV sub-cell limiting and p-refinement is applied. Shocks and the phase interface are detected 
and treated with the refined FV sub-cell grid. At vortical structures in the wake of the droplet, a higher polynomial degree is applied to reduce 
numerical dissipation. The upper half of figure \ref{fig:shock_drop_hp_t100} highlights the partition of the computational domain due to DLB. 
With processor units receiving between $1$ and $64$ elements, a substantial load imbalance between elements can be constituted. Elements discretized 
by the FV sub-cell scheme or the DG scheme with an increased polynomial degree have the highest cost.  

As a final test case, we consider the extended setup with a vapor-filled cavity inside the droplet. At $t=10\mu \mathrm{s}$ a very similar result 
is observed. However, the transmitted incident shock wave is now reflected at the surface of the cavity inside the droplet. Here, a temperature drop 
can be observed due to the evaporation of the liquid droplet into the vapor of the cavity. Between $t=10\mu \mathrm{s}$ and $t=40\mu \mathrm{s}$, the bubble 
undergoes a significant deformation caused by the formation of a high-speed liquid jet along the $y=0$ axis. This liquid jet leads to a bubble 
collapse at $t=45\mu \mathrm{s}$, with shock waves of the collapse reflected at the back of the droplet. The high pressure in the primary cavity and in 
the secondary cavities, formed after the collapse, causes condensation. This leads to a slight increase in the liquid temperature at the surface 
of the primary and secondary cavities during and after the collapse. 

In summary, the two-dimensional shock-droplet interactions demonstrate the applicability of the proposed $\text{HLLP}_{mq}$ Riemann solver to complex 
two-phase simulations with phase transition in the presence of surface tension and significant interface deformations.

\section{Conclusion}
\label{sec:Conclusion}
In this paper, we presented a sharp interface approach for the simulation of compressible two-phase flows with phase transition. 
It employs the first-order hyperbolic continuum model of Godunov, Peshkov and Romenski to describe  
compressible, inviscid heat-conducting fluid flow in the bulk phases. The main contribution of this work is the construction of two novel
interfacial Riemann solvers that provide a thermodynamically consistent coupling at the phase boundary in the presence of
phase transition. The developed Riemann solvers address two key challenges related to the modeling of phase transition 
in the sharp interface context: a loss of self-similarity due to irreversible effects at the interface like entropy production
and heat conduction and avoiding the breakdown of the continuum assumption across the phase boundary. Using the hyperbolic GPR model, 
irreversible effects are treated as relaxation processes and confined to a source term. To obtain a unique and thermodynamically consistent
entropy solution, we employ a local phase transition model to determine the source term and thus predict the entropy production and heat 
dissipation associated with phase transition.

The novel interfacial solvers are constructed from integral jump relations across a simplified wave fan, analogously to the established HLLC methodology.
The equation system is closed by a kinetic relation that employs an entropy estimate from a kinetic theory-based phase transition model. 
With phenomenological force flux relations of Onsager theory, the entropy production is related to the interfacial mass and heat flux.
The resulting non-linear equation system can be solved iteratively. We propose two approximate solvers for the two-phase Riemann problem 
denoted $\text{HLLP}_{mq}$ and $\text{HLLP}_{m}$. While the $\text{HLLP}_{mq}$ solver relies on an iterative solution in both the mass and heat flux,
the simplified $\text{HLLP}_{m}$ solver requires only an iteration in the mass flux.
In addition, we discussed the treatment of the thermal relaxation source term of the GPR model. 
We implemented a recently developed semi-analytical scheme \cite{Chiocchetti2023} that allows to reproduce the Fourier law in the 
stiff relaxation limit. 

To validate the presented method, we investigated a range of representative test cases. First, we compared the inviscid, heat-conducting GPR model 
against Euler-Fourier computations for heat-driven single-phase flows. An excellent agreement between both methods is observed while the GPR 
model demonstrates a superior computational efficiency in the presence of large heat conductivities due to the
lack of a parabolic time step constraint. The two-phase Riemann solvers are validated 
against molecular dynamics data for an evaporating shock tube simulation with the Lennard-Jones shifted and truncated potential. Both the $\text{HLLP}_{mq}$ and $\text{HLLP}_{m}$ 
solvers yield near-identical results and match the molecular dynamics reference data well. Further, a close agreement with an analog Euler-Fourier computation 
is reported. Due to the high thermal conductivity of the LJTS-fluid, the GPR model is computationally advantageous, since it is not restricted
by a parabolic time step constraint. As an additional test case, we studied an evaporating n-Dodecane shock tube. 
While the $\text{HLLP}_{mq}$ and $\text{HLLP}_{m}$ solver produce again near-identical solutions, the GPR results deviated 
significantly from the Euler-Fourier reference. The authors contribute this to a lack of resolution, required for the thermal impulse in case of 
low thermal conductivities. This assumption is supported by mesh convergence studies, suggesting a slow convergence against the Euler-Fourier result. 

Finally, we applied the framework to two-dimensional shock-droplet interactions with phase transition. 
The setups feature surface tension, severe interface deformations and topological changes of the phase boundary. 
To meet the high local resolution requirement at the interface, the developed phase transition kernel was combined with the hp-adaptive 
discretization scheme, introduced in \cite{Mossier2022,Mossier2023}. A robust performance of the proposed interface Riemann solvers
was demonstrated for complex simulations and a good qualitative agreement with similar shock-droplet investigations with the Euler-Fourier 
method \cite{Joens2023} was achieved. 

In the future, we extend the presented framework towards multi-component flows. This allows for validation against experiments with an
evaporating liquid in an inert gaseous atmosphere. Therefore, the bulk phases need to accommodate species transport and diffusion,
while the phase transition models need to be extended to the presence of a multi-component mixture.

\section*{Declarations}

\textbf{Funding} We gratefully acknowledge the support of the German Research Foundation (DFG) for the research reported in this publication through the project
"Droplet Interaction Technologies", grant GRK 2160/1 and GRK 2160/2 (project 270852890), the framework of the research unit FOR 2895 (grant BE 6100/3-1)
and through Germany's Excellence Strategy EXC 2075 (project 390740016). 

Further, S. Chiocchetti acknowledges funding from the European Union's
Horizon Europe research and innovation program under project MoMeNTUM, Marie Skłodowska-Curie grant (agreement No. 101109532).

All simulations were performed on the national supercomputer 
HPE Apollo Systems \textit{HAWK} at the High Performance Computing Center Stuttgart (HLRS) under the grant number \textit{hpcmphas/44084}.
\\\\
\textbf{Conflict of interest} The corresponding author states on behalf of all authors, that there is no conflict of interest. 
\\\\
\textbf{Code availability} The open-source code FLEXI, on which all extensions are based, is available at www.flexi-project.org under the GNU GPL v3.0 license.
\\\\
\textbf{Availability of data and material} All data generated or analyzed during this study are included in this published article.

\bibliographystyle{spmpsci}     
\bibliography{references.bib}   

\end{document}